\documentclass[10pt]{article}

\usepackage[a4paper,margin=1in]{geometry}
\usepackage[T1]{fontenc}
\usepackage[utf8]{inputenc}
\usepackage{lmodern}
\usepackage{microtype}
\usepackage{cite}
\usepackage{mathrsfs,needspace,tikz}

\usepackage{amsmath,amssymb,amsthm,mathtools}
\usepackage[dvipsnames]{xcolor}

\usepackage[colorlinks=true,linkcolor=BrickRed,citecolor=RoyalBlue,urlcolor=RoyalBlue]{hyperref}
\usepackage[nameinlink,capitalize,noabbrev]{cleveref}

\numberwithin{equation}{section}

\newtheorem{theorem}{Theorem}[section]
\newtheorem{proposition}[theorem]{Proposition}
\newtheorem{lemma}[theorem]{Lemma}
\newtheorem{assumption}[theorem]{Assumption}
\theoremstyle{definition}
\newtheorem{remark}[theorem]{Remark}
\newtheorem{corollary}[theorem]{Corollary}

\crefname{assumption}{Assumption}{Assumptions}
\Crefname{assumption}{Assumption}{Assumptions}
\crefname{openproblem}{Open problem}{Open problems}
\Crefname{openproblem}{Open problem}{Open problems}

\newcommand{\ii}{\mathrm i}

\newcommand{\Env}{\boldsymbol{E}}

\newcommand{\spec}{\operatorname{spec}}

\hypersetup{pdftitle={Accidental Threefold Dirac Points at the Origin of the Brillouin Zone and Wave-Packet Dynamics},pdfauthor={Habib Ammari, Xin Fu and Bowen Li}}
\title{Accidental Threefold Dirac Points at the Origin of the Brillouin Zone and Wave-Packet Dynamics}
\author{Habib Ammari\thanks{Department of Mathematics, ETH Z\"urich, R\"amistrasse 101, 8092 Z\"urich, Switzerland, and Hong Kong Institute for Advanced Study, City University of Hong Kong, Kowloon Tong, Hong Kong SAR. Email: \texttt{habib.ammari@math.ethz.ch}.}
\and Xin Fu\thanks{School of Science, Institute for Theoretical Sciences, Westlake University, Hangzhou, Zhejiang 310030, P.R. China. Email: \texttt{1550862645cf@gmail.com}.}
\and Bowen Li\thanks{Department of Mathematics, City University of Hong Kong, Kowloon Tong, Hong Kong SAR. Email: \texttt{bowen.li@cityu.edu.hk}.}}
\date{}

\begin{document}
\maketitle
\begin{abstract}
We study accidental threefold degeneracies at zero quasimomentum for scalar periodic wave operators on the square lattice. At an isolated positive triple eigenvalue, the middle Bloch band has zero linear term. The two outer bands open linearly in every direction precisely when a Gram matrix of coupling coefficients is positive definite. Square symmetry makes the resulting cone isotropic. For weak perturbations of a constant medium, we construct two transverse tuning curves on which a singlet coincides with a doublet away from their common origin. The central band has positive quadratic curvature along each coordinate axis. For a dilute inclusion of size $\eta$, a variational gap condition gives a unique transverse crossing in a contrast interval of order $\eta^{-2}\log(1/\eta)$. The crossing frequency is asymptotic to $\sqrt{2\pi/\log(1/\eta)}$, and the cone speed vanishes linearly with $\eta$, with a coefficient determined by the far-field dipole moment. A three-component Dirac equation approximates positive-frequency wave packets with an error bound explicit in the material coefficients and the spectral separation at zero quasimomentum. When the singlet has nonzero cell mean, the response to a uniform cell source yields effective permittivity and permeability with simple zeros at the crossing frequency. We express their slopes in terms of the singlet mean and the cone speed. The resulting bulk impedance tends to zero in a weak-contrast example and grows inversely with the inclusion size in the dilute high-contrast crystal.
\end{abstract}

\noindent\textbf{Keywords.} Threefold Dirac point, Bloch bands, square lattice, wave-packet dynamics, high-contrast inclusions, periodic capacity, double-zero effective material properties.

\smallskip\noindent\textbf{2020 Mathematics Subject Classification.} Primary 35P05. Secondary 35B20, 35C20, 35L05, 47A55.

\setcounter{tocdepth}{1}
\hypersetup{bookmarksdepth=3}
\tableofcontents
\clearpage

\section{Introduction}\label{sec:introduction}

At the center of the Brillouin zone, three Bloch bands can meet in a cone, whereas the middle band has zero slope. In a medium with square symmetry, two periodic modes can share an eigenvalue by symmetry. We tune a material parameter to bring a third mode to the same eigenvalue. Under what conditions does the triple eigenvalue open into a cone? Can one determine the speed and the time scale of the resulting wave motion from the coefficients of the medium? When does the crossing give rise to simultaneous zeros of the effective permittivity and permeability?

We consider the scalar wave equation
\[
 b(x)\partial_t^2u-\nabla\cdot(A(x)\nabla u)=0,\qquad x\in\mathbb R^2,
\]
where the stiffness $A$ is real, symmetric, and uniformly elliptic, and the mass density $b$ is bounded and positive almost everywhere. Both coefficients are periodic with the unit square cell $Y=(-1/2,1/2)^2$. We represent its Brillouin zone by $Y^*=[-\pi,\pi)^2$, write $\alpha=(\alpha_1,\alpha_2)$ for the quasimomentum, and denote its origin by $\Gamma$. The eight rotations and reflections of the square form the dihedral group $D_4$. We call a one-dimensional symmetry type a \emph{singlet}. The standard two-dimensional representation is denoted by $E$, and a copy of it is called a \emph{doublet}. Two singlet types occur in our constructions. The fully invariant type is denoted by $A_1$. The type $B_1$ is even under reflection in either coordinate axis and changes sign under a quarter-turn. An eigenvalue in either singlet sector can coincide with a doublet eigenvalue at $\Gamma$.

Huang et al.\ observed an accidental threefold degeneracy in a square-lattice photonic crystal and related it to simultaneous zeros of the effective permittivity and permeability \cite{HuangLaiHangZhengChan2011}. Ashraf and Faryad examined how this interpretation depends on the symmetry of the modes at the crossing \cite{AshrafFaryad2016}. Sakoda calculated the isotropic cone formed by an $A_1$ singlet and an $E$ doublet in a tight-binding model \cite{Sakoda2012cone}, and used first-order perturbation of Maxwell's equations to identify the symmetry conditions for such crossings \cite{Sakoda2012symmetry}. Mei et al.\ computed the local dispersion and the zero Berry phase \cite{MeiWuChanZhang2012}. This configuration is also called a \emph{Dirac-like cone}. Its three modes admit a pseudospin-one description \cite{Fang2016}.

For honeycomb structures, Fefferman and Weinstein proved the existence of twofold Dirac points and justified the associated wave-packet dynamics \cite{FeffermanWeinstein2012,FeffermanWeinstein2014}. For the wave equation in a slowly modulated honeycomb medium, Xie and Zhu
justified a two-component Dirac approximation with a spatially varying
mass term \cite{xie}. Cassier and Weinstein studied the corresponding spectral problem for divergence-form operators with high contrast in the principal part \cite{CassierWeinstein2021}. Ammari, Fu and Jing obtained a two-component effective equation for subwavelength wave packets in a high-contrast honeycomb crystal \cite{ammari2026wave}. On the square lattice, Keller et al.\ studied quadratic touching at the Brillouin-zone corner \cite{Keller2018,Keller2020}. Chaban and Weinstein showed how a small deformation can split such a touching into two Dirac points \cite{ChabanWeinstein2025}. In the high-contrast limit, Kao, Lin and Osting used reflection symmetry to prove the existence of a twofold crossing of opposite-parity modes away from the high-symmetry points \cite{Kao2026}. Here, we construct crossings at $\Gamma$ by varying the coefficients of the wave equation. We determine the tuning parameters and the propagation speed from those coefficients.

Our first result concerns an isolated positive triple eigenvalue. 
For periodic media with Dirichlet or Neumann exclusions, Oudghiri-Idrissi,
Guzina, and Meng expressed the first-order splitting at a triple eigenvalue
in terms of three directional coupling coefficients \cite{guzina1}.
We express the linear part of the Bloch operator on its eigenspace through two real skew-symmetric coupling matrices. Their Gram matrix gives a basis-independent criterion for conical splitting in every direction. When the coefficients have $D_4$ symmetry, the cone is isotropic. Figure~\ref{fig:threefold-cone} illustrates the local frequency dispersion. We then approximate the positive-frequency evolution of a wave packet formed from the three periodic modes by a three-component first-order system. The error bound tracks both the relative separation from the remaining eigenvalues at $\Gamma$ and a bound on the coefficients. This dependence allows us to pass from a fixed medium to families in which the gap closes or the contrast diverges.

We give two constructions. The first positive periodic eigenvalue of the constant medium has multiplicity four. Its eigenspace decomposes into two singlets and a doublet. Two symmetry-preserving material perturbations vary the differences between their eigenvalues. A nondegeneracy condition on their Fourier moments yields two transverse analytic tuning curves. At the origin, the four modes coincide. Away from the origin on either curve, one singlet meets the doublet and the gap to the other singlet is linear to leading order in the distance to the origin. Restriction to an odd reflection sector shows that the central band has a positive quadratic term along each coordinate axis.

\begin{figure}[tbp]
\centering
\begin{tikzpicture}[x=1cm,y=1cm,>=stealth,font=\small,line cap=round,line join=round]
\definecolor{crystalink}{HTML}{174E60}
\definecolor{crystalfill}{HTML}{BADDE1}
\definecolor{crystalgrid}{HTML}{9CA9AB}
\begin{scope}
  \node[font=\small\bfseries] at (0,1.99) {(a) Square lattice};
  \fill[crystalfill!22] (-.52,-.52) rectangle (.52,.52);
  \foreach \u in {-1.56,-.52,.52,1.56}{
    \draw[crystalgrid,dashed,thin] (\u,-1.56)--(\u,1.56);
    \draw[crystalgrid,dashed,thin] (-1.56,\u)--(1.56,\u);
  }
  \foreach \x in {-1.04,0,1.04}
    \foreach \y in {-1.04,0,1.04}
      \draw[crystalink,fill=crystalfill] (\x,\y) circle[radius=.18];
  \draw[crystalink,thick] (-.52,-.52) rectangle (.52,.52);
  \node[crystalink] at (-.31,.31) {$Y$};
  \draw[<->] (-.52,-1.82)--node[fill=white,inner sep=1pt] {$1$} (.52,-1.82);
\end{scope}
\begin{scope}[xshift=4.45cm]
  \node[font=\small\bfseries] at (0,1.99) {(b) One resonator in $Y$};
  \draw[crystalink,thick,fill=crystalfill!15] (-1.56,-1.56) rectangle (1.56,1.56);
  \draw[crystalink,fill=crystalfill,thick] (0,0) circle[radius=.54];
  \node at (0,0) {$b=\tau$};
  \draw[crystalink,thin] (.38,.38)--(.73,.81)--(.98,.81);
  \node[anchor=south] at (.72,.82) {$D_\eta=\eta B$};
  \node at (-.90,-1.04) {$b=1$};
  \draw[thin] (-1.56,-1.56)--(-1.56,-1.64) (1.56,-1.56)--(1.56,-1.64);
  \node[below,font=\scriptsize] at (-1.56,-1.62) {$-1/2$};
  \node[below,font=\scriptsize] at (1.56,-1.62) {$1/2$};
\end{scope}
\begin{scope}[xshift=8.95cm]
  \node[font=\small\bfseries] at (0,1.99) {(c) Brillouin zone $Y^*$};
  \draw[crystalink,thick,fill=crystalfill!10] (-1.38,-1.38) rectangle (1.38,1.38);
  \draw[->,black!72] (-1.66,0)--(1.79,0) node[right,inner sep=1pt] {$\alpha_1$};
  \draw[->,black!72] (0,-1.61)--(0,1.67) node[right,inner sep=1.5pt] {$\alpha_2$};
  \fill[crystalink] (0,0) circle[radius=2pt];
  \node[below right,inner sep=3pt] at (0,0) {$\Gamma$};
  \node[above right,inner sep=2pt,font=\scriptsize] at (-1.38,0) {$-\pi$};
  \node[above left,inner sep=2pt,font=\scriptsize] at (1.38,0) {$\pi$};
\end{scope}
\end{tikzpicture}
\caption{Each square unit cell $Y$ contains one connected resonator $D_\eta=\eta B$. Here the reference shape $B$ is a disk. The array is invariant under quarter-turns and reflections about the cell center. Panel (a) highlights one cell and shows its lattice translates. Panel (b) enlarges the highlighted cell. Panel (c) shows the Brillouin zone of $\mathbb Z^2$ with $\Gamma$ at zero quasimomentum.}
\label{fig:crystal-geometry}
\end{figure}
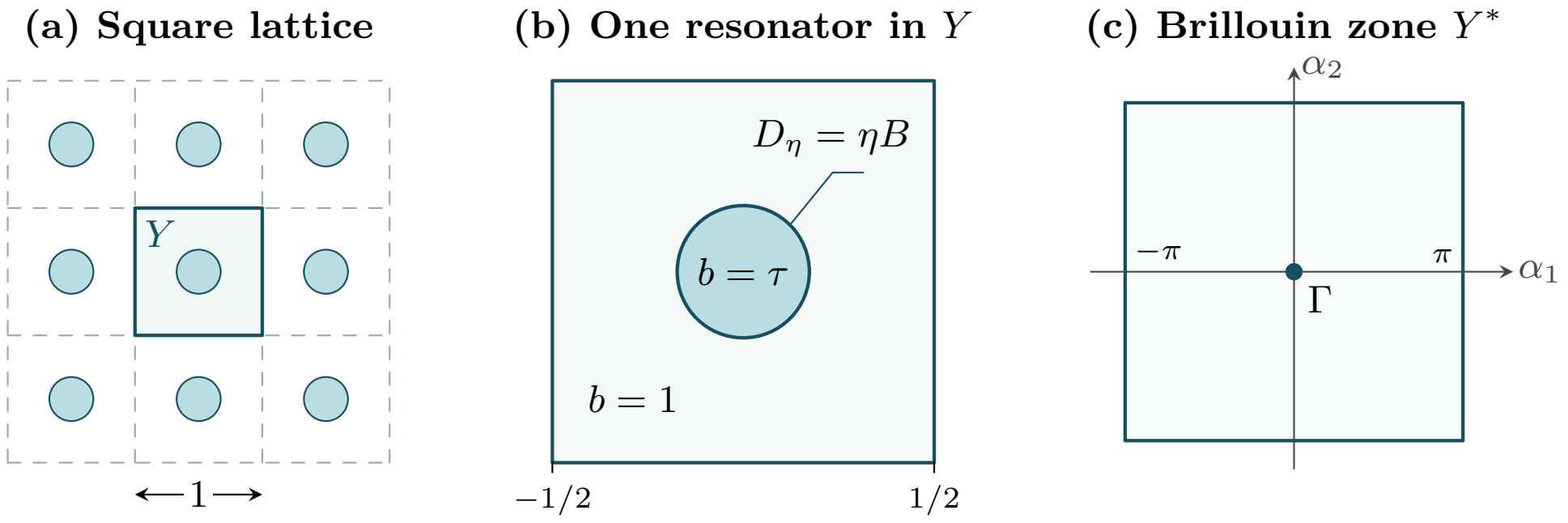

Our second construction uses one inclusion $D_\eta=\eta B$ in each unit cell, where $B$ is a fixed square-symmetric reference shape and $\eta>0$ is its scale. The stiffness is the identity matrix $A=I_2$, and the mass density is $\tau>1$ in each inclusion and one in the background. Figure~\ref{fig:crystal-geometry} shows the periodic array, the inclusion and the corresponding Brillouin zone. A centered disk or square inclusion gives a $D_4$-invariant material distribution. The relevant periodic modes are a capacitary singlet and two dipole modes related by a quarter-turn. These three modes belong to the same crystal with one inclusion per cell.

Potential operators describe localized modes of isolated high-index dielectric particles \cite{AmmariDabrowskiFitzpatrickMillienSini2019,AmmariLiZou2023}. In two dimensions, the logarithmic Green function distinguishes the principal mode from modes of zero mean \cite{AmmariLi2025}. For periodic arrays of thin dielectric rods, Holzmann and Lotoreichik showed that the material parameters can be chosen to open spectral gaps around prescribed frequencies for transverse magnetic waves \cite{holzmann}. In a dilute dielectric honeycomb lattice, Huang, Li, and Liu derived the slope of a symmetry-forced double eigenvalue at a Brillouin-zone vertex \cite{HuangLiLiu2026}. In our square cell, the dipole concentrates in the inclusion in the weighted $L^2$ norm. It competes at $\Gamma$ with a periodic singlet constrained by weighted orthogonality to the constant zero mode. The capacity of the small inclusion governs the singlet energy. Related logarithmic scales occur in perforated media and periodic small-trap problems \cite{allaire1991homogenization,Paquin2022}. Spectral representations for periodic media with contrast in the principal part were developed in \cite{LiptonViator2016}.

Put $\ell_\eta=\log(1/\eta)$. Under a gap condition on the reference shape, the singlet and doublet eigenvalues meet at a unique contrast in an interval of order $\ell_\eta/\eta^2$. The second, third and fourth bands then meet at $\Gamma$ at a frequency asymptotic to $\sqrt{2\pi/\ell_\eta}$. At $\Gamma$ the remaining eigenvalues stay at a uniformly positive relative distance from the triple. The cone speed is of order $\eta$, and its leading coefficient is determined by the far-field dipole moment of the localized mode. For a packet of quasimomentum bandwidth $\delta$, the approximation is valid over a propagation distance of order one in the slow spatial variable provided $\delta\sqrt{\ell_\eta}=o(\eta)$. We verify the gap condition for a square, its small symmetric deformations, and a disk. For a disk the speed coefficient can be evaluated explicitly.

In this dilute limit, the background wavelength becomes large compared with the lattice period, while the interior wavelength remains comparable to the inclusion size at the crossing. Thus, nonconstant internal modes contribute to the subwavelength spectrum. Remark~\ref{rem:bubble-comparison} compares this scaling with the bubble framework of \cite{AmmariDaviesHiltunen2026}, in which one connected high (acoustic) contrast resonator per cell contributes one subwavelength band. In our model the frequency of every fixed Bloch band tends uniformly to zero as $\eta\to0$, as shown in Remark~\ref{rem:dilute-band-count}. The three-component effective equation describes the isolated cluster at the crossing.

For the axial electric polarization, we derive the \emph{double-zero effective material properties} from a forced Bloch problem. Write $\omega_*$ for the crossing frequency and $v_D$ for the cone speed. The Dirac equation determines the product of the frequency slopes of the effective permittivity $\varepsilon_{\mathrm{eff}}$ and permeability $\mu_{\mathrm{eff}}$ at $\omega_*$. To determine the slopes separately, we take the macroscopic electric field to be the unweighted cell mean and apply a source with constant periodic amplitude. Its coupling to the real singlet $f_0$, normalized by $\int_Y b f_0^2\,dx=1$, is $m_*=\int_Y f_0\,dx$. When $m_*\ne0$, Proposition~\ref{prop:dz-bulk} gives simple zeros of both parameters and the limiting bulk impedance
\[
 Z_*:=\lim_{\omega\to\omega_*}
       \sqrt{\frac{\mu_{\mathrm{eff}}(\omega)}{\varepsilon_{\mathrm{eff}}(\omega)}}
     =\frac{m_*^2}{2v_D}.
\]
The two constructions give opposite impedance limits. Along the $A_1$ tuning curve in the piecewise-constant example of Section~\ref{subsec:weak-piecewise}, the singlet mean is of order $|p_1|$ and the cone speed tends to a positive constant as the perturbation parameter $p_1$ tends to zero. Hence $Z_*\asymp p_1^2$. In the dilute crystal, the positive-mean capacitary singlet satisfies $m_*\to1$, whereas $v_D$ is of order $\eta$. Corollary~\ref{cor:dilute-double-zero} gives the coefficient of the resulting $\eta^{-1}$ impedance law in terms of the reference inclusion.

Section~\ref{subsec:representation-preliminaries} introduces the symmetry projections. Proposition~\ref{prop:generalcriteria} in Section~\ref{sec:bloch} gives the conical splitting criterion. In Section~\ref{sec:effective}, Theorem~\ref{thm:effective-dynamics} bounds the wave-packet error and Section~\ref{subsec:double-zero} relates the Dirac equation to the forced cell problem. Sections~\ref{sec:weakperturbation} and~\ref{sec:dilute-verification} give the weak-contrast and dilute constructions in Theorems~\ref{thm:weak-dirac} and~\ref{thm:dilute-dirac}, respectively. Section~\ref{sec:conclusion} discusses the resulting propagation scales and open questions. Appendix~\ref{app:rep-proofs} collects the general representation-theoretic results.

\section{Symmetry spaces on the square lattice}\label{subsec:representation-preliminaries}

The symmetry decomposition separates the modes that enter the crossing. For coefficients with square symmetry, the operator at $\Gamma$ commutes with the rotations and reflections of the square. Modes of type $E$ (defined below) occur in doublets, and symmetry-preserving changes of the coefficients preserve their common eigenvalue. A singlet belongs to a different symmetry sector, so its eigenvalue can be tuned to meet the doublet without mixing the two sectors at $\Gamma$. The constructions below also establish spectral isolation and nonzero linear coupling, which are separate requirements for a threefold Dirac point.

In this section,
we describe the real symmetry types of the square and the projections onto the corresponding spaces. A representation is an action of the group by invertible linear maps. A nonzero real representation is irreducible when its only closed invariant subspaces are zero and the whole space. An irreducible representation is absolutely irreducible when every commuting real endomorphism is a scalar multiple of the identity. We call a one-dimensional real representation a character. For a finite group $G$, let $\widehat G$ contain one representative of each equivalence class of irreducible real representations. Appendix~\ref{app:rep-proofs} gives the general definitions and decomposition results used here.

We write $I_V$ for the identity on a vector space $V$, and $I_m$ for the identity matrix of order $m$. A superscript $\top$ denotes transpose and $*$ denotes the Hilbert-space adjoint. The groups $\mathrm O(H)$ and $\mathrm{SO}(m)$ denote, respectively, the surjective linear isometries of a real Hilbert space $H$ and the real orthogonal matrices of determinant one in dimension $m$.

\subsection{\texorpdfstring{The irreducible representations of $D_4$}{The irreducible representations of D4}}

We consider the dihedral group
\begin{equation}\label{rep-prelim:D4-presentation}
    D_4:=\langle r,s\mid r^4=e,\ s^2=e,\ srs=r^{-1}\rangle.
\end{equation}
Set the rotation matrix $R$ and the reflection matrix $S$ as follows:
\begin{equation}\label{rep-prelim:D4-standard-matrices}
    R:=\begin{pmatrix}0&-1\\1&0\end{pmatrix}, \qquad S:=\begin{pmatrix}1&0\\0&-1\end{pmatrix}.
\end{equation}
These matrices determine an orthogonal representation $E$ of $D_4$ on $\mathbb R^2$ by $E(r)=R$ and $E(s)=S$. We can identify $D_4$ with the matrix group generated by $R$ and $S$:
\begin{equation}\label{D4matrixgroup}
    D_4=\{R^kS^m:0\le k\le3,\ 0\leq m\leq 1\}.
\end{equation}

We classify the irreducible real representations of $D_4$.

\begin{lemma}\label{rep-prelim:D4-classification}
Up to equivalence, the irreducible real representations of $D_4$ are the four characters
\begin{equation}\label{rep-prelim:D4-characters}
    \chi_{\varepsilon,\delta}(R^kS^m):=\varepsilon^k\delta^m, \qquad \varepsilon,\delta\in\{1,-1\},
\end{equation}
and the standard two-dimensional representation $E$ given above. All five representatives are orthogonal and absolutely irreducible in the sense of \eqref{eq:scalar-commutant-hypothesis}. Thus, we may take
\begin{equation}\label{classficationofhatD4}
    \widehat{D_4}=\{\chi_{1,1},\chi_{1,-1},\chi_{-1,1},\chi_{-1,-1},E\}.
\end{equation}
\end{lemma}

\begin{proof}
Let $\tau$ be any irreducible real representation on a nonzero finite-dimensional space $V$, and write $A:=\tau(r)$ and $B:=\tau(s)$. The defining relations give
\begin{equation}
    A^4=I_V,\qquad B^2=I_V,\qquad BAB=A^{-1}.
\end{equation}
It follows that the operator $C:=A^2$ is an involution commuting with both $A$ and $B$, \textit{i.e.},
\begin{equation}
    C^2 = I_V, \qquad CA=AC, \qquad CB=BC.
\end{equation}
The eigenspaces of $C$ associated with the eigenvalues $1$ and $-1$ are therefore both invariant subspaces of the representation $\tau$, and
\begin{equation}
    V=\ker(C-I_V)\oplus\ker(C+I_V).
\end{equation}
Irreducibility implies $A^2=I_V$ or $A^2=-I_V$.

If $A^2=I_V$, then $AB=BA$. Choose a nonzero eigenspace of $A$, which is preserved by $B$. Since the restriction of $B$ to this space is an involution, it has a nonzero eigenvector $v$ there. Thus
\begin{equation}
    Av=\varepsilon v,\qquad Bv=\delta v, \qquad \varepsilon,\delta\in\{1,-1\}.
\end{equation}
The line $\mathbb Rv$ is invariant, so irreducibility gives $V=\mathbb Rv$ and $\tau\simeq\chi_{\varepsilon,\delta}$.

If $A^2=-I_V$, then $AB=-BA$. It follows that $\ker(B-I_V)$ is nonzero. Otherwise, $B=-I_V$. This contradicts $BA=-AB$ and invertibility of $A$. Choose $v\ne0$ with $Bv=v$, and put $w:=Av$. Then
\begin{equation}
    Av=w,\qquad Aw=-v,\qquad Bv=v,\qquad Bw=-w.
\end{equation}
The vectors $v,w$ are linearly independent because they belong to opposite eigenspaces of $B$. Their span is
invariant under both generators and therefore equals $V$. In the ordered basis $(v,w)$, the matrices of $A,B$ are
$R,S$, respectively. Hence $\tau\simeq E$.

Conversely, the four characters are irreducible and pairwise inequivalent, as seen from their values on $R,S$.
The representation $E$ is irreducible over $\mathbb R$: if a real line were invariant under $R$, some $v\ne0$
would satisfy $Rv=cv$ with $c\in\mathbb R$, giving $-v=R^2v=c^2v$, a contradiction. Its dimension distinguishes it from the four characters. This proves the classification \eqref{classficationofhatD4}.

All five displayed representatives are orthogonal. The four characters are absolutely irreducible because each acts on a one-dimensional real space.
For $E$, a matrix commuting with $S$ is diagonal, and a diagonal matrix commuting with $R$ has equal diagonal entries. Consequently,
\begin{equation}
    \{L\in\operatorname{End}(\mathbb R^2): LR=RL,\ LS=SL\}=\mathbb RI_2,
\end{equation}
which proves the absolute irreducibility of $E$.
\end{proof}

The four character labels are
\begin{equation}\label{rep-prelim:D4-type-labels}
    A_1:=\chi_{1,1},\qquad A_2:=\chi_{1,-1},\qquad B_1:=\chi_{-1,1},\qquad B_2:=\chi_{-1,-1}.
\end{equation}
Thus, $A_1$ is the trivial character. An individual one-dimensional irreducible copy is called a singlet, and an individual $E$-copy is called a doublet.

For every character $\chi$, we also have
\begin{equation}\label{rep-prelim:D4-character-twist}
    \chi\otimes E\simeq E.
\end{equation}
Indeed, multiplication of each $E(Q)$ by the nonzero scalar $\chi(Q)$ does not change the invariant subspaces.
Thus, $\chi\otimes E$ is irreducible of dimension two, and Lemma~\ref{rep-prelim:D4-classification} identifies its type with $E$.

\subsection{\texorpdfstring{Irreducible decompositions for $D_4$}{Irreducible decompositions for D4}}

Let $T:D_4\to\mathrm O(H)$ be an orthogonal representation on a separable real Hilbert space $H$, and write $T_Q=T(Q)$. In the orthogonal decomposition \eqref{decompofHeq}, $H_{\rho,a}$ denotes the $a$th copy of the symmetry type $\rho$, $A_\rho$ its index set, and $\Pi_{\rho,a}$ the projection onto that copy.
Every irreducible representation of $D_4$ is absolutely irreducible by Lemma~\ref{rep-prelim:D4-classification}, so Lemma~\ref{lem:matrix-coefficient-projections} applies to all five types. For a character $\chi$, and for $E$ with the standard orthonormal basis $e_1,e_2$ of $\mathbb R^2$, define
\begin{equation}\label{rep-prelim:D4-projection-averages}
\begin{aligned}
    P_\chi&:=\frac18\sum_{Q\in D_4}\chi(Q)T_Q,\\
    P_{E,i}&:=\frac{2}{8}\sum_{Q\in D_4}\langle Qe_i,e_i\rangle_{\mathbb R^2}T_Q,  \qquad i=1,2.
\end{aligned}
\end{equation}

\begin{lemma}\label{rep-prelim:D4-six-projections}
The averages \eqref{rep-prelim:D4-projection-averages} are
mutually orthogonal projections. The character projection
$P_\chi$ has range
\begin{equation}\label{rep-prelim:D4-character-range}
    \operatorname{Ran}P_\chi
       =\{u\in H:T_Qu=\chi(Q)u\ \mathrm{ for\ every\ }Q\in D_4\}.
\end{equation}
Set
\begin{equation}
    S_1:=R^2S=\begin{pmatrix}
        -1 & 0 \\
        0 & 1
    \end{pmatrix}, \qquad S_2:=S=\begin{pmatrix}
        1 & 0 \\
        0 & -1
    \end{pmatrix}.
\end{equation}
The two coordinate projections have the explicit forms
\begin{equation}\label{rep-prelim:D4-axial-projections}
\begin{aligned}
    P_{E,1}&=\frac14(I_H-T_{S_1})(I_H+T_{S_2}),\\
    P_{E,2}&=\frac14(I_H+T_{S_1})(I_H-T_{S_2}),
\end{aligned}
\end{equation}
and their ranges are
\begin{equation}\label{rep-prelim:D4-axial-ranges}
\begin{aligned}
    \operatorname{Ran}P_{E,1}&=\{u\in H:T_{S_1}u=-u,\ T_{S_2}u=u\},\\
    \operatorname{Ran}P_{E,2}&=\{u\in H:T_{S_1}u=u,\ T_{S_2}u=-u\}.
\end{aligned}
\end{equation}
Moreover,
\begin{equation}\label{rep-prelim:complete-six-decomposition}
\begin{aligned}
    P_E&:=P_{E,1}+P_{E,2}=\frac12(I_H-T_{R^2})=\sum_{a\in A_E}\Pi_{E,a},\\
    I_H&=\sum_{\varepsilon,\delta\in\{1,-1\}}P_{\chi_{\varepsilon,\delta}}+P_{E,1}+P_{E,2}.
\end{aligned}
\end{equation}
\end{lemma}

\begin{proof}
Lemma~\ref{lem:matrix-coefficient-projections} gives mutual
orthogonality and
\[
    P_\chi=\sum_{a\in A_\chi}\Pi_{\chi,a},
    \qquad
    P_{E,1}+P_{E,2}=\sum_{a\in A_E}\Pi_{E,a}.
\]
On each $\chi$-copy, every vector transforms by $\chi$,
so the closed span of the $\chi$-copies is contained in the right-hand
side of \eqref{rep-prelim:D4-character-range}.
Conversely, a vector $u$ in the right-hand side of \eqref{rep-prelim:D4-character-range} satisfies
\begin{equation}
    P_\chi u=\frac18\sum_{Q\in D_4}\chi(Q)^2u=u.
\end{equation}
This proves \eqref{rep-prelim:D4-character-range}.

For the representation $E$, a direct computation of $\langle Qe_i,e_i\rangle$ for $i=1,2$ and $Q\in D_4$ yields
\begin{equation}
    \begin{aligned}
    P_{E,1}&=\frac14(I_H-T_{R^2}-T_{S_1}+T_{S_2}),\\
    P_{E,2}&=\frac14(I_H-T_{R^2}+T_{S_1}-T_{S_2}).
\end{aligned}
\end{equation}
Since $S_1S_2=S_2S_1=R^2$, these expressions factor as \eqref{rep-prelim:D4-axial-projections}.

The operators $T_{S_1},T_{S_2}$ are commuting involutions.
The first expression in \eqref{rep-prelim:D4-axial-projections} consequently implies
\begin{equation}
    T_{S_1}P_{E,1}=-P_{E,1},\qquad T_{S_2}P_{E,1}=P_{E,1}.
\end{equation}
Conversely, if $T_{S_1}u=-u$ and $T_{S_2}u=u$, then $P_{E,1}u=u$. This identifies its range \eqref{rep-prelim:D4-axial-ranges}. The proof for $P_{E,2}$ follows the same argument.

Finally, \eqref{rep-prelim:complete-six-decomposition} follows from a direct computation and Lemma \ref{lem:matrix-coefficient-projections}.
\end{proof}

\section{\texorpdfstring{Threefold Dirac points at $\Gamma$ and square symmetry}{Threefold Dirac points at Gamma and square symmetry}}
\label{sec:bloch}

We first reduce the three Bloch bands at an isolated triple eigenvalue to a finite-dimensional spectral problem. We then use the square symmetry to identify the singlet and doublet and to determine the angular dependence of the cone.

\subsection{Periodic coefficients and the actions of square symmetry}

Let $Y=(-1/2,1/2)^2$ be the unit square cell. Let $A:\mathbb R^2\to\mathbb R^{2\times2}$ and $b:\mathbb R^2\to\mathbb R$ be measurable and $\mathbb Z^2$-periodic. Throughout this paper, we assume that 
\begin{itemize}
    \item[(a)] there are positive constants $a_*$ and $ a^* $ such that, for almost every $x\in\mathbb R^2$ and every $\xi\in\mathbb R^2$,
\begin{equation}\label{assumpforA}
    A(x) = A(x)^\top  \quad \mathrm{and} \quad a_*|\xi|^2\le \xi^\top A(x)\xi\le a^* |\xi|^2 ;
\end{equation}
    \item[(b)] there exists $p\in (1,\infty]$ such that
    \begin{equation}
        b|_Y \in L^p(Y)
    \end{equation}
    and $b(x)>0$ for almost every $x\in\mathbb R^2$.
\end{itemize}

We introduce the group action of $D_4$ on the fields over $Y$. Recall that $D_4$ can be realized as the matrix group \eqref{D4matrixgroup}. Let $V$ be a finite-dimensional linear space, and let $\rho$ be a representation of $D_4$ on $V$. For a field $F: Y\rightarrow V$, the induced action on $F$ is defined by
\begin{equation}\label{inducedaction}
    (Q\cdot_\rho F)(x):= \rho(Q) \big[ F(Q^\top x) \big], \qquad Q\in D_4,\, x\in Y.
\end{equation}
We say that $F$ is $\rho$-equivariant if $F$ is invariant under the action $\cdot_\rho$, namely, $Q\cdot_\rho F=F$ for every $Q\in D_4$. Equivalently, for every $Q\in D_4$ and almost every $x\in Y$,
\begin{equation}
    \rho(Q) F(x) =F(Qx).
\end{equation}

In particular, if $V=(\mathbb{C}^2)^{\otimes k}$ for $k\in \mathbb{N}$, and $\rho=\rho_k$ is the natural tensor power representation
\begin{equation}\label{tensorrepk}
    \rho_k(Q):= Q^{\otimes k}, \qquad Q\in D_4,
\end{equation}
then we say $F:Y\rightarrow V$ is $D_4$-\textit{invariant} if $F$ is $\rho_k$-equivariant. Throughout this paper, we only use the cases $k=0,1,2$, summarized as follows.
\begin{enumerate}
    \item If $k=0$, then $V=\mathbb{C}$ and $\rho_0$ is the trivial representation. Thus, a function $f:Y\rightarrow \mathbb{C}$ is $D_4$-invariant if
    \begin{equation}\label{symmetryk=0}
    f(x)= f(Q x) , \qquad Q\in D_4,\, x\in Y.
\end{equation}
    \item If $k=1$, then $V=\mathbb{C}^2$, and $\rho_1(Q) v = Qv$ for every $Q\in D_4$ and $v\in \mathbb{C}^2$. Thus, a vector field $F:Y\rightarrow \mathbb{C}^2$ is $D_4$-invariant if
    \begin{equation}\label{symmetryk=1}
    QF(x)= F(Qx) , \qquad Q\in D_4,\, x\in Y.
\end{equation}
    \item If $k=2$, then $V=\mathbb{C}^{2\times 2}$ and $\rho_2(Q) M = QMQ^\top$ for every $Q\in D_4$ and $M\in \mathbb{C}^{2\times 2}$. Thus, a matrix field $\mathbf{F}:Y\rightarrow \mathbb{C}^{2\times 2}$ is $D_4$-invariant if
    \begin{equation}\label{symmetryk=2}
    Q\mathbf{F}(x)Q^\top = \mathbf{F}(Qx) , \qquad Q\in D_4,\, x\in Y.
\end{equation}
\end{enumerate}

\subsection{\texorpdfstring{A general criterion for the existence of threefold Dirac points at $\Gamma$}{Criteria for threefold Dirac points at Gamma}}

Let
\[
Y^*:= \mathbb{R}^2/(2\pi \mathbb{Z}^2) \simeq [-\pi,\pi)^2
\]
be the Brillouin zone. The origin $(0,0) $ of the Brillouin zone $Y^*$ is denoted by $\Gamma$. For $\alpha\in Y^*$, consider the sesquilinear form
\begin{equation}\label{defaalphauv}
    a_\alpha(u,v):=\int_Y A(\nabla+\mathrm{i}\alpha)u\cdot\overline{(\nabla+\mathrm{i}\alpha)v}\,dx,
    \qquad u,v\in H^1(Y).
\end{equation}
Define the $b$-weighted $L^2$ space on $Y$ by
\begin{equation}
    L^2_b(Y):=L^2(Y,b\,dx), \qquad \langle u,v\rangle_{L^2_b(Y)}:=\int_Y bu\overline v\,dx.
\end{equation}

\begin{lemma}
    Let
    \[
    m_b:=\int_Y\min\{b(x),1\}\,dx>0.
    \]
    The following hold.
    \begin{itemize}
        \item[\emph{(a)}] For every $u\in H^1(Y)$,
        \begin{equation}\label{eq:cell-control}
            \|u\|_{L^2(Y)}^2\leq \frac{1}{\pi^2} \left( 1+\frac{2}{m_b}\right) \|\nabla u\|_{L^2(Y)}^2
            + \frac{2}{m_b} \|u\|_{L_b^2(Y)}^2.
        \end{equation}
        \item[\emph{(b)}] The embedding $H^1(Y)\hookrightarrow L_b^2(Y)$ is compact. The periodic subspace $H^1_\#(Y)$ of $H^1(Y)$ is dense in $L_b^2(Y)$.
        \item[\emph{(c)}] There are constants $c=c(a_*,m_b)>0$ and $C=C(a^*,p,\|b\|_{L^p(Y)})>0$ such that, for every $\alpha \in Y^*$ and $u\in H^1(Y)$, 
\[
    c\|u\|_{H^1(Y)}^2\le a_\alpha(u,u)+\|u\|_{L_b^2(Y)}^2\le C\|u\|_{H^1(Y)}^2.
\]
    \end{itemize}
\end{lemma}
\begin{proof}
    Set $\widetilde b:=\min\{b,1\}$ and $u_Y:=\int_Y u\,dx$. Since $0<\widetilde b\le1$, we have
\[
\begin{aligned}
    m_b|u_Y|^2&\le2\int_Y\widetilde b|u|^2\,dx+2\int_Y\widetilde b|u-u_Y|^2\,dx\\
    &\le2\|u\|_{L_b^2(Y)}^2+2\|u-u_Y\|_{L^2(Y)}^2.
\end{aligned}
\]
Combining this with
\[
    \|u\|_{L^2(Y)}^2=\|u-u_Y\|_{L^2(Y)}^2+|u_Y|^2
\]
and the Poincar\'e--Wirtinger inequality on $Y$,
\[
    \|u-u_Y\|_{L^2(Y)}^2\le\frac{1}{\pi^2}\|\nabla u\|_{L^2(Y)}^2,
\]
proves~(a).

Let
\[
    q:=\begin{cases}
        \dfrac{2p}{p-1}, & 1<p<\infty,\\[1ex]
        2, & p=\infty.
    \end{cases}
\]
By H\"older's inequality,
\[
    \|g\|_{L_b^2(Y)}^2\le\|b\|_{L^p(Y)}\|g\|_{L^q(Y)}^2, \qquad g\in L^q(Y).
\]
This estimate, combined with the compact embedding $H^1(Y)\hookrightarrow L^q(Y)$, yields the compact embedding $H^1(Y)\hookrightarrow L_b^2(Y)$.

Since $b\,dx$ is a finite Radon measure on the torus, smooth periodic functions are dense in $L_b^2(Y)$.
Hence,  $H^1_\#(Y)$ is dense in $L_b^2(Y)$. This proves (b).

For (c), let $\alpha\in Y^*$ and $u\in H^1(Y)$. The function $w(x):=\mathrm e^{\mathrm i\alpha\cdot x}u(x)$ belongs to $H^1(Y)$ and satisfies
\[
    |w|=|u|, \qquad \nabla w =\mathrm e^{\mathrm i\alpha\cdot x}(\nabla+\mathrm i\alpha)u.
\]
Applying~(a) to $w$ and using the ellipticity $a_*I \leq A\leq a^* I$, we obtain
\[
    \|u\|_{L^2(Y)}^2\le\frac{1}{\pi^2a_*}\left(1+\frac{2}{m_b}\right)
    a_\alpha(u,u)+\frac{2}{m_b}\|u\|_{L_b^2(Y)}^2.
\]
Moreover, $|\alpha|^2\le2\pi^2$, so
\[
\begin{aligned}
    \|\nabla u\|_{L^2(Y)}^2&\le2\|(\nabla+\mathrm i\alpha)u\|_{L^2(Y)}^2
    +2|\alpha|^2\|u\|_{L^2(Y)}^2\\
    &\le\frac{2}{a_*}a_\alpha(u,u)+4\pi^2\|u\|_{L^2(Y)}^2.
\end{aligned}
\]
These two inequalities yield
\[
    c\|u\|_{H^1(Y)}^2\le a_\alpha(u,u)+\|u\|_{L_b^2(Y)}^2
\]
with $c=c(a_*,m_b)>0$.

For the upper bound, H\"older's inequality and the Sobolev embedding give
\[
    \|u\|_{L_b^2(Y)}^2\le\|b\|_{L^p(Y)}\|u\|_{L^q(Y)}^2
    \le C_p\|b\|_{L^p(Y)}\|u\|_{H^1(Y)}^2,
\]
where $C_p$ depends only on $p$. Also,
\begin{equation}
    a_\alpha(u,u)\le
    2a^*\bigl(\|\nabla u\|_{L^2(Y)}^2+|\alpha|^2\|u\|_{L^2(Y)}^2\bigr)
    \le4\pi^2a^*\|u\|_{H^1(Y)}^2.
\end{equation}
Adding the last two estimates proves the upper bound with $C=C(a^*,p,\|b\|_{L^p(Y)})>0$.
\end{proof}

By this lemma, for every $\alpha\in Y^*$, restricting $a_\alpha$ to $H^1_\#(Y) \times H^1_\#(Y)$ gives a densely defined, closed, nonnegative symmetric sesquilinear form on $L_b^2(Y)$.
By the representation theorems for closed forms
\cite[Theorems~2.1 and~2.6]{Kato1995}, there exists a unique nonnegative self-adjoint operator $L(\alpha)$ on $L_b^2(Y)$ such that
\begin{equation}\label{blocheigenprobL}
    a_\alpha(u,v)=\langle L(\alpha)u,v\rangle_{L_b^2(Y)},
    \qquad u\in\operatorname{Dom}L(\alpha) \subset H^1_\#(Y),\quad v\in H^1_\#(Y).
\end{equation}
Formally, 
\begin{equation}
    L(\alpha)= - \frac{1}{b} (\nabla+\mathrm{i}\alpha) \cdot A(\nabla+\mathrm{i}\alpha)  .
\end{equation}
The compact embedding $H^1_\#(Y)\hookrightarrow L_b^2(Y)$ implies that $L(\alpha)$ has compact resolvent.
The Bloch eigenvalue problem for $L(\alpha)$ is precisely
\begin{equation}\label{prob:blocheigenpair}
    -(\nabla+\mathrm{i}\alpha) \cdot \{A(\nabla+\mathrm{i}\alpha)  u \} = \lambda b u, \qquad u\in H^1_\#(Y).
\end{equation}

We denote the Bloch eigenvalues, listed in nondecreasing order
and counted with multiplicity, by
\begin{equation}
    0\leq\lambda_1(\alpha)\leq\lambda_2(\alpha)
    \leq\cdots \leq \lambda_n(\alpha) \leq \cdots .
\end{equation}
Each ordered eigenvalue $\lambda_n(\alpha)$ depends continuously on $\alpha$; see \cite{kuchment2016overview}. 

To study threefold Dirac points at $\Gamma$, we make the following assumption.

\begin{assumption}\label{assump:threefold}
    At $\alpha=0$, \eqref{prob:blocheigenpair} has an eigenvalue
    $\lambda_*>0$ of multiplicity three.
\end{assumption}

Under Assumption \ref{assump:threefold}, let $n_*$ be the smallest index such that $\lambda_{n_*}(0)=\lambda_*$.
Then
\begin{equation}
    \lambda_{n_*}(0)=\lambda_{n_*+1}(0)=\lambda_{n_*+2}(0)=\lambda_*.
\end{equation}
Since the eigenvalues at $\alpha=0$ are discrete, by continuity of the Bloch eigenvalues, there exists a neighborhood $U$ of the origin such that, for every $\alpha \in U$, the three eigenvalues $\lambda_{n_*}(\alpha)$, $\lambda_{n_*+1}(\alpha)$, and $\lambda_{n_*+2}(\alpha)$ are isolated from the remaining Bloch eigenvalues. For such $\alpha$, define
\begin{equation}\label{eq:threefold-branches}
    \lambda_-(\alpha):=\lambda_{n_*}(\alpha), \qquad \lambda_c(\alpha):=\lambda_{n_*+1}(\alpha),
    \qquad \lambda_+(\alpha):=\lambda_{n_*+2}(\alpha).
\end{equation}
Thus,
\begin{equation}
    \lambda_-(\alpha)\leq\lambda_c(\alpha) \leq\lambda_+(\alpha), \qquad
    \lambda_-(0)=\lambda_c(0)=\lambda_+(0)=\lambda_*.
\end{equation}

Write $\delta_{rs}$ for the Kronecker symbol. Choose a real basis $\{f_0,f_1,f_2\}\subset H^1_\#(Y;\mathbb{R})$ of the eigenspace corresponding to $\lambda_*$ at $\alpha=0$, normalized by
\begin{equation}\label{b-orthogonality}
        \int_Y b f_r f_s \,dx = \delta_{rs}, \qquad r,s\in \{0,1,2\}.
\end{equation}
Let $e_1,e_2$ be the coordinate vectors in $\mathbb R^2$. Define the $3\times 3$ coupling matrix $C_j$, $j=1,2$, by
    \begin{equation}\label{defofCj}
        (C_j)_{rs}:= \int_Y (f_r A \nabla f_s - f_s A \nabla f_r) \cdot e_j \,dx,\qquad r,s\in \{0,1,2\},
    \end{equation}
    and the $2\times 2$ Gram matrix $G$ by
    \begin{equation}
        G_{jk}:= \frac{1}{2}\operatorname{tr}(C_j^\top C_k) =\sum_{0\leq r<s\leq 2} (C_j)_{rs} (C_k)_{rs}.
    \end{equation}

Under the change of basis
\begin{equation}\label{fOtop}
    \begin{pmatrix}
        f_0 \\ f_1 \\ f_2
    \end{pmatrix}\mapsto
    O^\top\begin{pmatrix}
        f_0 \\ f_1 \\ f_2
    \end{pmatrix}
\end{equation}
for $O\in \mathrm{O}(3)$, the matrices $C_j$ and $G$ transform as
\begin{equation}\label{CjGtransform}
    C_j \mapsto O^\top C_j O, \qquad G\mapsto G.
\end{equation}
Thus, the matrix $G$ is independent of the choice of real basis that satisfies \eqref{b-orthogonality}. 

Write $\mathcal L(U,V)$ for the space of bounded linear maps between normed spaces $U,V$. Let $F:\mathbb C^3\to L^2_b(Y)$ be the map
\begin{equation}\label{defFC3toL2b}
    Fz:=\sum_{r=0}^2z_rf_r.
\end{equation}
The normalization \eqref{b-orthogonality} of $\{f_0,f_1,f_2\}$ gives $F^*F=I_3$. Set $P(0)=FF^*$, the orthogonal projection in $L^2_b(Y)$ onto the complexified eigenspace. Then $P(0)F=F$. 

We call the crossing a \emph{threefold Dirac point} if the outer branches open linearly in every direction and the middle branch has vanishing linear dispersion.
Building on the three-mode dispersion formula in \cite{guzina1}, we give a basis-independent criterion for conical splitting in every direction and determine the
leading Bloch eigenfunctions without assuming square symmetry.

\begin{proposition}\label{prop:generalcriteria}
    Under Assumption~\ref{assump:threefold}, the crossing at $\Gamma$ is a threefold Dirac point if and only if $G$ is positive definite. If $G>0$, then near $\alpha=0$,
        \begin{equation}\label{diracconegeneralcriteria}
        \begin{aligned}
            \lambda_{\pm}(\alpha)&=\lambda_* \pm \sqrt{\alpha^\top G\alpha} + O(|\alpha|^2), \\
            \lambda_c(\alpha)&=\lambda_*+ O(|\alpha|^2).
        \end{aligned}
        \end{equation}
    Under the same positivity condition, for any real basis normalized by \eqref{b-orthogonality}, there exists a real orthogonal matrix $O\in \mathrm{SO}(3)$, independent of $\alpha$, such that
    \begin{equation}\label{normalform}
        O^\top(\alpha_1C_1+\alpha_2C_2)O = \begin{pmatrix}
            0&\theta^\top\\
            -\theta&0_{2\times 2}
        \end{pmatrix}, \qquad \theta=\theta(\alpha):=G^{1/2}\alpha.
    \end{equation}
    Define the real basis
    \begin{equation}\label{ftildetransform}
        \begin{pmatrix}
        \widetilde{f}_0 \\ \widetilde{f}_1 \\ \widetilde{f}_2
    \end{pmatrix} := O^\top\begin{pmatrix}
        f_0 \\ f_1 \\ f_2
    \end{pmatrix}.
    \end{equation}
    We use $O_{\mathcal X}(h)$ for a remainder of norm at most $C|h|$ in a normed space $\mathcal X$, with the medium fixed. At $\alpha \neq 0$ sufficiently small, the eigenfunctions corresponding to $\lambda_\pm(\alpha)$ and $\lambda_c(\alpha)$ can be normalized by
    \begin{equation}
        \int_Y b|u_s(\alpha)|^2\,dx=1, \qquad s\in\{+,c,-\},
    \end{equation}
    and their phases can be chosen so that
    \begin{equation}\label{eigenfunctionexpansion}
        \begin{aligned}
            u_\pm(\alpha) &= \frac{1}{\sqrt{2}} \widetilde f_0 \pm \frac{\mathrm{i}}{\sqrt{2}} \left( \frac{\theta_1}{|\theta|} \widetilde f_1 + \frac{\theta_2}{|\theta|} \widetilde f_2 \right)  +\mathcal{O}_{H^1(Y)}(|\alpha|),\\
            u_c(\alpha)&= - \frac{\theta_2}{|\theta|}\widetilde f_1 + \frac{\theta_1}{|\theta|} \widetilde f_2
            + \mathcal{O}_{H^1(Y)}(|\alpha|).
        \end{aligned}
    \end{equation}
\end{proposition}

\begin{proof}

\textit{Step 1.} We first reduce the isolated three-dimensional spectral subspace to a Hermitian matrix. Expanding the form $a_\alpha$ gives
\begin{equation}
    \begin{aligned}
    a_\alpha(u,v)&=a_0(u,v)+\sum_{j=1}^2\alpha_j a_j^{(1)}(u,v)+\int_Y (\alpha^\top A\alpha)u\overline v\,dx,\\
    a_j^{(1)}(u,v)&:=\mathrm{i}\int_Y\bigl(uA\nabla\overline v-\overline v A\nabla u\bigr)\cdot e_j\,dx.
\end{aligned}
\end{equation}
In particular, the expanded expression in real $\alpha \in Y^*$ admits a holomorphic extension to $\mathbb{C}^2$.

Choose a positively oriented circle $\mathcal C$ centered at $\lambda_*$ that encloses no other eigenvalue of $L(0)$ and lies in its resolvent set. For sufficiently small $\alpha$, define
\begin{equation}
    P(\alpha):=\frac{1}{2\pi\mathrm{i}} \int_{\mathcal C}(z-L(\alpha))^{-1}\,dz.
\end{equation}
To justify analyticity in the form norm, set $V=H^1_\#(Y)$ and let $V'$ be its antidual. The pencil $z\langle\cdot,\cdot\rangle_b-a_0$ maps $V$ isomorphically onto $V'$ on $\mathcal C$, as follows by expansion in an eigenbasis of $L(0)$. The inverse is uniformly bounded there. The polynomial form perturbation and a Neumann series give a holomorphic inverse $V'\to V$ for complex $\alpha$ near zero, uniformly on $\mathcal C$ \cite[Ch.~VII, \S4]{Kato1995}. Composing with the embedding $L_b^2(Y)\hookrightarrow V'$ and integrating shows that $\alpha\mapsto P(\alpha)$ is holomorphic with values in $\mathcal L(L_b^2(Y),V)$. For real $\alpha$, the projection is orthogonal and its rank is three by continuity. Its range is the spectral subspace of the three indicated bands.

For sufficiently small real $\alpha$, the matrix $S(\alpha):=F^*P(\alpha)F$ is positive definite. Define
\begin{equation}\label{eq:generalcriteria:frame-definition}
    \mathcal W(\alpha):=P(\alpha)F S(\alpha)^{-1/2}.
\end{equation}
Then,
\begin{equation}\label{eq:generalcriteria:frame-proof}
\begin{aligned}
    \mathcal W(\alpha)^*\mathcal W(\alpha)&=I,\\
    \operatorname{range}\mathcal W(\alpha)&=\operatorname{range}P(\alpha),\\
    \mathcal W(\alpha)&=F+\mathcal O_{\mathcal L(\mathbb C^3;H^1_\#(Y))}(|\alpha|).
\end{aligned}
\end{equation}
The frame $\mathcal W$ is real analytic in $\alpha$. The first two identities in \eqref{eq:generalcriteria:frame-proof} imply that 
\begin{equation}\label{WWstarP}
    \mathcal W(\alpha) \mathcal W(\alpha)^* =P(\alpha).
\end{equation}

Define the Hermitian matrix $M(\alpha)$ by
\begin{equation}
    M_{rs}(\alpha):=a_\alpha\bigl(\mathcal W(\alpha)\mathbf e_s,\mathcal W(\alpha)\mathbf e_r\bigr),
    \qquad r,s\in\{0,1,2\},
\end{equation}
where $(\mathbf e_0,\mathbf e_1,\mathbf e_2)$ denotes the standard basis of $\mathbb C^3$. Equivalently,
$M(\alpha)=\mathcal W(\alpha)^*L(\alpha)\mathcal W(\alpha)$. Moreover, since the range of $\mathcal{W}(\alpha)$ is in the range of $P(\alpha)$, and $\mathrm{range}\,P(\alpha)$ is an invariant subspace of $L(\alpha)$, we get
\begin{equation}\label{LWequalWM}
    \begin{aligned}
        L(\alpha)\mathcal{W}(\alpha) & =P(\alpha) L(\alpha)\mathcal{W}(\alpha) \\
        &=\mathcal W(\alpha) \mathcal W(\alpha)^* L(\alpha)\mathcal{W}(\alpha) \\
        &=\mathcal W(\alpha) M(\alpha),
    \end{aligned}
\end{equation}
where we used \eqref{WWstarP} in the second line. Thus, $M(\alpha)$ represents the restriction of $L(\alpha)$ to the spectral subspace associated with $\lambda_-(\alpha),\lambda_c(\alpha),\lambda_+(\alpha)$, and the eigenvalues of $M(\alpha)$ are exactly $\lambda_-(\alpha),\lambda_c(\alpha),\lambda_+(\alpha)$.

\textit{Step 2.} We claim that, with $C(\alpha):=\alpha_1C_1+\alpha_2C_2$,
\begin{equation}\label{eq:generalcriteria:reduction-proof}
    M(\alpha)=\lambda_*I-\mathrm{i}C(\alpha)+R(\alpha), \qquad |R(\alpha)|\leq K|\alpha|^2.
\end{equation}
Clearly, $M(0)=\lambda_*I$. To calculate its first derivatives, set
\begin{equation}
    w_{r,j}:=\left.\partial_{\alpha_j}\mathcal W(\alpha)\mathbf e_r \right|_{\alpha=0}, 
    \qquad r\in \{0,1,2\},\,j\in \{1,2\}.
\end{equation}
Differentiating $\mathcal W(\alpha)^*\mathcal W(\alpha)=I$ gives
\begin{equation}
    \langle w_{s,j},f_r\rangle_{L^2_b(Y)}+\langle f_s,w_{r,j}\rangle_{L^2_b(Y)}=0.
\end{equation}
Since each $f_r$ is an eigenfunction of $L(0)$ for the eigenvalue $\lambda_*$,
\begin{equation}
    a_0(w_{s,j},f_r)+a_0(f_s,w_{r,j})=\lambda_*\bigl(\langle w_{s,j},f_r\rangle_{L^2_b(Y)}+
    \langle f_s,w_{r,j}\rangle_{L^2_b(Y)}\bigr)=0.
\end{equation}
Thus,
\begin{equation}
    \begin{aligned}
    \left.\partial_{\alpha_j}M_{rs}(\alpha)\right|_{\alpha=0}&=a_j^{(1)}(f_s,f_r)\\
    &=\mathrm{i}\int_Y(f_sA\nabla f_r-f_rA\nabla f_s)\cdot e_j\,dx\\
    &=-\mathrm{i}(C_j)_{rs}.
\end{aligned}
\end{equation}
The estimate \eqref{eq:generalcriteria:reduction-proof} now follows from the real analyticity of $M$.

Next, we construct the normal form \eqref{normalform}. Each $C_j$ is real and skew-symmetric. Define
\begin{equation}
    \mathcal{C}_j:=\begin{pmatrix}
        -(C_j)_{12}\\ (C_j)_{02}\\ -(C_j)_{01}
    \end{pmatrix} \in \mathbb{R}^3, \qquad j=1,2, \qquad \mathrm{and} \qquad \mathcal{C}(\alpha) = \alpha_1 \mathcal{C}_1 + \alpha_2 \mathcal{C}_2.
\end{equation}
Then, $C_jz=\mathcal{C}_j\times z$ for every $z\in\mathbb C^3$, and
\begin{equation}
    G_{jk}=\mathcal{C}_j\cdot \mathcal{C}_k.
\end{equation}
Since $G>0$, the vectors $\mathcal C_1,\mathcal C_2$ are linearly independent and
$|\mathcal{C}_1\times \mathcal{C}_2|^2=\det G>0$. Set
\begin{equation}
    p:=\frac{\mathcal{C}_1\times \mathcal{C}_2}{|\mathcal{C}_1\times \mathcal{C}_2|},\qquad
    B:=\begin{pmatrix} p \times \mathcal{C}_1&p\times \mathcal{C}_2\end{pmatrix} .
\end{equation}
Then, $|p|=1$, $p^\top B=0$, and $G=B^\top B$. Define
\begin{equation}
    O:=\begin{pmatrix}
        p&B
    \end{pmatrix} \begin{pmatrix}
        1&0 \\
        0&G^{-1/2}
    \end{pmatrix}= \begin{pmatrix}p&BG^{-1/2}\end{pmatrix}.
\end{equation}
It follows that $O$ is an orthogonal matrix. Moreover,
\begin{equation}
    \det\begin{pmatrix}p&B\end{pmatrix}=\sqrt{\det G},
\end{equation}
so $\det O=1$. In particular, $O\in\mathrm{SO}(3)$ and is independent of $\alpha$. With $\theta=G^{1/2}\alpha$, we have
\begin{equation}\label{OCO1}
    C(\alpha)p = \mathcal{C}(\alpha) \times p=-BG^{-1/2}\theta.
\end{equation}
Moreover, using $C_j(-C_kp)=G_{jk}p$, we get
\begin{equation}\label{OCO2}
    C(\alpha)BG^{-1/2}=p\theta^\top.
\end{equation}
Combining \eqref{OCO1}--\eqref{OCO2} gives
\begin{equation}\label{eq:generalcriteria:normalform-proof}
    O^\top C(\alpha)O=\begin{pmatrix}
        0&\theta^\top\\
        -\theta&0_{2\times2}
      \end{pmatrix}.
\end{equation}

We now prove the eigenvalue expansion \eqref{diracconegeneralcriteria}. By \eqref{eq:generalcriteria:normalform-proof}, the Hermitian matrix $-\mathrm{i}C(\alpha)$ has three eigenvalues
\begin{equation}
    \pm|\theta|= \pm \sqrt{\alpha^\top G\alpha}, \quad 0.
\end{equation}
By \eqref{eq:generalcriteria:reduction-proof} and the min--max principle for Hermitian matrices, it follows that
\begin{equation}
    |\lambda_c(\alpha)-\lambda_*| + \left|\lambda_\pm(\alpha)-\left(\lambda_*\pm \sqrt{\alpha^\top G\alpha} \right) \right| \leq 2K|\alpha|^2.
\end{equation}

\textit{Step 3.} It remains to prove the eigenfunction expansion \eqref{eigenfunctionexpansion}. For $\alpha\neq0$, define
\begin{equation}\label{defofzs}
    z_\pm(\alpha):=\frac{1}{\sqrt2}
    \begin{pmatrix}
       1 \\ \pm\mathrm{i}\theta_1/|\theta| \\ \pm\mathrm{i}\theta_2/|\theta|
    \end{pmatrix},
    \qquad z_c(\alpha):=\begin{pmatrix}0\\-\theta_2/|\theta|\\\theta_1/|\theta|\end{pmatrix}.
\end{equation}
By \eqref{eq:generalcriteria:normalform-proof}, they are orthonormal eigenvectors
of $-\mathrm{i}O^\top C(\alpha)O$ with eigenvalues $\pm |\theta|$ and $0$, respectively. Let
\begin{equation}
    \widetilde M(\alpha):=O^\top M(\alpha)O, \qquad \widetilde R(\alpha):=O^\top R(\alpha)O,
    \qquad \mu_s(\alpha):= \left\{ 
    \begin{aligned}
        & \lambda_* \pm |\theta| ,&& s=\pm, \\
        & \lambda_*, && s=c.
    \end{aligned}\right.
\end{equation}
For each $s\in\{+,c,-\}$, let $\Pi_s(\alpha)$ be the orthogonal projection onto the one-dimensional eigenspace of
$\widetilde M(\alpha)$ corresponding to $\lambda_s(\alpha)$. The residual equation is
\begin{equation}
    \bigl(\widetilde M(\alpha)-\mu_s(\alpha)I\bigr) z_s(\alpha) =\widetilde R(\alpha)z_s(\alpha).
\end{equation}
The separation of the three linear eigenvalues is at least $\sqrt{\lambda_{\min}(G)}|\alpha|$. By the quadratic remainder, the other two eigenvalues remain at distance at least $\tfrac12\sqrt{\lambda_{\min}(G)}|\alpha|$ from $\mu_s(\alpha)$ for small $|\alpha|$. Restricting the residual equation to the range of $I-\Pi_s(\alpha)$ therefore yields
\begin{equation}
    |(I-\Pi_s(\alpha))z_s(\alpha)| \leq\frac{C}{|\alpha|}|\widetilde R(\alpha)z_s(\alpha)| \leq C|\alpha|.
\end{equation}
In particular, $\Pi_s(\alpha)z_s(\alpha)\neq0$ for sufficiently small $|\alpha|>0$. Choose the phase by setting
\begin{equation}
    \widetilde{z}_s(\alpha):= \frac{\Pi_s(\alpha)z_s(\alpha)} {|\Pi_s(\alpha)z_s(\alpha)|}.
\end{equation}
Then $\widetilde{z}_s(\alpha)$ is a unit eigenvector of $\widetilde M(\alpha)$ and
\begin{equation}
    \widetilde{z}_s(\alpha)=z_s(\alpha)+O(|\alpha|).
\end{equation}

Finally, define the Bloch eigenfunctions by
\begin{equation}
    u_s(\alpha):=\mathcal W(\alpha)O \widetilde{z}_s(\alpha).
\end{equation}
The property \eqref{eq:generalcriteria:frame-proof} of $\mathcal{W}(\alpha)$ and the identity \eqref{LWequalWM} show that
\begin{equation}
    L(\alpha)u_s(\alpha)=\lambda_s(\alpha)u_s(\alpha),\qquad \|u_s(\alpha)\|_{L^2_b(Y)}=1,
\end{equation}
and
\begin{equation}
    \begin{aligned}
    \|u_s(\alpha)-FOz_s(\alpha)\|_{H^1(Y)}
    &\leq \|(\mathcal W(\alpha)-F)O \widetilde{z}_s(\alpha)\|_{H^1(Y)}+\|FO(\widetilde{z}_s(\alpha)-z_s(\alpha))\|_{H^1(Y)}\\
    &\leq C|\alpha|.
\end{aligned}
\end{equation}
By definition, $FO\mathbf e_r=\widetilde f_r$. Substituting the expression \eqref{defofzs} of $z_s(\alpha)$ into $FOz_s(\alpha)$ gives
\begin{equation}
    \begin{aligned}
    u_\pm(\alpha)&=\frac{\widetilde f_0}{\sqrt2}\pm\frac{\mathrm{i}}{\sqrt2}
    \left(\frac{\theta_1}{|\theta|}\widetilde f_1+\frac{\theta_2}{|\theta|}\widetilde f_2\right)
    +\mathcal O_{H^1(Y)}(|\alpha|),\\
    u_c(\alpha)&=-\frac{\theta_2}{|\theta|}\widetilde f_1+\frac{\theta_1}{|\theta|}\widetilde f_2
    +\mathcal O_{H^1(Y)}(|\alpha|).
\end{aligned}
\end{equation}
This proves the assertions when $G>0$. Conversely, if $G$ is singular, choose a unit vector $e$ in its kernel. The definition of $G$ gives $\tfrac12\operatorname{tr}(C(e)^\top C(e))=e^\top Ge=0$, so $C(e)=0$. By \eqref{eq:generalcriteria:reduction-proof}, all three eigenvalues along $\alpha=te$ are equal to $\lambda_*+O(t^2)$. Thus, the outer branches do not open linearly in that direction.
\end{proof}

Proposition~\ref{prop:generalcriteria} and the wave-packet estimate in Theorem~\ref{thm:effective-dynamics} do not require $D_4$ symmetry. Under the additional symmetry of $A$ and $b$, the next result makes the cone isotropic and identifies its singlet and doublet. The positivity of $G$ remains a separate hypothesis. The constitutive interpretation in Proposition~\ref{prop:dz-bulk} further requires the singlet to have nonzero cell mean.

\begin{proposition}\label{prop:d4criteria}
Under Assumption~\ref{assump:threefold}, suppose that $A$ and $b$ are $D_4$-invariant in the sense of \eqref{symmetryk=2} and \eqref{symmetryk=0}, respectively, and that $G>0$. Then,
\begin{equation}\label{eq:d4criteria:gamma}
    G=\gamma^2I_2,\qquad \gamma:=\sqrt{G_{11}}>0.
\end{equation}
Moreover, let $\{f_0,f_1,f_2\}\subset H^1_\#(Y;\mathbb{R})$ be a real basis of the eigenspace corresponding to $\lambda_*$ at $\alpha=0$, normalized by \eqref{b-orthogonality}, and let $\{\widetilde f_0,\widetilde f_1,\widetilde f_2\}$ be the real basis given in \eqref{ftildetransform}. There exists a one-dimensional real representation
\begin{equation}\label{eq:d4criteria:character}
    \chi:D_4\rightarrow\{1,-1\}
\end{equation}
such that $\widetilde f_0$ is $\chi$-equivariant, and $(\widetilde f_1, \widetilde f_2)^\top$ is $\chi \otimes \rho_1$-equivariant, where $\rho_1$ is the natural tensor representation defined in \eqref{tensorrepk} for $k=1$.
\end{proposition}

\begin{remark}
    Since the representation $\chi \otimes \rho_1$ is equivalent to $\rho_1$, we may find an $O\in \mathrm{O}(2)$ such that
    \begin{equation}
        \begin{pmatrix}
            g_1 \\g_2
        \end{pmatrix} = O \begin{pmatrix}
            \widetilde f_1 \\ \widetilde f_2
        \end{pmatrix} 
    \end{equation}
    and $(g_1, g_2)^\top$ is $\rho_1$-equivariant. The reduced matrices in this new basis must be conjugated by the same basis change.
\end{remark}

\begin{proof}[Proof of Proposition~\ref{prop:d4criteria}]
Let $\mathcal E_*=\operatorname{span}_{\mathbb{R}}\{f_0,f_1,f_2\}$ be the real eigenspace corresponding to $\lambda_*$, equipped with the $\langle\cdot,\cdot\rangle_{L^2_b(Y)}$ inner product. The trivial representation $\rho_0$ induces the action 
\begin{equation}
    \pi:=\cdot_{\rho_0}:D_4\rightarrow O(\mathcal{E}_*), \qquad [\pi(Q) u](x) :=[Q\cdot_{\rho_0} u](x)= u(Q^\top x),
\end{equation}
where $O(\mathcal{E}_*)$ denotes the group of linear isometries of $\mathcal{E}_*$ with respect to $\langle\cdot,\cdot\rangle_{L^2_b(Y)}$; see \eqref{inducedaction}. Restrict $F$ in \eqref{defFC3toL2b} to $\mathbb R^3$. The restriction is an isometry onto $\mathcal E_*$, which we still denote by $F$. Then, the pullback
\begin{equation}
    \widehat \pi :=F^*\circ\pi\circ F:D_4\rightarrow O(3), \qquad \widehat \pi(Q):= F^*\circ\pi(Q)\circ F
\end{equation}
gives an orthogonal representation of $D_4$. 

By the definition of $F$, one has $F^*F=\mathrm{Id}_{\mathbb{R}^3}$ and $FF^*= \mathrm{Id}_{\mathcal{E}_*}$. Thus,
\begin{equation}
    F\widehat \pi(Q)=FF^*\pi (Q)F =\pi (Q)F.
\end{equation}
Applying this identity to the standard basis vector $\mathbf{e}_s$, and using $F\mathbf{e}_s = f_s$, gives
\begin{equation}
    f_s(Q^\top x)= [\pi(Q)F\mathbf{e}_s](x) = [F\widehat \pi(Q)\mathbf{e}_s](x)= \sum_{r=0}^2 \widehat \pi(Q)_{rs} f_r(x).
\end{equation}
This is equivalent to
\begin{equation}\label{transformofQ}
    \begin{pmatrix}
        f_0(Qx) \\ f_1(Qx) \\ f_2(Qx)
    \end{pmatrix}=\widehat \pi(Q)\begin{pmatrix}
        f_0(x) \\ f_1(x) \\ f_2(x)
    \end{pmatrix}.
\end{equation}
Changing variables $x=Qy$ in \eqref{defofCj} and using \eqref{transformofQ} and $A(Qy)=QA(y)Q^\top$ gives
\begin{equation}\label{covarianceC}
    C(Q\alpha)=\widehat \pi(Q)C(\alpha) \widehat \pi(Q)^\top, \qquad C(\alpha):=\alpha_1 C_1 + \alpha_2 C_2.
\end{equation}
This, combined with $\alpha^\top G\alpha = \frac{1}{2}\operatorname{Tr}(C(\alpha)^\top C(\alpha))$, implies
\begin{equation}
    G=Q^\top GQ, \qquad Q\in D_4.
\end{equation}
A symmetric $2\times2$ matrix commuting with the quarter-turn is scalar. Since $G$ is a Gram matrix, the scalar is nonnegative, and it is positive under the hypothesis $G>0$. This gives \eqref{eq:d4criteria:gamma}.

Let $O\in\mathrm{SO}(3)$ be constructed in Proposition~\ref{prop:generalcriteria}, and set
\begin{equation}
    \widetilde \pi(Q):=O^\top\widehat\pi(Q)O, \qquad \widetilde C_j:=O^\top C_jO, \qquad
    \widetilde C(\alpha):=O^\top C(\alpha)O.
\end{equation}
Then $\widetilde\pi$ is the representation in the transformed basis $\{\widetilde f_0,\widetilde f_1, \widetilde f_2\}$. The normal form of Proposition~\ref{prop:generalcriteria} and $G=\gamma^2I_2$ give
\begin{equation}\label{eq:d4criteria:normalform-proof}
    \widetilde C(\alpha)=\gamma\begin{pmatrix}
              0&\alpha^\top\\
              -\alpha&0_{2\times2}
          \end{pmatrix}.
\end{equation}
Moreover, the covariance identity \eqref{covarianceC} becomes
\begin{equation}\label{eq:d4criteria:covariance-proof}
    \widetilde C(Q\alpha)=\widetilde\pi(Q)\widetilde C(\alpha)\widetilde\pi(Q)^\top.
\end{equation}

We now give the irreducible decomposition of $\widetilde\pi$. Set
\begin{equation}
    L:=\operatorname{span}_{\mathbb{R}}\{\mathbf{e}_0\}, \qquad H:=\operatorname{span}_{\mathbb{R}}\{\mathbf{e}_1,\mathbf{e}_2\}.
\end{equation}
By the normal form \eqref{eq:d4criteria:normalform-proof}, $H$ is spanned by the real kernel lines of $\widetilde C(\alpha)$ for $\alpha \neq 0$, \emph{i.e.},
\begin{equation}
    H = \operatorname{span}_{\mathbb R}\biggl(\bigcup_{\alpha\ne0}\ker_{\mathbb R} \widetilde C(\alpha)\biggr).
\end{equation}
On the other hand, \eqref{eq:d4criteria:covariance-proof} implies $\widetilde\pi(Q) \ker_{\mathbb R} \widetilde C(\alpha)=\ker_{\mathbb R} \widetilde C(Q\alpha)$. Thus, $\widetilde\pi(Q)$ permutes the kernel lines spanning $H$, so $H$ is invariant. Orthogonality then implies that $L=H^\perp$ is invariant as well.

Consequently,
\begin{equation}
    \widetilde\pi(Q) = \begin{pmatrix}
        \chi(Q) & 0 \\ 0 & D(Q)
    \end{pmatrix}, \qquad \chi(Q) \in \{\pm1\}, \quad D(Q)\in \mathrm{O}(2).
\end{equation}
Substituting this block decomposition into \eqref{eq:d4criteria:covariance-proof} and comparing the blocks yields
$D(Q)=\chi(Q)Q$. Finally, the representation $Q\mapsto\chi(Q)Q$ on $H$ is irreducible. Indeed, any invariant real line would be invariant under $\chi(R)R$, and hence under the rotation $R$ through $\pi/2$, which is impossible. This irreducible decomposition and \eqref{transformofQ} give the $\chi$-equivariance of $\widetilde f_0$ and the $\chi\otimes \rho_1$-equivariance of $(\widetilde f_1, \widetilde f_2)$, and hence complete the proof.
\end{proof}

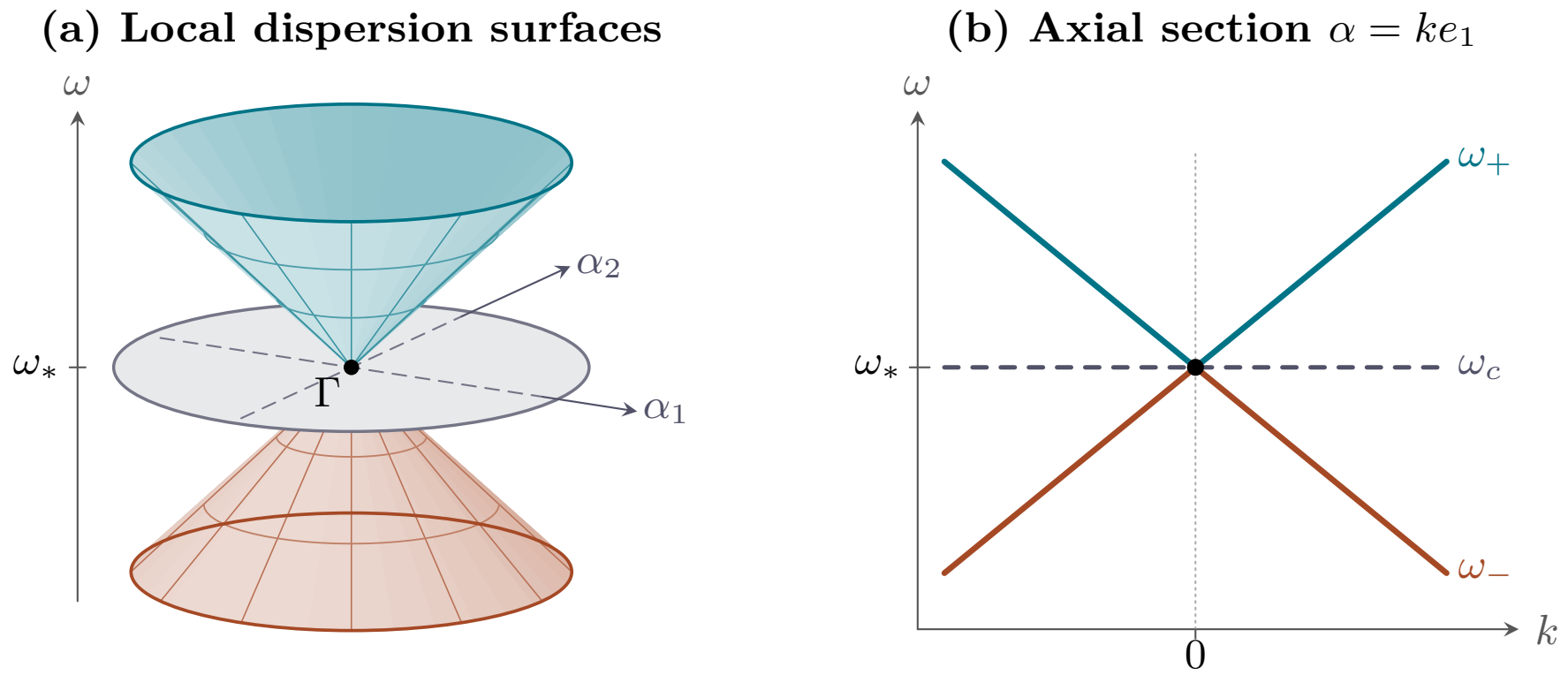
\begin{figure}[tbp]
\centering
\begin{tikzpicture}[x=1cm,y=1cm,font=\small,>=stealth,line join=round,line cap=round]
\definecolor{upperband}{HTML}{007487}
\definecolor{lowerband}{HTML}{A54925}
\definecolor{centerband}{HTML}{515268}
\begin{scope}
  \node[font=\small\bfseries] at (0,2.52) {(a) Local dispersion surfaces};
  \draw[->,black!66] (-2.05,-1.75)--(-2.05,1.92) node[above] {$\omega$};
  \draw[black!66] (-2.11,0)--(-1.99,0);
  \node[left,inner sep=2pt] at (-2.10,0) {$\omega_*$};
  \def\plotcone#1#2{%
    \foreach \j in {0,...,31}{%
      \pgfmathsetmacro{\angA}{11.25*\j}%
      \pgfmathsetmacro{\angB}{11.25*(\j+1)}%
      \pgfmathtruncatemacro{\mix}{33+12*cos(\angA-50)}%
      \path[fill=#2!\mix!white,draw=#2!\mix!white,line width=.32pt]
        (0,0)--({1.65*cos(\angA)},{.44*sin(\angA)+1.53*(#1)})
        --({1.65*cos(\angB)},{.44*sin(\angB)+1.53*(#1)})--cycle;%
    }%
    \foreach \ang in {180,210,...,360}{%
      \draw[#2!78,line width=.3pt] (0,0)--
            ({1.65*cos(\ang)},{.44*sin(\ang)+1.53*(#1)});%
    }%
    \foreach \rad in {.34,.67}{%
      \draw[#2!70,line width=.35pt]
       plot[domain=180:360,samples=31,variable=\a]
         ({1.65*\rad*cos(\a)},{\rad*(.44*sin(\a)+1.53*(#1))});%
    }%
    \draw[#2,line width=.7pt] (0,{1.53*(#1)}) ellipse[x radius=1.65,y radius=.44];%
  }
  \plotcone{-1}{lowerband}
  \path[fill=centerband!13,draw=centerband!80,line width=.65pt]
        (0,0) ellipse[x radius=1.78,y radius=.48];
  \draw[densely dashed,centerband!80] (-1.43,.22)--(1.43,-.22) (-.825,-.381)--(.825,.381);
  \draw[->,centerband] (1.43,-.22)--(2.14,-.33) node[right,inner sep=1pt] {$\alpha_1$};
  \draw[->,centerband] (.825,.381)--(1.64,.757) node[right,inner sep=1pt] {$\alpha_2$};
  \plotcone{1}{upperband}
  \fill (0,0) circle[radius=1.65pt];
  \node[below left,inner sep=2pt] at (0,0) {$\Gamma$};
\end{scope}
\begin{scope}[xshift=6.32cm]
  \node[font=\small\bfseries] at (.12,2.52) {(b) Axial section $\alpha=k e_1$};
  \draw[->,black!66] (-2.08,-1.96)--(-2.08,1.92) node[above] {$\omega$};
  \draw[->,black!66] (-2.08,-1.96)--(2.42,-1.96) node[right] {$k$};
  \draw[black!66] (-2.14,0)--(-1.99,0);
  \node[left,inner sep=2pt] at (-2.12,0) {$\omega_*$};
  \draw[black!32,densely dotted] (0,-1.96)--(0,1.60);
  \draw[black!66] (0,-1.96)--(0,-2.03);
  \node[below,inner sep=1pt] at (0,-2.01) {$0$};
  \draw[upperband,line width=1.25pt] (-1.88,1.54)--(0,0)--(1.88,1.54);
  \draw[lowerband,line width=1.25pt] (-1.88,-1.54)--(0,0)--(1.88,-1.54);
  \draw[centerband,line width=.95pt,dashed] (-1.88,0)--(1.88,0);
  \node[anchor=west,upperband,inner sep=2pt] at (1.89,1.54) {$\omega_+$};
  \node[anchor=west,centerband,inner sep=2pt] at (1.89,0) {$\omega_c$};
  \node[anchor=west,lowerband,inner sep=2pt] at (1.89,-1.54) {$\omega_-$};
  \fill (0,0) circle[radius=1.8pt];
\end{scope}
\end{tikzpicture}
\caption{The square-symmetric threefold Dirac point at $(\Gamma,\omega_*)$ in the first-order frequency approximation \eqref{eq:square-frequency-cones}. The conical sheets have radial slopes $+v_D$ and $-v_D$. The central sheet is horizontal to first order. On the right, $k$ is signed quasimomentum along $e_1$.}
\label{fig:threefold-cone}
\end{figure}

At a square-symmetric conical crossing, let $\omega_*=\sqrt{\lambda_*}$ and $\omega_j(\alpha)=\sqrt{\lambda_j(\alpha)}$ for $j\in\{-,c,+\}$. Taking the square root in \eqref{diracconegeneralcriteria} gives
\begin{equation}\label{eq:square-frequency-cones}
 \omega_\pm(\alpha)=\omega_*\pm v_D|\alpha|+O(|\alpha|^2),\qquad
 \omega_c(\alpha)=\omega_*+O(|\alpha|^2),\qquad v_D=\frac{\gamma}{2\omega_*}.
\end{equation}
Figure~\ref{fig:threefold-cone} shows the three first-order frequency surfaces and their restriction to a coordinate axis.

\section{\texorpdfstring{Effective dynamics near threefold Dirac points at $\Gamma$}{Effective dynamics near threefold Dirac points at Gamma}}\label{sec:effective}

Effective systems near Dirac-like cones have also been derived by high-frequency homogenization for periodic structures with constrained inclusions \cite{AntonakakisCrasterGuenneau2013}. Here we estimate the error in the positive-frequency evolution near the isolated triple eigenvalue in Assumption~\ref{assump:threefold}. When the Gram matrix $G$ in Section~\ref{sec:bloch} is positive definite, the effective system describes a threefold Dirac point. The error estimate below applies to any isolated positive triple eigenvalue. For $x,\xi\in\mathbb R^2$, we use the unitary Fourier transform
\begin{equation}
    \widehat{f}(\xi)=\frac{1}{2\pi} \int_{\mathbb{R}^2} \mathrm{e}^{-\mathrm{i}x\cdot \xi} f(x)\,dx, \qquad f(x)=
    \frac{1}{2\pi} \int_{\mathbb{R}^2} \mathrm{e}^{\mathrm{i}x\cdot \xi} \widehat{f}(\xi)\,d\xi.
\end{equation}
For $k\geq 0$, we define the Sobolev (semi-)norms by 
\begin{equation}
    \|f\|_{\dot{H}^k}^2:= \int_{\mathbb{R}^2} |\xi|^{2k} |\widehat{f}(\xi)|^2\,d\xi, \qquad \|f\|_{H^k}^2:= \int_{\mathbb{R}^2} (1+|\xi|^{2})^k |\widehat{f}(\xi)|^2\,d\xi.
\end{equation}

Fix a real basis $\{f_0,f_1,f_2\}\subset H^1_\#(Y;\mathbb{R})$ of the eigenspace corresponding to $\lambda_*$ at $\alpha=0$, normalized by \eqref{b-orthogonality}.

Given a bandwidth parameter $\delta>0$ and a profile $\mathbf a=(a_0,a_1,a_2)^\top$ defined on $(X,T) \in \mathbb{R}^2\times [0,\infty)$, where each $a_r(\cdot,T)$ is a Schwartz function, the wave packet concentrated near $\lambda_*$ is defined by
\begin{equation}
    \mathrm{WP}_\delta^{\mathbf a}(x,t):= \delta \sum_{r=0}^2 f_r(x) a_r(\delta x,\delta t),
\end{equation}
where $X=\delta x$ and $T=\delta t$ are the slow variables. The prefactor $\delta$ is the two-dimensional $L^2$ normalization.

Let the Bloch frequency $\omega_n(\alpha)$ be given by
\begin{equation}
    \omega_n(\alpha):=\sqrt{\lambda_n(\alpha)},
\end{equation}
where $\lambda_n(\alpha)$ denotes the $n$-th Bloch eigenvalue at $\alpha$. In particular, we set $\omega_*:=\sqrt{\lambda_*}$. Define the effective Dirac symbol $\mathsf D(\xi)$ by
\begin{equation}
    \mathsf D(\xi):= -\frac{\mathrm{i}}{2 \omega_*} C(\xi) = -\frac{\mathrm{i}}{2 \omega_*}(\xi_1 C_1 + \xi_2 C_2),
\end{equation}
where $C_j$ is defined by \eqref{defofCj}.

Given $\Env_0:\mathbb{R}^2\rightarrow \mathbb{C}^3$ in the Schwartz class $\mathcal S(\mathbb R^2;\mathbb C^3)$, let $\Env=(E^{(0)},E^{(1)},E^{(2)})^\top$ denote the evolving vector envelope with initial profile $\Env_0$. It solves
\begin{equation}\label{effectiveequationofE}
    \mathrm{i}\partial_T \Env = \mathsf D(D_X) \Env, \qquad \Env(X,0)=\Env_0(X), \qquad D_X:= - \mathrm{i}\nabla_X.
\end{equation}
Equivalently,
\begin{equation}\label{evolutionequationE}
    \widehat{\Env}(\xi,T)= \mathrm{e}^{-\mathrm{i}T\mathsf D(\xi)} \widehat{\Env_0}(\xi).
\end{equation}
Since the symbol is Hermitian, this evolution preserves every Fourier Sobolev norm. 

On $L_b^2(\mathbb R^2):=L^2(\mathbb R^2,b\,dx)$ the sesquilinear form $\int_{\mathbb{R}^2} A\nabla u\cdot \overline{\nabla v}\,dx$, for $u,v \in H^1(\mathbb{R}^2)$, defines a nonnegative self-adjoint operator $L$. Closedness follows from ellipticity and \eqref{eq:cell-control} summed over the translated unit cells. The corresponding wave equation is
\begin{equation}
    \partial_t^2 u + Lu=0, \qquad \mathrm{equivalently,} \quad b\partial_t^2 u -\nabla \cdot (A\nabla u)=0 \quad \mathrm{in} \ \mathbb{R}^2.
\end{equation}
Write $\mathrm{WP}^{\Env_0}_\delta:=\mathrm{WP}^{\Env}_\delta(\cdot,0)$. We take the initial displacement to be this packet and the initial velocity to be $-\mathrm i\sqrt L\,\mathrm{WP}^{\Env_0}_\delta$. The solution is consequently given by the positive-frequency propagator $e^{-\mathrm it\sqrt L}$. Define its approximation by
\begin{equation}
    u_\delta^{\mathrm{app}}(x,t):=\mathrm{e}^{-\mathrm{i}t\omega_* } \mathrm{WP}^{\Env}_\delta (x,t),
\end{equation}
where $\Env$ satisfies the effective Dirac equation \eqref{effectiveequationofE} with initial data $\Env_0$. Define the approximation error by
\begin{equation}
    u_\delta(\cdot,t):=\mathrm{e}^{-\mathrm{i}t\sqrt L}\mathrm{WP}^{\Env_0}_\delta,\qquad e_\delta(\cdot,t):=u_\delta(\cdot,t)-u_\delta^{\mathrm{app}}(\cdot,t).
\end{equation}

The next theorem bounds this error on time intervals that may grow as $\delta$ tends to zero.

To track the dependence on the material coefficients, we introduce
\begin{equation}\label{coefficientbound}
    s^2:= \sup_{0\neq v \in H^1_\#(Y)} \frac{\int_Y \lambda_{\max}(A(x))|v(x)|^2\,dx}{ \lambda_*^{-1} a_0(v,v)+\|v\|_{L^2_b(Y)}^2 }, 
\end{equation}
where $\lambda_{\max}(A(x))$ denotes the maximal eigenvalue of $A(x)$ and $a_0$ is defined in \eqref{defaalphauv}. Let $\Delta$ be the distance between $\lambda_*$ and the remaining Bloch eigenvalues at $\alpha=0$. Define its relative value by
\begin{equation}
    d:= \min \left\{ 1, \frac{\Delta}{\lambda_*} \right\}. 
\end{equation}

Fix $c=1/100$ and introduce the quasimomentum radius
\begin{equation}\label{radiusassump}
 r:=\min\{1,c\,\omega_*d/s\}.
\end{equation}

\begin{theorem}\label{thm:effective-dynamics}
Under Assumption~\ref{assump:threefold}, there is an absolute constant $C>0$ and, for each integer $N\ge2$, a constant $C_N>0$ depending only on $N$, such that for every Schwartz initial profile $\Env_0$, $\delta>0$ and $t\ge0$, 
\begin{equation}\label{eq:effective-dynamics-error}
 \|e_\delta(\cdot,t)\|_{L_b^2(\mathbb R^2)}\le C\frac{s\delta}{\omega_*d}\|\Env_0\|_{\dot H^1}
+C\frac{s^2\delta^2t}{\omega_*d}\|\Env_0\|_{\dot H^2}+C_N\left(\frac{\delta}{r}\right)^N\|\Env_0\|_{\dot H^N}.
\end{equation}
\end{theorem}

\begin{remark}
For a family of media and fixed initial profile $\Env_0$, the right-hand side of \eqref{eq:effective-dynamics-error} converges to zero in the time range $0\le t\le t_\delta$, provided
\begin{equation}
    \delta = o (r), \qquad \mathrm{and} \qquad t_\delta=o\left(\frac{\omega_*d}{s^2\delta^2}\right).
\end{equation}
Indeed, the definition of $r$ gives $s\delta/(\omega_* d)\leq c\delta/r$.
\end{remark}

\subsection{Floquet--Bloch representation of the discrepancy}

The weighted inner product for the Floquet--Bloch decomposition is
\begin{equation}
    \langle u,v\rangle_{L^2_b(\mathbb{R}^2)} := \int_{\mathbb{R}^2} bu\overline{v}\,dx.
\end{equation}

For $u\in C_c^\infty(\mathbb R^2)$, define the two-dimensional Floquet--Bloch transform by
\begin{equation}
    (\mathcal Bu)(x,\alpha):=\frac1{2\pi}\sum_{n\in\mathbb Z^2}\mathrm{e}^{-\mathrm{i}\alpha\cdot(x+n)}u(x+n),
    \qquad x\in \mathbb{R}^2,\quad \alpha\in Y^*\simeq [-\pi,\pi)^2.
\end{equation}
For every $\alpha\in Y^*$, $(\mathcal Bu)(\cdot,\alpha)$ is $\mathbb{Z}^2$-periodic, and
\begin{equation}\label{Blochgradient}
    \mathcal{B}(\nabla u)=(\nabla_x +\mathrm{i}\alpha) \mathcal{B}u.
\end{equation}
By the Floquet--Bloch theory (see, for instance, \cite{ammari2026wave}), one has the inverse formula
\begin{equation}\label{FBinversion}
    u(x)=\frac{1}{2\pi} \int_{Y^*} \mathrm{e}^{\mathrm{i}\alpha\cdot x} (\mathcal Bu)(x,\alpha)\,d\alpha,
\end{equation}
and the Plancherel identity
\begin{equation}\label{Plancherel}
    \|u\|_{L^2_b(\mathbb{R}^2)}^2 = \int_{Y^*} \|(\mathcal Bu)(\cdot,\alpha)\|_{L^2_b(Y)}^2\,d\alpha.
\end{equation}

By density, the transform extends unitarily to $L_b^2(\mathbb R^2)$, with fibers defined for almost every $\alpha$. For such $\alpha$, choose a complete $L_b^2(Y)$-orthonormal system 
$\{u_n(\cdot,\alpha)\}_{n\geq1}$ such
that
\begin{equation}
    L(\alpha)u_n(\cdot,\alpha)=\lambda_n(\alpha)u_n(\cdot,\alpha).
\end{equation}
The eigenfunctions are chosen measurably in $\alpha$ and extended periodically in the spatial variable.

For almost every $x$, write $F_xz=\sum_{r=0}^2 f_r(x)z_r$ for the row obtained by evaluating the periodic modes in \eqref{defFC3toL2b}. We first compute the transform of the wave packet.

\begin{lemma}\label{lem:FBtransofWP}
    We have
    \begin{equation}
        \mathcal{B}(\mathrm{WP}_\delta^{\Env})(x,\alpha,t)= \delta^{-1}\sum_{k\in \mathbb{Z}^2} \mathrm{e}^{2\pi \mathrm{i}k\cdot x}
    F_x\left[ \mathrm{e}^{-\mathrm{i}t\mathsf D(\alpha+2\pi k)}\widehat{\Env_0} \left( \frac{\alpha+2\pi k}{\delta}\right)\right].
    \end{equation}
\end{lemma}

\begin{proof}
    The periodicity of each $f_r$ gives, term by term in the Floquet--Bloch transform,
    \begin{equation}
        (\mathcal B\,\mathrm{WP}_\delta^{\Env})(x,\alpha,t)
        =F_x\big[(\mathcal B\Env^\delta)(x,\alpha,t)\big],
    \end{equation}
    where the transform acts componentwise on the rescaled envelope
    \begin{equation}
        \Env^\delta(x,t):=\delta \Env(\delta x,\delta t).
    \end{equation}
    As a $\mathbb{Z}^2$-periodic function, $\mathcal{B}(\Env^\delta)(\cdot,\alpha,t)$ can be expanded in a Fourier series. Its $k$-th Fourier coefficient is computed by
    \begin{equation}
    \begin{aligned}
        &\int_Y \mathcal{B}(\Env^\delta)(x,\alpha,t) \mathrm{e}^{-2\pi \mathrm{i}k\cdot x}\,dx \\
        &\qquad  = \frac{\delta}{2\pi} \sum_{n\in \mathbb{Z}^2} \int_Y \mathrm{e}^{-\mathrm{i}\alpha\cdot(x+n)}
        \Env(\delta (x+n),\delta t)\mathrm{e}^{-2\pi \mathrm{i}k\cdot x}\,dx \\
        &\qquad =\frac{\delta}{2\pi} \int_{\mathbb{R}^2}  \mathrm{e}^{-\mathrm{i}(\alpha+2\pi k)\cdot z}
        \Env(\delta z,\delta t) \,dz \\
        &\qquad = \delta^{-1} \widehat{\Env} \left( \frac{\alpha+2\pi k}{\delta},\delta t \right).
    \end{aligned}
    \end{equation}

Using \eqref{evolutionequationE} and 
\begin{equation}
    \delta t\mathcal{D}\left( \frac{\alpha+2\pi k}{\delta} \right) = t\mathsf D(\alpha+2\pi k),
\end{equation}
we get
\begin{equation}
    \mathcal{B}(\Env^\delta)(x,\alpha,t) = \delta^{-1}\sum_{k\in \mathbb{Z}^2} \mathrm{e}^{2\pi \mathrm{i}k\cdot x}
    \mathrm{e}^{-\mathrm{i}t\mathsf D(\alpha+2\pi k)}\widehat{\Env_0} \left( \frac{\alpha+2\pi k}{\delta}\right).
\end{equation}
Since $\Env_0$ is in the Schwartz class, the above series absolutely converges.
\end{proof}

The following proposition gives a representation of the discrepancy $e_\delta$ based on the Floquet--Bloch transform, which is the starting point of the subsequent error analysis.

\begin{proposition}\label{lem:discrepancy-bloch}
For $n\ge1$ and $k\in\mathbb Z^2$, we define the row $B_{n,k}(\alpha)\in\mathbb C^{1\times3}$ by
\begin{equation}\label{eq:lemma31-rows}
 B_{n,k}(\alpha)z := \int_Y b(x) \mathrm{e}^{2\pi \mathrm{i} k\cdot x} (Fz)(x) \overline{u_n(x,\alpha)}\,dx .
\end{equation}
Then, we have
\begin{equation}\label{eq:lemma31-physical-expansion}
 e_\delta(x,t)=\frac{1}{2\pi} \sum_{n\geq 1} \int_{Y^*} \mathrm{e}^{\mathrm{i}\alpha\cdot x}   c_{\delta,n}(\alpha,t)u_n(x,\alpha)\,d\alpha,
\end{equation}
where
\begin{equation}\label{defcdeltan}
\begin{aligned}
    c_{\delta,n}(\alpha,t)&:=\big\langle (\mathcal Be_\delta)(\cdot,\alpha,t), u_n(\cdot,\alpha)\big\rangle_{L^2_b(Y)} \\
    &\ =\frac{\mathrm{e}^{-\mathrm{i}t\omega_n(\alpha)}}{\delta} \sum_{k\in \mathbb{Z}^2} B_{n,k}(\alpha)
        \left( I_3- \mathrm{e}^{\mathrm{i}t\{\omega_n(\alpha)I_3-\omega_*I_3-\mathsf D(\alpha+2\pi k)\} } \right) \widehat{\Env_0} \left( \frac{\alpha+2\pi k}{\delta}\right) .
\end{aligned}
\end{equation}
The expansion \eqref{eq:lemma31-physical-expansion} converges in $L^2_b(\mathbb{R}^2)$. Moreover,
\begin{equation}\label{eq:lemma31-error-norm}
 \|e_\delta(\cdot,t)\|_{L_b^2(\mathbb R^2)}^2 =\int_{Y^*}\sum_{n\geq1}|c_{\delta,n}(\alpha,t)|^2\,d\alpha.
\end{equation}
\end{proposition}

\begin{proof}
    By the Floquet--Bloch inverse formula \eqref{FBinversion} and the Plancherel identity \eqref{Plancherel}, we get
    \begin{equation}
        \begin{aligned}
            e_\delta (x,t)&= \frac{1}{2\pi} \int_{Y^*} \mathrm{e}^{\mathrm{i}\alpha\cdot x} (\mathcal Be_\delta)(x,\alpha,t)\,d\alpha \\
            &=\frac{1}{2\pi} \sum_{n\geq 1} \int_{Y^*} \mathrm{e}^{\mathrm{i}\alpha\cdot x}   c_{\delta,n}(\alpha,t)u_n(x,\alpha)\,d\alpha
        \end{aligned}
    \end{equation}
and \eqref{eq:lemma31-error-norm}. It remains to compute $c_{\delta,n}$. By \eqref{Blochgradient} and functional calculus, we have
    \begin{equation}
        \mathcal{B}u_\delta= \mathcal{B}\mathrm{e}^{-\mathrm{i}t\sqrt{L}}(\mathrm{WP}_\delta^{\Env_0}) =\mathrm{e}^{-\mathrm{i}t\sqrt{L(\alpha)}} \mathcal{B}(\mathrm{WP}_\delta^{\Env_0}).
    \end{equation}
Hence,
\begin{equation}
    \begin{aligned}
        c_{\delta,n}(\alpha,t)&= \big\langle \mathrm{e}^{-\mathrm{i}t\sqrt{L(\alpha)}} \mathcal{B}(\mathrm{WP}_\delta^{\Env_0})(\cdot,\alpha) - \mathrm{e}^{-\mathrm{i}t\omega_*} \mathcal{B}(\mathrm{WP}_\delta^{\Env})(\cdot,\alpha,t) , u_n(\cdot,\alpha)\big\rangle_{L^2_b(Y)} \\
        &= \frac{\mathrm{e}^{-\mathrm{i}t\omega_n(\alpha)}}{\delta} \sum_{k\in \mathbb{Z}^2} B_{n,k}(\alpha)
        \left( I_3- \mathrm{e}^{\mathrm{i}t\{\omega_n(\alpha)I_3-\omega_*I_3-\mathsf D(\alpha+2\pi k)\} } \right) \widehat{\Env_0} \left( \frac{\alpha+2\pi k}{\delta}\right) .
    \end{aligned}
\end{equation}
Here, we used Lemma \ref{lem:FBtransofWP}.
\end{proof}

\subsection{Quantitative estimates for the Bloch coefficients}

Let $\mathcal I_*=\{n_*,n_*+1,n_*+2\}$ denote the indices of the triple at $\Gamma$. We estimate the component of the wave packet in this cluster and its projection onto the remaining bands separately.

In the series expansion \eqref{defcdeltan} of $c_{\delta,n}$, the terms with $k\neq0$ decay faster than any power of $\delta$ as $\delta\to0$, since $\widehat{\Env_0}$ is in the Schwartz class. The following lemma estimates the principal contribution of the series, corresponding to the component $k=0$.

\begin{lemma}\label{componentk0est}
    There is an absolute constant $C>0$ such that, for $|\alpha|<r$, $t\ge0$
    and every $z\in \mathbb{C}^3$, 
    \begin{equation}\label{estimateBnnotin}
        \sum_{n\notin \mathcal I_*} |B_{n,0}(\alpha)z|^2 \leq C\frac{s^2}{\lambda_*d^2} |\alpha|^2|z|^2,
    \end{equation}
    and
    \begin{equation}\label{estimateBnin}
        \sum_{n\in \mathcal I_* } \left| B_{n,0}(\alpha)
        \left( I_3- \mathrm{e}^{\mathrm{i}t\{\omega_n(\alpha)I_3-\omega_*I_3-\mathsf D(\alpha)\} } \right) z\right|^2 \leq C\frac{s^4}{\lambda_* d^2} t^2|\alpha|^4 |z|^2.
    \end{equation}
\end{lemma}

\begin{proof}
The weak Bloch eigenvalue equations \eqref{blocheigenprobL} imply
\begin{equation}
    a_\alpha(Fz,u_n(\cdot,\alpha))=\lambda_n(\alpha)\langle Fz,u_n(\cdot,\alpha)\rangle_{L_b^2(Y)},
 \qquad a_0(Fz,u_n(\cdot,\alpha))=\lambda_*\langle Fz,u_n(\cdot,\alpha)\rangle_{L_b^2(Y)}.
\end{equation}
Therefore,
\begin{equation}\label{eq:qbc-weak-overlap}
 (a_\alpha-a_0)(Fz,u_n(\cdot,\alpha)) =(\lambda_n(\alpha)-\lambda_*)B_{n,0}(\alpha)z.
\end{equation}
Next, we estimate $(a_\alpha-a_0)(Fz,u_n(\cdot,\alpha))$ and $\lambda_n(\alpha)-\lambda_*$.

\textit{Step 1.} Introduce the energy norm
\begin{equation}
    \|v\|_{\mathcal E}^2 :=a_0(v,v)+\lambda_*\|v\|_{L_b^2(Y)}^2, \qquad v\in H_\#^1(Y).
\end{equation}
The coefficient bound \eqref{coefficientbound} becomes
\begin{equation}
\label{eq:qbc-energy-multiplication}
 \int_Y\lambda_{\max}(A)|v|^2\,dx \le\frac{s^2}{\lambda_*}\|v\|_{\mathcal E}^2, \qquad v\in H_\#^1(Y).
\end{equation}

Expand the form $a_\alpha$ as
\begin{equation}
    a_\alpha-a_0=\ell_\alpha+p_\alpha,
\end{equation}
where
\begin{equation}
    \begin{aligned}
      \ell_\alpha(u,v) &=\mathrm{i}\int_Y\left[u(A\alpha)\cdot\nabla\overline v
      -\overline v(A\alpha)\cdot\nabla u\right]dx,\\
      p_\alpha(u,v) &=\int_Y(\alpha^\top A\alpha)u\overline v\,dx.
    \end{aligned}
\end{equation}
The energy norm and \eqref{eq:qbc-energy-multiplication} give
\begin{equation}
    \|A^{1/2}\nabla v\|_{L^2(Y;\mathbb C^2)} \le\|v\|_{\mathcal E}, \qquad
 \|A^{1/2}\alpha v\|_{L^2(Y;\mathbb C^2)} \le\frac{s|\alpha|}{\omega_*}\|v\|_{\mathcal E}.
\end{equation}
Therefore,
\begin{equation}\label{eq:qbc-form-parts}
|\ell_\alpha(u,v)|\le\frac{2s|\alpha|}{\omega_*}\|u\|_{\mathcal E}\|v\|_{\mathcal E},
\qquad |p_\alpha(u,v)|\le\frac{s^2|\alpha|^2}{\lambda_*} \|u\|_{\mathcal E}\|v\|_{\mathcal E}.
\end{equation}
Set
\begin{equation}
    h:=\frac{2s|\alpha|}{\omega_*} +\frac{s^2|\alpha|^2}{\lambda_*}.
\end{equation}
We consequently have
\begin{equation}\label{eq:qbc-form-perturbation}
 |(a_\alpha-a_0)(u,v)| \le h\|u\|_{\mathcal E}\|v\|_{\mathcal E},
\end{equation}
and
\begin{equation}\label{eq:qbc-form-comparison}
 (1-h)\|v\|_{\mathcal E}^2\le a_\alpha(v,v)+\lambda_*\|v\|_{L_b^2(Y)}^2\le(1+h)\|v\|_{\mathcal E}^2.
\end{equation}
Applying the min--max principle to \eqref{eq:qbc-form-comparison} gives, for every index $n\geq 1$,
\begin{equation}\label{eq:qbc-minmax}
 (1-h)(\lambda_n(0)+\lambda_*) \le\lambda_n(\alpha)+\lambda_* \le(1+h)(\lambda_n(0)+\lambda_*).
\end{equation}

With $c=1/100$, the radius bound \eqref{radiusassump} implies
\begin{equation}\label{eq:form-smallness}
    \frac{s|\alpha|}{\omega_*}<cd,\qquad h\le\frac{3s|\alpha|}{\omega_*}<d/16.
\end{equation}

For fixed $z\in\mathbb C^3$, since $L(0)Fz=\lambda_*Fz$ and $F^*F=I_3$,
\begin{equation}\label{eq:qbc-F-energy}
 \|Fz\|_{\mathcal E}=\sqrt{2}\omega_*|z|.
\end{equation}
This, combined with \eqref{eq:qbc-form-perturbation} and the estimate $h\leq 3s|\alpha|/\omega_*$ in \eqref{eq:form-smallness}, yields
\begin{equation}
    |(a_\alpha-a_0)(Fz,v)| \le Cs|\alpha||z|\|v\|_{\mathcal E}.
\end{equation}
By \eqref{eq:qbc-form-comparison} and $h<1/16$, this implies
\begin{equation}\label{eq:qbc-source-dual-bound}
 |(a_\alpha-a_0)(Fz,v)| \le Cs|\alpha||z|  \left[a_\alpha(v,v)+\lambda_*\|v\|_{L_b^2(Y)}^2\right]^{1/2}.
\end{equation}

The family
\begin{equation}
    \left\{\frac{u_n(\cdot,\alpha)}{\sqrt{\lambda_n(\alpha)+\lambda_*}} :n\geq 1 \right\}
\end{equation}
forms an orthonormal basis of $H^1_\#(Y)$ equipped with the inner product
\begin{equation}
    a_\alpha(u,v)+\lambda_*\langle u,v\rangle_{L_b^2(Y)}.
\end{equation}
Applying Bessel's inequality to the continuous antilinear functional $v\mapsto(a_\alpha-a_0)(Fz,v)$, using
\eqref{eq:qbc-source-dual-bound}, gives
\begin{equation}\label{eq:qbc-source-bessel}
 \sum_{n\ge1} \frac{|(a_\alpha-a_0)(Fz,u_n)|^2}{\lambda_n(\alpha)+\lambda_*} \le Cs^2|\alpha|^2|z|^2.
\end{equation}

\textit{Step 2.} In this step, we assume that $n\notin\mathcal I_*$, and prove \eqref{estimateBnnotin}. For index $n\notin\mathcal I_*$, either $\lambda_n(0)\le\lambda_*(1-d)$, or $\lambda_n(0)\ge\lambda_*(1+d)$. Thus, \eqref{eq:qbc-minmax}--\eqref{eq:form-smallness} imply that
\begin{equation}
\label{eq:qbc-outside-separation}
 |\lambda_n(\alpha)-\lambda_*|\ge \lambda_*\frac{13d}{16}, \qquad n\notin\mathcal I_*.
\end{equation}
Moreover, \eqref{eq:qbc-outside-separation} and $\lambda_n(\alpha)+\lambda_*\le |\lambda_n(\alpha)-\lambda_*|+2\lambda_*$ give
\begin{equation}\label{eq:qbc-off-quotient}
 \frac{\lambda_n(\alpha)+\lambda_*}{|\lambda_n(\alpha)-\lambda_*|^2}\le\frac{C}{\lambda_*d^2}.
\end{equation}
Combining \eqref{eq:qbc-weak-overlap} and \eqref{eq:qbc-source-bessel}--\eqref{eq:qbc-off-quotient}, we obtain the desired estimate \eqref{estimateBnnotin}:
\begin{equation}
    \begin{aligned}
\sum_{n\notin\mathcal I_*}|B_{n,0}(\alpha)z|^2 &\le\frac{C}{\lambda_*d^2}
\sum_{n\geq 1}\frac{|(a_\alpha-a_0)(Fz,u_n(\cdot,\alpha))|^2}{\lambda_n(\alpha)+\lambda_*}\\
 &\le C\frac{s^2}{\lambda_*d^2}|\alpha|^2|z|^2.
    \end{aligned}
\end{equation}

\textit{Step 3.} In this step, we assume that $n\in\mathcal I_*$, and prove \eqref{estimateBnin}. For $n\in\mathcal I_*$, $\lambda_n(0)=\lambda_*$, and hence \eqref{eq:qbc-minmax}--\eqref{eq:form-smallness} imply
\begin{equation}\label{eq:qbc-cluster-shift}
 |\lambda_n(\alpha)-\lambda_*| \le2h\lambda_* <\lambda_*\frac{d}{8}.
\end{equation}
Moreover, \eqref{eq:qbc-form-comparison} gives
\begin{equation}
    (1-h)\|u_n(\cdot,\alpha)\|_{\mathcal E}^2 \le a_\alpha(u_n(\cdot,\alpha),u_n(\cdot,\alpha))+\lambda_*
 =\lambda_n(\alpha)+\lambda_*.
\end{equation}
Consequently,
\begin{equation}\label{eq:qbc-cluster-energy}
 \|u_n(\cdot,\alpha)\|_{\mathcal E}\le C\omega_*,
\end{equation}
and
\begin{equation}\label{eq:qbc-frequency-shift}
 |\omega_n(\alpha)-\omega_*| \leq \frac{|\lambda_n(\alpha)-\lambda_*|}{\omega_*}\le Cs|\alpha|.
\end{equation}

Let $P(0)=FF^*$ be the base triple projection and set
\begin{equation}
    \chi_n(\cdot,\alpha) :=(I-P(0))u_n(\cdot,\alpha).
\end{equation}
Decompose $u_n(\cdot,\alpha)=P(0)u_n(\cdot,\alpha)+\chi_n(\cdot,\alpha)$. By
\begin{equation}
    \ell_\alpha(Fz,P(0)u_n(\cdot,\alpha))=-\mathrm{i} B_{n,0}(\alpha)C(\alpha) z = 2\omega_* B_{n,0}(\alpha)\mathsf D(\alpha) z,
\end{equation}
and \eqref{eq:qbc-weak-overlap}, we obtain the following:
\begin{equation}\label{eq:qbc-cancellation-equation}
 (\lambda_n(\alpha)-\lambda_*)B_{n,0}(\alpha)z=2\omega_* B_{n,0}(\alpha)\mathsf D(\alpha) z +\ell_\alpha(Fz,\chi_n(\cdot,\alpha))+p_\alpha(Fz,u_n(\cdot,\alpha)).
\end{equation}
As a consequence, it follows that
\begin{equation}\label{eq:qbc-frequency-cancellation}
 \begin{aligned}
 &B_{n,0}(\alpha)\bigl[(\omega_n(\alpha)-\omega_*)I_3-\mathsf D(\alpha)\bigr]z\\
 &\quad=
 \frac{\ell_\alpha(Fz,\chi_n(\cdot,\alpha))+p_\alpha(Fz,u_n(\cdot,\alpha))}{\omega_n(\alpha)+\omega_*}
 + \frac{\omega_*-\omega_n(\alpha)}{\omega_*+\omega_n(\alpha)}
 B_{n,0}(\alpha) \mathsf D(\alpha) z.
 \end{aligned}
\end{equation}
We then have
\begin{equation}\label{eq:qbc-frequency-cancellation-norm}
 \begin{aligned}
 &\left| B_{n,0}(\alpha)\bigl[(\omega_n(\alpha)-\omega_*)I_3-\mathsf D(\alpha)\bigr]z \right| \\
 &\quad \leq \frac{2s|\alpha|}{\lambda_*} \|Fz\|_{\mathcal{E}} \|\chi_n(\cdot,\alpha)\|_{\mathcal{E}} + 
 \frac{s^2|\alpha|^2}{\omega_*^3} \|Fz\|_{\mathcal{E}} \|u_n(\cdot,\alpha)\|_{\mathcal{E}} + \frac{|\omega_*-\omega_n(\alpha)|}{\omega_*} | B_{n,0}(\alpha) \mathsf D(\alpha) z |\\
 &\quad\leq C\frac{s|\alpha| |z| }{\omega_*} \|\chi_n(\cdot,\alpha)\|_{\mathcal{E}} + 
C \frac{s^2|\alpha|^2|z|}{\omega_*} + C\frac{s|\alpha|}{\omega_*} | B_{n,0}(\alpha) \mathsf D(\alpha) z | ,
 \end{aligned}
\end{equation}
where we used \eqref{eq:qbc-form-parts} in the first inequality, \eqref{eq:qbc-F-energy}, \eqref{eq:qbc-cluster-energy}, and \eqref{eq:qbc-frequency-shift} in the second inequality.

By the definition \eqref{eq:lemma31-rows} of $B_{n,0}(\alpha)$, we have $|B_{n,0}(\alpha)|\leq 1$. Moreover, the definition \eqref{defofCj} of $C_j$ gives $\ell_\alpha(f_s,f_r)=-\mathrm{i}(C(\alpha))_{rs}$, hence
\begin{equation}\label{eq:qbc-linear-compression}
 \ell_\alpha(Fz,Fw)=-\mathrm{i}\langle C(\alpha)z,w\rangle_{\mathbb C^3},
 \qquad |C(\alpha)|\le4s\omega_*|\alpha|.
\end{equation}
The norm bound follows from \eqref{eq:qbc-form-parts} and \eqref{eq:qbc-F-energy}. Therefore, 
\begin{equation}
    \left| B_{n,0}(\alpha)\bigl[(\omega_n(\alpha)-\omega_*)I_3-\mathsf D(\alpha)\bigr]z \right|  \leq 
    C\frac{s|\alpha| |z| }{\omega_*} \|\chi_n(\cdot,\alpha)\|_{\mathcal{E}}  + C\frac{s^2|\alpha|^2}{\omega_*} |z | .
\end{equation}

We claim that
\begin{equation}\label{eq:qbc-complement-energy}
 \|\chi_n(\cdot,\alpha)\|_{\mathcal E}\le C\frac{s|\alpha|}{d}.
\end{equation}
Assuming this bound for the moment, we obtain
\begin{equation}\label{eq:qbc-inline-mismatch}
 \left|B_{n,0}(\alpha)\bigl[(\omega_n(\alpha)-\omega_*)I_3-\mathsf D(\alpha)\bigr]z\right|
 \le C\frac{s^2|\alpha|^2}{\omega_*d}|z|.
\end{equation}
Finally, since the matrix $\mathsf D(\alpha)$ is Hermitian, the Duhamel formula gives 
\begin{equation}
    \begin{aligned}
  &B_{n,0}(\alpha) \left[I_3-\mathrm{e}^{\mathrm{i}t\{(\omega_n(\alpha)-\omega_*)I_3-\mathsf D(\alpha)\}}\right]z\\
 &\quad=-\mathrm{i}\int_0^t \mathrm{e}^{\mathrm{i}\sigma(\omega_n(\alpha)-\omega_*)} B_{n,0}(\alpha)
 \bigl[(\omega_n(\alpha)-\omega_*)I_3-\mathsf D(\alpha)\bigr] \mathrm{e}^{-\mathrm{i}\sigma\mathsf D(\alpha)}z\,d\sigma.
    \end{aligned}
\end{equation}
Since $|\mathrm{e}^{-\mathrm{i}\sigma\mathsf D(\alpha)}z|=|z|$, applying \eqref{eq:qbc-inline-mismatch} yields
\begin{equation}
    \left|B_{n,0}(\alpha) \left[I_3-\mathrm{e}^{\mathrm{i}t\{(\omega_n(\alpha)-\omega_*)I_3-\mathsf D(\alpha)\}}\right]z\right|
 \le Ct\frac{s^2|\alpha|^2}{\omega_*d}|z|.
\end{equation}
Squaring and summing over $n\in\mathcal I_*$ proves \eqref{estimateBnin}.

It remains to prove \eqref{eq:qbc-complement-energy}. Choose a complete $L_b^2(Y)$-orthonormal eigenbasis $e_m$ of $L(0)$, with eigenvalues $\lambda_m(0)$ and cluster index set $\mathcal I_*$. For every $n,m\geq 1$, the weak Bloch eigenvalue equations \eqref{blocheigenprobL} give
\begin{equation}\label{eq:qbc-base-coefficients}
 (\lambda_n(\alpha)-\lambda_m(0))\langle u_n(\cdot,\alpha),e_m\rangle_{L_b^2(Y)}=(a_\alpha-a_0)(u_n(\cdot,\alpha),e_m).
\end{equation}
For $n\in\mathcal I_*$, \eqref{eq:qbc-form-perturbation}, \eqref{eq:form-smallness}, and \eqref{eq:qbc-cluster-energy} give
\begin{equation}
    |(a_\alpha-a_0)(u_n(\cdot,\alpha),v)| \le Cs|\alpha|\|v\|_{\mathcal E}, \qquad v\in H^1_\#(Y).
\end{equation}
Bessel's inequality yields
\begin{equation}\label{eq:qbc-base-source-bessel}
 \sum_{m\geq 1}\frac{|(a_\alpha-a_0)(u_n,e_m)|^2}{\lambda_m(0)+\lambda_*} \le Cs^2|\alpha|^2.
\end{equation}
For $m\notin\mathcal I_*$, and $n\in \mathcal I_*$,
\begin{equation}
     |\lambda_m(0)-\lambda_*|\ge\lambda_*d, \qquad |\lambda_n(\alpha)-\lambda_*|<\lambda_*\frac{d}{8};
\end{equation}
see \eqref{eq:qbc-cluster-shift}. Hence,
\begin{equation}
     |\lambda_m(0)-\lambda_n(\alpha)| \ge\frac78|\lambda_m(0)-\lambda_*|.
\end{equation}
Since
$\lambda_m(0)+\lambda_*\le|\lambda_m(0)-\lambda_*|+2\lambda_*$, it follows that
\begin{equation}\label{eq:qbc-base-quotient}
 \frac{\lambda_m(0)+\lambda_*}{|\lambda_m(0)-\lambda_n(\alpha)|} \le\frac C d.
\end{equation}
Using the spectral expansion of $\chi_n$ and \eqref{eq:qbc-base-coefficients}--\eqref{eq:qbc-base-quotient},
\begin{equation}
    \begin{aligned}
        \|\chi_n\|_{\mathcal E}^2
 &=\sum_{m\notin\mathcal I_*}(\lambda_m(0)+\lambda_*)|\langle u_n(\cdot,\alpha),e_m\rangle_{L_b^2(Y)}|^2\\
 &=\sum_{m\notin\mathcal I_*}\frac{(\lambda_m(0)+\lambda_*)^2}{|\lambda_m(0)-\lambda_n(\alpha)|^2}
    \frac{|(a_\alpha-a_0)(u_n(\cdot,\alpha),e_m)|^2}{\lambda_m(0)+\lambda_*}\\
 &\le C\frac{s^2|\alpha|^2}{d^2}.
    \end{aligned}
\end{equation}
This proves \eqref{eq:qbc-complement-energy}.
\end{proof}

\subsection{Proof of the wave-packet estimate}
We combine the cluster estimate with a bound for the Fourier components outside the central reciprocal cell.
\begin{proof}[Proof of Theorem~\ref{thm:effective-dynamics}]
\medskip\noindent
\textit{Decomposition of the Bloch coefficients.}

For $\alpha\in Y^*=[-\pi,\pi)^2$ and $k\in\mathbb Z^2$, set
\[
 \alpha_k:=\alpha+2\pi k,
 \qquad
 z_{\delta,k}(\alpha)
 :=\widehat{\Env_0}\left(\frac{\alpha_k}{\delta}\right).
\]
The coefficient $c_{\delta,n}$ in Proposition \ref{lem:discrepancy-bloch} can be written as
\begin{equation}
\label{eq:ed-coefficient-decomposition}
c_{\delta,n} =c_{\delta,n}^{(0)}+c_{\delta,n}^{(\ne0)}, \qquad c_{\delta,n}^{(\ne0)}:=\sum_{k\in \mathbb{Z}^2\setminus \{0\} }c_{\delta,n}^{(k)},
\end{equation}
where
\begin{equation}\label{eqsinglekequalzero}
    c_{\delta,n}^{(0)}(\alpha,t):=\frac{\mathrm{e}^{-\mathrm{i}t\omega_n(\alpha)}}{\delta}  B_{n,0}(\alpha)
        \left( I_3- \mathrm{e}^{\mathrm{i}t\{\omega_n(\alpha)I_3-\omega_*I_3-\mathsf D(\alpha)\} } \right) \widehat{\Env_0} \left( \frac{\alpha}{\delta}\right) ,
\end{equation}
and
\begin{equation}\label{eq:ed-single-k}
\begin{aligned}
    &c_{\delta,n}^{(k)}(\alpha,t):= \\
    &\quad \frac{\mathrm{e}^{-\mathrm{i}t\omega_n(\alpha)}}{\delta} B_{n,k}(\alpha)
        \left( I_3- \mathrm{e}^{\mathrm{i}t\{\omega_n(\alpha)I_3-\omega_*I_3-\mathsf D(\alpha+2\pi k)\} } \right) \widehat{\Env_0} \left( \frac{\alpha+2\pi k}{\delta}\right) .
\end{aligned}
\end{equation}

Recall that $F^*F=I_3$. By Parseval's identity in $L_b^2(Y)$ and the definition \eqref{eq:lemma31-rows} of $B_{n,k}(\alpha)$, for every $z\in\mathbb C^3$ and $k\in \mathbb{Z}^2$,
\begin{equation}\label{eq:ed-B-parseval}
 \begin{aligned}
 \sum_{n\ge1}|B_{n,k}(\alpha)z|^2&=\sum_{n\ge1}\left|\int_Y b(x) \mathrm{e}^{2\pi \mathrm{i} k\cdot x} (Fz)(x) \overline{u_n(x,\alpha)}\,dx\right|^2\\
 &=\|\mathrm{e}^{2\pi \mathrm{i} k\cdot x} Fz\|_{L_b^2(Y)}^2=|z|^2.
 \end{aligned}
\end{equation}
Moreover, $\mathsf D(\xi)$ is Hermitian for every real $\xi$, so $|\mathrm{e}^{-\mathrm{i}t\mathsf D(\xi)}z|=|z|$. The identity \eqref{eq:ed-B-parseval} therefore gives
\begin{equation}\label{eq:ed-single-k-global}
 \left( \sum_{n\ge1}|c_{\delta,n}^{(k)}(\alpha,t)|^2 \right)^{1/2}
 \le\frac2\delta \left|\widehat{\Env_0}\left(\frac{\alpha+2\pi k}{\delta}\right) \right|, \qquad k\in\mathbb Z^2.
\end{equation}

\textit{Step 1: the terms with $k\ne0$.} By Minkowski's inequality and \eqref{eq:ed-single-k-global},
\begin{equation}\label{estikneqzero1}
    \begin{aligned}
         \left(\int_{Y^*}\sum_{n\ge1}|c_{\delta,n}^{(\ne0)}(\alpha,t)|^2\,d\alpha \right)^{1/2}
 &\le \sum_{k\ne0} \left( \int_{Y^*}\sum_{n\ge1}|c_{\delta,n}^{(k)}(\alpha,t)|^2\,d\alpha\right)^{1/2}\\
 &\le 2\sum_{k\ne0} \left( \delta^{-2}\int_{Y^*} \left| \widehat{\Env_0}\left(\frac{\alpha+2\pi k}{\delta}\right)
 \right|^2d\alpha \right)^{1/2}\\
 &= 2\sum_{k\ne0} \left( \int_{\delta^{-1} (Y^*+2\pi k)}|\widehat{\Env_0}(q)|^2\,dq \right)^{1/2}.
    \end{aligned}
\end{equation}
For $\alpha\in[-\pi,\pi)^2$, each coordinate satisfies $|\alpha_i+2\pi k_i|\ge\pi|k_i|$. Hence the shifted cells satisfy
\begin{equation}
    |q|\ge\frac{\pi|k|}{\delta} \quad\text{for }q\in \delta^{-1} (Y^*+2\pi k).
\end{equation}
Consequently,
\begin{equation}\label{estikneqzero2}
    \begin{aligned}
         &\sum_{k\ne0} \left( \int_{\delta^{-1} (Y^*+2\pi k)}|\widehat{\Env_0}(q)|^2\,dq \right)^{1/2}\\
 &\quad\le \left(\frac{\delta}{\pi}\right)^N \sum_{k\ne0}|k|^{-N} \left( \int_{\delta^{-1} (Y^*+2\pi k)}
  |q|^{2N}|\widehat{\Env_0}(q)|^2\,dq \right)^{1/2}\\
 &\quad\le \left(\frac{\delta}{\pi}\right)^N \left(\sum_{k\ne0}|k|^{-2N}\right)^{1/2}
 \left( \sum_{k\ne0}\int_{\delta^{-1} (Y^*+2\pi k)}|q|^{2N}|\widehat{\Env_0}(q)|^2\,dq \right)^{1/2}\\
 &\quad\le C_N\delta^N\|\Env_0\|_{\dot H^N}.
    \end{aligned}
\end{equation}
Here, we used the Cauchy--Schwarz inequality in the second inequality and
the convergence of
\begin{equation}
    \sum_{k\ne0}|k|^{-2N}
\end{equation}
in dimension two when $N>1$. Therefore, \eqref{estikneqzero1} and \eqref{estikneqzero2} imply that
\begin{equation}\label{eq:ed-nonzero-k-norm}
 \left(\int_{Y^*}\sum_{n\ge1}|c_{\delta,n}^{(\ne0)}(\alpha,t)|^2\,d\alpha \right)^{1/2}
 \le C_N\delta^N\|\Env_0\|_{\dot H^N}.
\end{equation}

\textit{Step 2: the principal term for $|\alpha|<r$.} 
For $n\in\mathcal I_*$, the estimate \eqref{estimateBnin} in Lemma \ref{componentk0est} gives
\begin{equation}\label{eq:ed-principal-triple}
 \sum_{n\in\mathcal I_*} |c_{\delta,n}^{(0)}(\alpha,t)|^2 \le
 C\frac{s^4t^2|\alpha|^4}{\delta^2\lambda_*d^2} \left|\widehat{\Env_0} \left( \frac{\alpha}{\delta}\right)\right|^2.
\end{equation}
For the complementary indices, set $z=\widehat{\Env_0}(\alpha/\delta)$ and $U=\exp(-\mathrm it\mathsf D(\alpha))$. The scalar factor $\exp(\mathrm it(\omega_n(\alpha)-\omega_*))$ has modulus one. Applying \eqref{estimateBnnotin} separately to $z$ and $Uz$, with $|Uz|=|z|$, gives
\begin{equation}\label{eq:ed-principal-complement}
 \sum_{n\notin\mathcal I_*}|c_{\delta,n}^{(0)}(\alpha,t)|^2
 \le \frac{2}{\delta^2}\sum_{n\notin\mathcal I_*}
 \bigl(|B_{n,0}(\alpha)z|^2+|B_{n,0}(\alpha)Uz|^2\bigr)
 \le C\frac{s^2|\alpha|^2}{\delta^2\lambda_*d^2}|z|^2 .
\end{equation}

Combining \eqref{eq:ed-principal-triple}--\eqref{eq:ed-principal-complement}, and integrating over
$|\alpha|<r$, we find
\begin{equation}
    \begin{aligned}
         &\int_{|\alpha|<r}\sum_{n\ge1} |c_{\delta,n}^{(0)}(\alpha,t)|^2\,d\alpha\\
 &\quad\le C\frac{s^2}{\delta^2\lambda_*d^2} \int_{|\alpha|<r}|\alpha|^2 \left| \widehat{\Env_0} \left( \frac{\alpha}{\delta}\right) \right|^2\,d\alpha + C\frac{s^4t^2}{\delta^2\lambda_*d^2}
       \int_{|\alpha|<r}|\alpha|^4 \left| \widehat{\Env_0} \left( \frac{\alpha}{\delta}\right) \right|^2\,d\alpha\\
 &\quad\le C\frac{s^2\delta^2}{\lambda_*d^2}\|\Env_0\|_{\dot H^1}^2 +C\frac{s^4\delta^4t^2}{\lambda_*d^2}
    \|\Env_0\|_{\dot H^2}^2.
    \end{aligned}
\end{equation}
Using $\lambda_*=\omega_*^2$, we conclude
\begin{equation}\label{eq:ed-principal-near-norm}
\left( \int_{|\alpha|<r}\sum_{n\ge1}|c_{\delta,n}^{(0)}(\alpha,t)|^2\,d\alpha \right)^{1/2}\le
 C\frac{s\delta}{\omega_*d}\|\Env_0\|_{\dot H^1} +C\frac{s^2\delta^2t}{\omega_*d}\|\Env_0\|_{\dot H^2}.
\end{equation}
Here, $|\alpha|<r$ is contained in $Y^*$ because $r\le1$.

\textit{Step 3: the principal term for $|\alpha|\ge r$.} By \eqref{eq:ed-single-k-global} with $k=0$,
\begin{equation}
    \begin{aligned}
    \int_{Y^*\cap\{|\alpha|\ge r\}}\sum_{n\ge1}|c_{\delta,n}^{(0)}(\alpha,t)|^2\,d\alpha &\le \frac{4}{\delta^2}
       \int_{Y^*\cap\{|\alpha|\ge r\}} \left| \widehat{\Env_0} \left( \frac{\alpha}{\delta}\right) \right|^2\,d\alpha\\
 &\le4\int_{|q|\ge \frac{r}{\delta}}|\widehat{\Env_0}(q)|^2\,dq \\
 & \leq 4 \left( \frac{\delta}{r} \right)^{2N} \int_{|q|\ge \frac{r}{\delta}} |q|^{2N}|\widehat{\Env_0}(q)|^2\,dq .
    \end{aligned}
\end{equation}
Hence,
\begin{equation}\label{eq:ed-principal-far-norm}
 \left( \int_{Y^*\cap\{|\alpha|\ge r\}}\sum_{n\ge1}|c_{\delta,n}^{(0)}(\alpha,t)|^2\,d\alpha \right)^{1/2}
 \le 2\left(\frac{\delta}{r}\right)^N\|\Env_0\|_{\dot H^N}.
\end{equation}

\textit{Step 4: conclusion.} The norm identity \eqref{eq:lemma31-error-norm} in Proposition \ref{lem:discrepancy-bloch} gives
\begin{equation}
    \|e_\delta(\cdot,t)\|_{L_b^2(\mathbb R^2)} = \left( \int_{Y^*}\sum_{n\ge1}
       |c_{\delta,n}(\alpha,t)|^2\,d\alpha \right)^{1/2}.
\end{equation}
Decomposing into the sum over $k\ne0$, the principal term on $|\alpha|<r$, and the principal term on $|\alpha|\ge r$, the triangle inequality and \eqref{eq:ed-nonzero-k-norm}, \eqref{eq:ed-principal-near-norm}, and \eqref{eq:ed-principal-far-norm} yield
\begin{equation}
    \begin{aligned}
        \|e_\delta(\cdot,t)\|_{L_b^2(\mathbb R^2)} &\le C\frac{s\delta}{\omega_*d}\|\Env_0\|_{\dot H^1}
       +C\frac{s^2\delta^2t}{\omega_*d}\|\Env_0\|_{\dot H^2}\\
 &\quad+ 2\left(\frac{\delta}{r}\right)^N\|\Env_0\|_{\dot H^N}+C_N\delta^N\|\Env_0\|_{\dot H^N}.
    \end{aligned}
\end{equation}
Since $r\le1$, the last two terms are bounded by
\begin{equation}
     C_N\left(\frac{\delta}{r}\right)^N\|\Env_0\|_{\dot H^N}.
\end{equation}
This proves \eqref{eq:effective-dynamics-error}.
\end{proof}

\subsection{Longitudinal and transverse envelopes}\label{sec:longitudinal}
Suppose that the medium has $D_4$ symmetry and $G=\gamma^2 I_2$ with $\gamma>0$. In the canonical basis of Proposition~\ref{prop:generalcriteria}, write the envelope as $\Env=(h,\mathrm i v_1,\mathrm i v_2)^\top$, with $v=(v_1,v_2)^\top$. With $v_D$ as in \eqref{eq:square-frequency-cones}, substitution of the normal form \eqref{normalform} in \eqref{effectiveequationofE} gives
\begin{equation}\label{eq:longitudinal-transverse}
 \partial_T h+v_D\operatorname{div}_X v=0,\qquad
 \partial_Tv+v_D\nabla_X h=0.
\end{equation}
Thus, $\partial_T^2 h-v_D^2\Delta_Xh=0$, while $\partial_T(\partial_{X_1}v_2-\partial_{X_2}v_1)=0$. For a nonzero Fourier variable $\xi$, the component of $\widehat v$ parallel to $\xi$ couples to $\widehat h$ and the perpendicular component is stationary. These are the two conical branches and the central branch of the linearized frequency problem. The latter conclusion concerns the effective system to first order in quasimomentum.

\subsection{Effective permittivity and permeability near the crossing}\label{subsec:double-zero}

We first relate \eqref{eq:longitudinal-transverse} to the planar Maxwell equations. We then use a forced Bloch problem to identify the two bulk parameters separately. Square symmetry makes the quadratic spatial response isotropic and forces the cell means of the doublet modes to vanish. The additional nonzero mean condition on the singlet determines whether the uniform source couples to the crossing.

Consider a plane envelope $e^{\ii(\xi\cdot X-\Omega T)}$ and denote its scalar and vector amplitudes by $h$ and $v$. Its physical quasimomentum and frequency are $\alpha=\delta\xi$ and $\omega=\omega_*+\delta\Omega$, respectively. With $\nu=\omega-\omega_*$, \eqref{eq:longitudinal-transverse} gives
\begin{equation}\label{eq:dz-plane}
 \alpha\cdot v=\frac{\nu}{v_D}h,\qquad \alpha h=\frac{\nu}{v_D}v.
\end{equation}
Take the polarization with axial electric field $E_z$ and in-plane magnetic field $H$. We use the time dependence $e^{-\ii\omega t}$ and units in which the vacuum permittivity and permeability are one. Let $q=RH$, where $R$ is the quarter-turn in \eqref{rep-prelim:D4-standard-matrices}. For axial permittivity $\varepsilon_{\mathrm{eff}}(\omega)$ and isotropic in-plane permeability $\mu_{\mathrm{eff}}(\omega)$, the source-free Maxwell equations read
\begin{equation}\label{eq:dz-maxwell}
 \alpha\cdot q=\omega\varepsilon_{\mathrm{eff}}(\omega)E_z,\qquad
 \alpha E_z=\omega\mu_{\mathrm{eff}}(\omega)q.
\end{equation}
For any $Z_*>0$, the choice $E_z=Z_*h$ and $q=v$ makes \eqref{eq:dz-plane} the linearization of \eqref{eq:dz-maxwell} at a common zero of the two parameters, with
\begin{equation}\label{eq:dz-slope-product}
 \omega_*\varepsilon_{\mathrm{eff}}'(\omega_*)=\frac{1}{v_D Z_*},\qquad
 \omega_*\mu_{\mathrm{eff}}'(\omega_*)=\frac{Z_*}{v_D}.
\end{equation}
Thus, the cone speed fixes the product of the two slopes. Their ratio fixes the limiting impedance
$Z_*=\lim_{\omega\to\omega_*}\sqrt{\mu_{\mathrm{eff}}(\omega)/\varepsilon_{\mathrm{eff}}(\omega)}$.
On the propagating branches the signed refractive index, chosen positive above the crossing, is $\nu/(\omega_*v_D)+O(\nu^2)$ for smooth constitutive laws. On the stationary branch, $h=0$ and $v\perp\alpha$ for $\alpha\ne0$, so the magnetic field in this representation is longitudinal. These relations give the double-zero effective material properties associated with the three-component envelope \cite{HuangLaiHangZhengChan2011,Fang2016}. In the pressure formulation of acoustics, the same two parameters are the effective compressibility (the reciprocal of the bulk modulus) and the effective mass density, respectively \cite{LiuHuangChan2012}.

Guzina, Meng, and Oudghiri-Idrissi developed a finite-frequency
homogenization theory for forced periodic wave equations near repeated
and closely spaced Bloch eigenvalues \cite{guzina2}. To fix the relative field normalization, we take the macroscopic electric field to be the unweighted cell mean of the periodic microscopic electric field amplitude. We compute this mean for a source with constant periodic amplitude \cite{Silveirinha2007}. Write $z=\omega^2$. For $z\notin\spec L(\alpha)$, let $U_{\alpha,z}\in H^1_\#(Y)$ solve
\begin{equation}\label{eq:dz-forced-cell}
 a_\alpha(U_{\alpha,z},w)-z\langle U_{\alpha,z},w\rangle_b
       =\int_Y\overline w\,dx,\qquad w\in H^1_\#(Y),
\end{equation}
where $a_\alpha$ is given by \eqref{defaalphauv}. The source is continuous on the form domain. Since $|Y|=1$, its response and inverse response are
\begin{equation}\label{eq:dz-response}
 \mathscr S(\alpha,z):=\int_Y U_{\alpha,z}\,dx,\qquad
 \mathscr K(\alpha,z):=\frac{1}{\mathscr S(\alpha,z)}
 \quad\text{when }\mathscr S(\alpha,z)\ne0 .
\end{equation}
For fixed $z$ such that \eqref{eq:dz-response} is regular at $\alpha=0$, the form resolvent is real analytic in $\alpha$ there. Both the source and the mean are invariant under $D_4$. Inversion removes the odd powers of $\alpha$, and the quarter-turn makes the quadratic term isotropic. Thus,
\begin{equation}\label{eq:dz-taylor}
 \mathscr K(\alpha,z)=\mathscr K_0(z)+\mathscr K_2(z)|\alpha|^2+O_z(|\alpha|^4).
\end{equation}
Comparison with the scalar Maxwell symbol $|\alpha|^2/\mu_{\mathrm{eff}}-\omega^2\varepsilon_{\mathrm{eff}}$ defines the bulk parameters, wherever the right-hand sides are finite, by
\begin{equation}\label{eq:dz-definition}
 \varepsilon_{\mathrm{eff}}(\omega):=-\mathscr K_0(\omega^2)/\omega^2,\qquad
 \mu_{\mathrm{eff}}(\omega):=1/\mathscr K_2(\omega^2).
\end{equation}
In particular, testing \eqref{eq:dz-forced-cell} at $\alpha=0$ with $w=1$ gives
$\varepsilon_{\mathrm{eff}}(\omega)\int_Y U_{0,z}\,dx=\int_Y bU_{0,z}\,dx$.

\begin{proposition}\label{prop:dz-bulk}
Suppose that Assumption~\ref{assump:threefold} holds, the medium has $D_4$ symmetry, and $G=\gamma^2 I_2$ with $\gamma>0$. Write $f_0,f_1,f_2$ for the real canonical basis in \eqref{ftildetransform}, normalized by \eqref{b-orthogonality}, with $f_0$ the singlet. If its cell mean $m_*:=\int_Y f_0\,dx$ is nonzero, both bulk parameters in \eqref{eq:dz-definition} extend analytically to $\omega_*$ and have simple zeros there. For a fixed medium, as $\nu=\omega-\omega_*\to0$,
\begin{equation}\label{eq:dz-parameters}
 \varepsilon_{\mathrm{eff}}(\omega)=\frac{2\nu}{\omega_*m_*^2}+O(\nu^2),\qquad
 \mu_{\mathrm{eff}}(\omega)=\frac{m_*^2\nu}{2\omega_*v_D^2}+O(\nu^2),\qquad
 Z_*=\frac{m_*^2}{2v_D}.
\end{equation}
\end{proposition}

\begin{proof}
Use the canonical basis to construct the frame $\mathcal W$ and the reduced matrix $M$ in \eqref{eq:generalcriteria:frame-definition} and \eqref{eq:generalcriteria:reduction-proof}. For a unit vector $n$, put
$B_n=-\ii C(n)$ in this basis. Then
$M(kn)=\lambda_*I_3+kB_n+O(k^2)$ and
$\mathbf e_0^* B_n^2\mathbf e_0=\gamma^2$ by \eqref{eq:d4criteria:normalform-proof}.
Set $c_j(k)=\int_Y\overline{\mathcal W(kn)\mathbf e_j}\,dx$ for $j=0,1,2$, and let $c(k)$ be the column with these entries. Integration over $Y$ is invariant under $D_4$, so it vanishes on the doublet. Thus, $c(0)=m_*\mathbf e_0$.

Splitting the form resolvent into its contribution from the isolated cluster and its regular part gives
\begin{equation}\label{eq:dz-reduced-response}
 \mathscr S(kn,z)=c(k)^*(M(kn)-zI_3)^{-1}c(k)+\mathscr S_{\mathrm{reg}}(k,z).
\end{equation}
Indeed, on a fixed contour enclosing the triple, the inverse of the form pencil in \eqref{eq:dz-forced-cell} is analytic in $k$ as an operator from the antidual of $H^1_\#(Y)$ to $H^1_\#(Y)$. Its Cauchy integral with weight $(\zeta-z)^{-1}$ gives the part regular in $z$ inside the contour. Hence, $\mathscr S_{\mathrm{reg}}$ and its first two $k$-derivatives are bounded near $(0,\lambda_*)$.

Write $s_0=\lambda_*-z$ and $\mathscr S(kn,z)=S_0(z)+k^2 S_2(z)+O_z(k^4)$. Expanding the inverse in \eqref{eq:dz-reduced-response} at $k=0$ gives
\[
 S_0(z)=\frac{m_*^2}{s_0}+O(1),\qquad
 S_2(z)=\frac{m_*^2\gamma^2}{s_0^3}+O(s_0^{-2}).
\]
The term $s_0^{-3}B_n^2$ is the only term of pole order three in the quadratic coefficient of the matrix inverse. The quadratic part of $M$ and the derivatives of $c$ contribute poles of order at most two. In particular, $S_0$ is nonzero in a punctured neighborhood of $\lambda_*$. Taking the reciprocal gives
\[
 \mathscr K_0=\frac{s_0}{m_*^2}+O(s_0^2),\qquad
 \mathscr K_2=-\frac{S_2}{S_0^2}=-\frac{\gamma^2}{m_*^2s_0}+O(1).
\]
These are Laurent expansions in $z$. Since $s_0=-2\omega_*\nu+O(\nu^2)$ and $\gamma=2\omega_*v_D$, \eqref{eq:dz-definition} yields \eqref{eq:dz-parameters}.
\end{proof}

With the cell mean as an electric field, the corresponding leading amplitudes in \eqref{eq:dz-maxwell} are $E_z=m_*h$ and $q=(2v_D/m_*)v$. The forced response also determines the approximation when quasimomentum and detuning have the same order. Set $\alpha=\delta\xi$ and $\omega=\omega_*+\delta\Omega$. In the canonical basis, if $\Omega\ne0$ and $\Omega^2\ne v_D^2|\xi|^2$, the limiting matrix
$-\ii C(\xi)-2\omega_*\Omega I_3$ is invertible. The $(0,0)$ entry of its inverse is $-\Omega/[2\omega_*(\Omega^2-v_D^2|\xi|^2)]$. Applying \eqref{eq:dz-reduced-response} along $\delta\xi$ gives
\begin{equation}\label{eq:dz-joint}
 \mathscr K\bigl(\delta\xi,(\omega_*+\delta\Omega)^2\bigr)
 =\frac{2\omega_*\delta}{m_*^2}
    \left(\frac{v_D^2|\xi|^2}{\Omega}-\Omega\right)+O(\delta^2),
\end{equation}
uniformly on compact subsets of this set of $(\xi,\Omega)$ for a fixed medium. The leading term agrees with the scalar Maxwell symbol obtained from \eqref{eq:dz-parameters}. On the characteristic branches the three-component system \eqref{eq:longitudinal-transverse} retains all three branches and specifies the polarization on each branch.

\paragraph{Quadratic term of the central band.}
Fix a unit vector $n=(n_1,n_2)$, and let $M$ denote the reduced matrix in the canonical basis. Write $M(kn)=\lambda_*I_3+kB_n+k^2M_2(n)+O_n(k^3)$. By \eqref{eq:d4criteria:normalform-proof}, the Hermitian matrix $B_n=-\ii C(n)$ has a simple zero eigenvalue with normalized eigenvector $w_c(n)=(0,-n_2,n_1)^\top$. The matrix $(M(kn)-\lambda_*I_3)/k$ extends analytically to $k=0$. Differentiating its simple eigenvalue issuing from zero gives the coefficient of $k^2$ in the central band. Complex conjugation takes the fiber at $kn$ to the fiber at $-kn$, so the ordered central eigenvalue is even along this line. Consequently,
\begin{equation}\label{eq:central-directional-quadratic}
 \lambda_c(kn)=\lambda_*+k^2 w_c(n)^*M_2(n)w_c(n)+O_n(k^4).
\end{equation}
This expansion is taken at a fixed medium and direction. Section~\ref{subsec:weakdynamics} evaluates the limiting axial coefficient in the weak-contrast family.

\section{The weak perturbation model}\label{sec:weakperturbation}

We apply the framework of Sections \ref{sec:bloch} and \ref{sec:effective} to weak perturbations of a constant medium. We first identify a fourfold eigenvalue at $\Gamma$ for the reference constant medium. $D_4$-invariant perturbations satisfying moment conditions then give exact tuning curves along which a singlet coincides with a doublet, producing an isolated threefold Dirac point at $\Gamma$. We then verify these moment conditions for concrete examples and establish the effective dynamics with explicit dependence on the perturbations.

\subsection{\texorpdfstring{Moment conditions for a threefold crossing at $\Gamma$}{Moment conditions for a threefold crossing at Gamma}} \label{subsec:weakconstantbase}

Let $a>0$ and take
\begin{equation}\label{eq:weakconstantbase}
 A_0=aI_2,\qquad b_0=1.
\end{equation}
At $\Gamma$, the Bloch eigenvalue problem is
\begin{equation}
     -a\Delta u=\lambda u,\qquad u\in H_\#^1(Y).
\end{equation}
The Fourier modes give the eigenvalues
\begin{equation}
    \lambda_k^{\mathrm{free}}= 4\pi^2a |k|^2, \qquad k\in \mathbb{Z}^2.
\end{equation}
Consequently,
\begin{equation} \label{eq:weakconstantlambda}
 \lambda_0:=4\pi^2a
\end{equation}
has exact multiplicity four, corresponding to $k= (\pm 1,0),  (0,\pm 1)$.  A real $L^2(Y)$-orthonormal basis of the corresponding eigenspace is
\begin{equation}\label{eq:weakconstantquartet}
 \begin{aligned}
 \phi_A&=\cos(2\pi x_1)+\cos(2\pi x_2), &&\phi_B=\cos(2\pi x_1)-\cos(2\pi x_2),\\
 \phi_{E_1}&=\sqrt{2}\sin(2\pi x_1), &&\phi_{E_2}=\sqrt{2}\sin(2\pi x_2).
 \end{aligned}
\end{equation}
To separate the $A_1$ and $B_1$ singlets from one axial component of the doublet, define
\begin{equation}\label{symmetryspaceAB1}
    \begin{aligned}
 \mathcal V_A&:=\{u\in H_\#^1(Y;\mathbb R):u(Rx)=u(x),\ u(Sx)=u(x),\ x\in Y\},\\
 \mathcal V_B&:=\{u\in H_\#^1(Y;\mathbb R):u(Rx)=-u(x),\ u(Sx)=u(x), \ x\in Y \},\\
 \mathcal V_{E_1}&:=\{u\in H_\#^1(Y;\mathbb R):
     u(-x_1,x_2)=-u(x),\ u(x_1,-x_2)=u(x),\ x=(x_1,x_2)\in Y\}.
 \end{aligned}
\end{equation}

The basis in \eqref{eq:weakconstantquartet} satisfies $\phi_A\in\mathcal V_A$, $\phi_B\in\mathcal V_B$, $\phi_{E_1}\in\mathcal V_{E_1}$ and $\phi_{E_2}(x)=\phi_{E_1}(R^\top x)$. The subscript $\#$ denotes periodicity, including in the coefficient space $L^\infty_\#(Y)$.

Let $p=(p_1,p_2)\in\mathbb R^2$. For $j=1,2$, let $\sigma_j,\beta_j\in L^\infty_\#(Y;\mathbb R)$ be $D_4$-invariant profiles for the stiffness and mass density, respectively. Define
\begin{equation}\label{eq:weakjointfamily}
  A_p:=a(1+ p_1 \sigma_1+ p_2 \sigma_2)I_2,\qquad b_p:= 1+ p_1 \beta_1+ p_2 \beta_2.
\end{equation}
The following lemma gives the four eigenvalue branches issuing from $\lambda_0$ and their first variations.

\begin{lemma}\label{lem:perturbationsplitweak}
    Set
\begin{equation}\label{defofepsp}
    \epsilon(p):=\sum_{j=1}^2|p_j|\bigl(\|\sigma_j\|_\infty+\|\beta_j\|_\infty\bigr).
\end{equation}
There exist universal constants $\epsilon_0,C>0$ and, on the open set $\{p\in\mathbb R^2:\epsilon(p)<\epsilon_0\}$, real analytic functions $\lambda_A(p),\lambda_B(p),\lambda_E(p)$ and real analytic
$H^1_\#(Y;\mathbb R)$-valued functions $u_A(p),u_B(p),u_{E_1}(p),u_{E_2}(p)$ with the following
properties:
\begin{itemize}
    \item[(a)] Set
\begin{equation}
    \lambda_{E_1}(p):=\lambda_E(p),\qquad\lambda_{E_2}(p):=\lambda_E(p).
\end{equation}
The four functions satisfy
\begin{equation}
    -\nabla \cdot (A_p\nabla u_M(p))=\lambda_M(p) b_p u_M(p), \qquad M\in \{A,B,E_1,E_2\}.
\end{equation}
Moreover,
\begin{equation}
    u_A(p)\in\mathcal V_A,\qquad u_B(p)\in\mathcal V_B,\qquad
u_{E_1}(p)\in\mathcal V_{E_1}, \qquad u_{E_2}(p)(x)=u_{E_1}(p)(R^\top x),
\end{equation}
and
\begin{equation}\label{normalizebpuMuN}
    \int_Y b_pu_M(p)u_N(p)\,dx=\delta_{MN}, \qquad M,N\in\{A,B,E_1,E_2\}.
\end{equation}
\item[(b)] For every $M\in \{A,B,E_1,E_2\}$,
\begin{equation}
    u_M(0)=\phi_M,\qquad \lambda_M(0)=\lambda_0,
\end{equation}
and
\begin{equation}
    \|u_M(p)-\phi_M\|_{H^1(Y)}\le C\epsilon(p),\qquad|\lambda_M(p)-\lambda_0|\le Ca\epsilon(p).
\end{equation}
The first-order expansion of the eigenvalue with respect to $p$ is given by
\begin{equation}\label{eq:weak-branch-expansion}
\lambda_M(p) =\lambda_0 +a\sum_{j=1}^2p_j
\left[ \int_Y\sigma_j|\nabla\phi_M|^2\,dx-4\pi^2\int_Y\beta_j\phi_M^2\,dx\right]
+O\bigl(a\epsilon(p)^2\bigr),
\end{equation}
where the remainder constant is universal.
\item[(c)] The eigenvalues of the perturbed eigenvalue problem
\begin{equation}\label{perturedpeigenprob}
    -\nabla \cdot (A_p\nabla u )=\lambda b_p u , \qquad u \in H^1_\#(Y;\mathbb{R})
\end{equation}
in $(\lambda_0/2,3\lambda_0/2)$, counted with multiplicity, are $\lambda_A(p),\lambda_B(p),\lambda_{E_1}(p),\lambda_{E_2}(p)$.
\end{itemize}
\end{lemma}
\begin{proof}
    All functions in this proof are real. For $u,v \in H^1_\#(Y)$, introduce the bilinear forms
\begin{equation}
\mathfrak{a}_p(u,v) =\int_Y  A_p \nabla u\cdot \nabla v\,dx, \qquad  \mathfrak{b}_p(u,v) =\int_Y b_p uv\,dx,
\end{equation}
and
\begin{equation}\label{defofqp}
    \mathfrak{q}_p :=\mathfrak{a}_p-\mathfrak{a}_0 - \lambda_0 (\mathfrak{b}_p-\mathfrak{b}_0).
\end{equation}

Let $\iota\in \{A,B,E_1\}$, and $\mathcal{V}_\iota^\perp$ be the subspace of $\mathcal{V}_\iota$ perpendicular to $\phi_\iota$ under the $L^2(Y)$ inner product. More precisely, 
\begin{equation}
    \mathcal{V}_\iota^\perp := \{w\in \mathcal{V}_\iota:\int_Y w\phi_\iota \,dx=0 \}.
\end{equation}
The symmetry projections commute with $-a\Delta$. In the $\iota$ sector, $\phi_\iota$ spans its eigenspace at $\lambda_0$. Choose an $L^2(Y)$-orthonormal eigenbasis $\{\psi_m\}_{m\ge1}$ of the $L^2(Y)$-closure of $\mathcal V_\iota^\perp$. Its members lie in $\mathcal V_\iota^\perp$ and satisfy
\begin{equation}\label{eigenfunctionsinV}
    -a\Delta \psi_m = \nu_m \psi_m, \qquad \nu_m = 4\pi^2 a|k_m|^2
\end{equation}
for some $k_m\in \mathbb{Z}^2$. We have
\begin{equation}\label{characH1forvVperp}
    \|v\|_{H^1(Y)}^2 = \sum_{m\geq 1} \left( 1+ \frac{\nu_m}{a} \right) |\langle v,\psi_m \rangle_{L^2(Y)}|^2, \qquad v\in \mathcal{V}_\iota^\perp.
\end{equation}

For a real linear functional $h:\mathcal{V}_\iota^\perp \rightarrow \mathbb{R}$ that satisfies
\begin{equation}\label{eq:dual-functional}
    |h(v)|\leq M \| v\|_{H^1(Y)}, \qquad v\in \mathcal{V}_\iota^\perp,
\end{equation}
we claim that 
\begin{equation}\label{estimateofw3Ma}
    \left\| \sum_{m\geq 1} \frac{h(\psi_m)}{\nu_m-\lambda_0} \psi_m \right\|_{H^1(Y)} \leq \frac{3M}{a}.
\end{equation}
Indeed, \eqref{characH1forvVperp} implies that
\begin{equation}\label{boundofwnumhpsim}
\begin{aligned}
    \left\| \sum_{m\geq 1} \frac{h(\psi_m)}{\nu_m-\lambda_0} \psi_m \right\|_{H^1(Y)}^2& =\sum_{m\geq 1} \left( 1+ \frac{\nu_m}{a} \right) \frac{|h(\psi_m)|^2}{|\nu_m-\lambda_0|^2} \\
    &=\sum_{m\geq 1} \frac{(1+\nu_m/a)^2}{|\nu_m-\lambda_0|^2}  \frac{|h(\psi_m)|^2}{1+\nu_m/a}.
\end{aligned}
\end{equation}
By the explicit values of $\nu_m$ and $\lambda_0$, we obtain
\begin{equation}\label{spectralgap3a}
    \frac{(1+\nu_m/a)^2}{|\nu_m-\lambda_0|^2}  \leq \frac{9}{a^2}, \qquad m\geq 1.
\end{equation}
For $N\ge1$, take $v_N=\sum_{m=1}^Nh(\psi_m)(1+\nu_m/a)^{-1}\psi_m$ in \eqref{eq:dual-functional}. Then
\begin{equation}\label{estimateofvhpsim}
 \sum_{m=1}^N\frac{|h(\psi_m)|^2}{1+\nu_m/a}
 =\|v_N\|_{H^1(Y)}^2=h(v_N)\le M\|v_N\|_{H^1(Y)} .
\end{equation}
Taking the supremum over $N$ bounds the full sum by $M^2$.
Thus, \eqref{estimateofw3Ma} follows from \eqref{boundofwnumhpsim}--\eqref{estimateofvhpsim}.

\textit{Step 1.} We first solve the perturbed eigenvalue problem in the space $\mathcal{V}^\perp_\iota$. Given $(p,\mu)\in \mathbb{R}^3$, we assume that there is a real function $w=w(p,\mu)\in \mathcal{V}_\iota^\perp$ such that
\begin{equation}
    \lambda_\iota:=\lambda_0+\mu
\end{equation}
and
\begin{equation}
    u_\iota := \phi_\iota + w
\end{equation}
solves the variational problem
\begin{equation}
    \mathfrak{a}_p( u_\iota , \psi_m)= \lambda_\iota \mathfrak{b}_p( u_\iota , \psi_m), \qquad m\geq 1.
\end{equation}
We may write $w=\sum_{m\geq 1}w_m \psi_m$. Using the definition \eqref{defofqp} of $\mathfrak{q}_p$,  \eqref{eigenfunctionsinV}, and orthogonality, we obtain
\begin{equation}\label{constructionofwm}
    (\nu_m-\lambda_0) w_m =-\mathfrak{q}_p(\phi_\iota+w,\psi_m)+\mu\mathfrak{b}_p(\phi_\iota+w,\psi_m).
\end{equation}
The right-hand side still depends on $w$, so we solve these coefficient equations by iteration. Motivated by \eqref{constructionofwm}, we define 
\begin{equation}\label{iterationdefinition}
\left\{
\begin{aligned}
    w^{(j+1)}&:= \sum_{m\geq 1} \frac{-\mathfrak{q}_p(\phi_\iota+w^{(j)},\psi_m)+\mu\mathfrak{b}_p(\phi_\iota+w^{(j)},\psi_m)}{\nu_m-\lambda_0} \psi_m, \qquad j\geq 0, \\
    w^{(0)}&:=0.
\end{aligned}\right.
\end{equation}
Recall that $\epsilon(p)$ is defined in \eqref{defofepsp}.
The Cauchy--Schwarz inequality gives, for $u,v\in H^1_\#(Y)$,
\begin{equation}\label{estimateofqb}
    \begin{aligned}
    |(\mathfrak b_p - \mathfrak{b}_0)(u,v)| &\leq \epsilon(p)\|u\|_{H^1(Y)} \|v\|_{H^1(Y)} ,\\
  |\mathfrak{q}_p(u,v)| &\leq 4\pi^2 a\epsilon(p) \|u\|_{H^1(Y)} \|v\|_{H^1(Y)} .
    \end{aligned}
\end{equation}
It follows from \eqref{estimateofw3Ma}, \eqref{estimateofqb}, and
\begin{equation}\label{formuleachiteration}
     w^{(j+1)} - w^{(j)} = \sum_{m\geq 1} \frac{-\mathfrak{q}_p(w^{(j)}-w^{(j-1)},\psi_m)+\mu\mathfrak{b}_p(w^{(j)}-w^{(j-1)},\psi_m)}{\nu_m-\lambda_0} \psi_m
\end{equation}
that
\begin{equation}\label{boundofeachiteration}
    \| w^{(j+1)} - w^{(j)} \|_{H^1(Y)} \leq \kappa(p,\mu) \| w^{(j)} - w^{(j-1)} \|_{H^1(Y)},
\end{equation}
where
\begin{equation}\label{boundofkappa}
    \kappa(p,\mu):= \frac{3}{a} \big[ 4\pi^2 a \epsilon (p) + |\mu| (1+\epsilon (p)) \big] .
\end{equation}

Choose $(p,\mu)\in \mathbb{R}^3$ such that
\begin{equation}\label{rangeofepsandmu}
    \epsilon(p)+\frac{|\mu|}{a}<\frac{1}{1000}.
\end{equation}
Then $\kappa(p,\mu)<\frac{1}{2}$, and the iterates converge to the unique  solution  $w=w(p,\mu)$ satisfying
\begin{equation}\label{eigenproblemforVperp}
    \mathfrak{a}_p(\phi_\iota+w, v)=(\lambda_0+\mu) \mathfrak{b}_p(\phi_\iota+w, v), \qquad v\in \mathcal{V}_\iota^\perp.
\end{equation}
For $v\in\mathcal V_\iota^\perp$ one has $\mathfrak b_0(\phi_\iota,v)=0$. Thus, the first iterate contains $\mu(\mathfrak b_p-\mathfrak b_0)(\phi_\iota,v)$ in its mass term. Since $\|\phi_\iota\|_{H^1(Y)}^2=1+4\pi^2$, the bound \eqref{estimateofw3Ma} applied to the first iterate, together with \eqref{estimateofqb}, \eqref{boundofeachiteration}, and \eqref{rangeofepsandmu}, gives
\begin{equation}\label{boundofwpmu}
    \| w(p,\mu) \|_{H^1(Y)}  \leq C\epsilon(p) \left(1+\frac{|\mu|}{a} \right) \le C_2\epsilon(p)
\end{equation}
with universal constants $C,C_2>0$. Each iterate $w^{(j)}$ is polynomial in the parameters $(p,\mu)$. The same estimate holds for sufficiently small complex values of these parameters, so the convergence is locally uniform there. Hence, the limit is analytic, and its restriction to real parameters is real analytic.

We estimate the dependence of $w(p,\mu)$ on $\mu$. For $v\in\mathcal V_\iota^\perp$, orthogonality, \eqref{estimateofqb}, and \eqref{boundofwpmu} imply
\begin{equation}\label{boundofbpv}
    \begin{aligned}
    &|\mathfrak b_p(\phi_\iota+w(p,\mu),v)| \\
    &\qquad = |\mathfrak b_0(w(p,\mu),v)+(\mathfrak b_p-\mathfrak b_0)(\phi_\iota+w(p,\mu),v) | \\
    & \qquad \leq C\epsilon(p)\|v\|_{H^1(Y)}.
\end{aligned}
\end{equation}
Let $\mu,\mu'$ satisfy both \eqref{rangeofepsandmu} for the same $p$. Subtracting equations \eqref{eigenproblemforVperp} for $\mu$ and $\mu'$ and applying \eqref{estimateofw3Ma} and \eqref{boundofbpv} give
\begin{equation}
    \begin{aligned}
    \|w(p,\mu)-w(p,\mu')\|_{H^1(Y)}&\le\kappa(p,\mu)\|w(p,\mu)-w(p,\mu')\|_{H^1(Y)}\\
    &\quad+\frac{C\epsilon(p)}{a}|\mu-\mu'|.
\end{aligned}
\end{equation}
The condition \eqref{rangeofepsandmu} ensures $\kappa(p,\mu)<\frac{1}{2}$. Absorbing the first term yields
\begin{equation}\label{eq:w-lipschitz-mu}
    \|w(p,\mu)-w(p,\mu')\|_{H^1(Y)}\le \frac{C\epsilon(p)}{a}|\mu-\mu'|.
\end{equation}

Next, we determine $\mu=\mu(p)$ by imposing the eigenvalue equation in the remaining direction $\phi_\iota$. Since $w\in\mathcal{V}_\iota^\perp$ and $\|\phi_\iota\|_{L^2(Y)}=1$, we have $\mathfrak b_0(\phi_\iota+w,\phi_\iota)=1$.
Moreover, using $-a\Delta\phi_\iota=\lambda_0\phi_\iota$, we obtain
\begin{equation}
    \mathfrak a_0(\phi_\iota+w,\phi_\iota)=\lambda_0\mathfrak b_0(\phi_\iota+w,\phi_\iota)=\lambda_0.
\end{equation}
Consequently, the equation
\begin{equation}\label{eigenproblemforphi}
    \mathfrak a_p(\phi_\iota+w,\phi_\iota)=(\lambda_0+\mu)\mathfrak b_p(\phi_\iota+w,\phi_\iota)
\end{equation}
is equivalent to finding the fixed point of
\begin{equation}\label{defofDp}
    \mathscr D_p(\mu):=\frac{\mathfrak q_p(\phi_\iota+w(p,\mu),\phi_\iota)}{\mathfrak b_p(\phi_\iota+w(p,\mu),\phi_\iota)}.
\end{equation}
For $(p,\mu)$ satisfying \eqref{rangeofepsandmu}, by \eqref{boundofwpmu}, there is a universal constant $C_D>0$ such that
\begin{equation}\label{eq:uniform-scalar-denominator}
\begin{aligned}
    |\mathfrak b_p(\phi_\iota+w(p,\mu),\phi_\iota)-1| &=|(\mathfrak b_p-\mathfrak b_0)(\phi_\iota+w(p,\mu),\phi_\iota)|\\
    &\le C_D\epsilon(p).
\end{aligned}
\end{equation}
Whenever
\begin{equation}
    C_D\epsilon(p)<\frac{1}{2},
\end{equation}
the denominator of \eqref{defofDp} is nonzero. Hence $\mathscr{D}_p(\mu)$ is well-defined. 

We now specify a range of $p$ on which the fixed point equation $\mu = \mathscr D_p(\mu)$ has a solution $\mu(p)$ compatible with \eqref{rangeofepsandmu}. The estimates \eqref{estimateofqb}, \eqref{boundofwpmu}, \eqref{eq:w-lipschitz-mu}, and \eqref{eq:uniform-scalar-denominator} give universal constants $C_0,C_1\ge1$ such that
\begin{equation}\label{eq:scalar-map-estimates}
\begin{aligned}
    |\mathscr D_p(\mu)|&\le C_0a\epsilon(p),\\
    |\mathscr D_p(\mu)-\mathscr D_p(\mu')|&\le C_1\epsilon(p)^2|\mu-\mu'|.
\end{aligned}
\end{equation}
Choose
\begin{equation}
    K:=\max\{1,2C_0\},\qquad
    \epsilon_0:=\min\left\{\frac{1}{2000(1+K)},\frac{1}{4C_D},\frac{1}{2\sqrt{C_1}}\right\},
\end{equation}
and define the parameter neighborhood
\begin{equation}\label{eq:admissible-p-range}
    \mathcal P_{\epsilon_0}:=\{p\in\mathbb R^2:\epsilon(p)<\epsilon_0\}.
\end{equation}
Fix any $p\in\mathcal P_{\epsilon_0}$, and consider the interval
\begin{equation}
    I_p:=[-Ka\epsilon(p),Ka\epsilon(p)].
\end{equation}
For every $\mu\in I_p$,
\begin{equation}
    \epsilon(p)+\frac{|\mu|}{a}\le(1+K)\epsilon(p)<\frac1{1000}.
\end{equation}
Thus, $w(p,\mu)$ is defined for $p\in \mathcal{P}_{\epsilon_0}$ and $\mu\in I_p$. Moreover, \eqref{eq:scalar-map-estimates} shows that $\mathscr D_p$ maps $I_p$ into itself and has a Lipschitz constant at most $1/4$. The contraction mapping theorem therefore gives a unique fixed point $\mu=\mu(p)\in I_p$. In particular,
\begin{equation}\label{eq:quantitative-mu-range}
    |\mu(p)|\le Ka\epsilon(p), \qquad \epsilon(p)+\frac{|\mu(p)|}{a}<\frac1{1000}, \qquad p\in\mathcal P_{\epsilon_0}.
\end{equation}

At the fixed point, $|\partial_\mu(\mu-\mathscr D_p(\mu))|\ge3/4$. The analytic implicit function theorem therefore gives a real analytic function $\mu(p)$ on $\mathcal P_{\epsilon_0}$. Set
\begin{equation}
    w(p):=w(p,\mu(p)),\qquad\lambda_\iota(p):=\lambda_0+\mu(p), \qquad p\in \mathcal P_{\epsilon_0}.
\end{equation}
These functions satisfy
\begin{equation}
    \mathfrak a_p(\phi_\iota+w(p),v)=\lambda_\iota(p)\mathfrak b_p(\phi_\iota+w(p),v),\qquad v\in\mathcal V_\iota,
\end{equation}
and
\begin{equation}\label{boundofwandmu}
    \|w(p)\|_{H^1(Y)}\le C\epsilon(p),\qquad|\mu(p)|\le Ca\epsilon(p).
\end{equation}

We identify the first-order expansion of $\mu(p)$. Rewrite $\mu = \mathscr{D}_p(\mu(p))$ as
\begin{equation}
\mu(p)-\mathfrak q_p(\phi_\iota,\phi_\iota)
=\mathfrak q_p(w(p),\phi_\iota)-\mu(p)(\mathfrak b_p-\mathfrak b_0)(\phi_\iota+w(p),\phi_\iota).
\end{equation}
The estimates \eqref{estimateofqb} and \eqref{boundofwandmu} yield
\begin{equation}
    |\mathfrak q_p(w(p),\phi_\iota)|+|\mu(p)(\mathfrak b_p-\mathfrak b_0)(\phi_\iota+w(p),\phi_\iota)|\le Ca\epsilon(p)^2.
\end{equation}
Therefore, the first-order expansion of $\lambda_\iota(p)$ is given by
\begin{equation}\label{eq:weak-sector-expansion}
\lambda_\iota(p)=\lambda_0+\mathfrak q_p(\phi_\iota,\phi_\iota)+O(a\epsilon(p)^2),
\end{equation}
where
\begin{equation}\label{eq:weak-first-order-shift}
\mathfrak q_p(\phi_\iota,\phi_\iota)=a\sum_{j=1}^2p_j 
\left[\int_Y\sigma_j|\nabla\phi_\iota|^2\,dx-4\pi^2\int_Y\beta_j\phi_\iota^2\,dx\right].
\end{equation}

Finally, define
\begin{equation}
    u_\iota(p):=\frac{\phi_\iota + w(p)}{\sqrt{\mathfrak b_p(\phi_\iota+w(p),\phi_\iota+w(p))}} .
\end{equation}
The denominator in the definition of $u_\iota(p)$ is positive for a sufficiently small $\epsilon$. Indeed, orthogonality gives
\begin{equation}
\begin{aligned}
\mathfrak b_p(\phi_\iota+w(p),\phi_\iota+w(p))={}&1+\|w(p)\|_{L^2(Y)}^2\\
&+(\mathfrak b_p-\mathfrak b_0)(\phi_\iota+w(p),\phi_\iota+w(p))=1+O(\epsilon(p)).
\end{aligned}
\end{equation}
It follows that $u_\iota(p)$ is real analytic and satisfies
\begin{equation}
    \mathfrak b_p(u_\iota(p),u_\iota(p))=1,\qquad\|u_\iota(p)-\phi_\iota\|_{H^1(Y)}\le C\epsilon(p).
\end{equation}
Combining \eqref{eigenproblemforVperp} and \eqref{eigenproblemforphi}, we get
\begin{equation}\label{eq:weak-sector-eigenpair}
\mathfrak a_p(u_\iota(p),v)=\lambda_\iota(p)\mathfrak b_p(u_\iota(p),v),\qquad v\in\mathcal V_\iota.
\end{equation}

\textit{Step 2.} The variational identity \eqref{eq:weak-sector-eigenpair} is only for test functions in the corresponding symmetry space $\mathcal{V}_\iota$. We first extend this equation to every $v\in H^1_\#(Y)$. Let $\chi_B=\chi_{-1,1}$, so $\chi_B(R)=-1$ and $\chi_B(S)=1$. For arbitrary $v\in H^1_\#(Y)$, define the symmetrizations
\begin{equation}\label{symmetrizationAB}
    v_A(x):= \frac{1}{8}\sum_{Q\in D_4} v(Q^\top x), \qquad v_B(x):= \frac{1}{8}\sum_{Q\in D_4} \chi_B (Q) v(Q^\top x),
\end{equation}
and
\begin{equation}\label{symmetrizationE1}
    v_{E_1}(x_1,x_2):=\frac{1}{4} \big[ v(x_1,x_2) - v(-x_1,x_2) + v(x_1,-x_2) - v(-x_1,-x_2) \big] .
\end{equation}
These symmetrizations satisfy $v_\iota \in \mathcal{V}_\iota$ for $\iota\in \{A,B,E_1\}$. Hence,
$D_4$-invariance of $\sigma_j,\beta_j$ yields
\begin{equation}
    \mathfrak{a}_p(u_\iota,v) = \mathfrak{a}_p(u_\iota,v_\iota), \qquad  \mathfrak{b}_p(u_\iota,v) = \mathfrak{b}_p(u_\iota,v_\iota).
\end{equation}
Therefore, \eqref{eq:weak-sector-eigenpair} holds for every $v\in H^1_\#(Y)$. Next, we define
\begin{equation}\label{eq:weak-second-doublet-mode}
   \lambda_E(p):=\lambda_{E_2}(p):=\lambda_{E_1}(p), \qquad  u_{E_2}(p)(x):=u_{E_1}(p)(R^\top x).
\end{equation}
By rotational invariance, 
\begin{equation}
    \begin{aligned}
    \mathfrak a_p(u_{E_2}(p),v)&=\mathfrak a_p(u_{E_1}(p),v(R\cdot))\\
    &=\lambda_E(p)\mathfrak b_p(u_{E_1}(p),v(R\cdot))=\lambda_E(p)\mathfrak b_p(u_{E_2}(p),v),  \qquad v\in H^1_\#(Y),
\end{aligned}
\end{equation}
and
\begin{equation}
    \mathfrak b_p(u_{E_2}(p),u_{E_2}(p))=\mathfrak b_p(u_{E_1}(p),u_{E_1}(p))=1.
\end{equation}
Thus, $u_{E_1}(p)$ and $u_{E_2}(p)$ are linearly independent normalized eigenfunctions corresponding to $\lambda_E(p)$. Since $\phi_{E_1}(R^\top x)=\phi_{E_2}(x)$, we also obtain
\begin{equation}\label{eq:weak-doublet-mode-estimate}
    \|u_{E_2}(p)-\phi_{E_2}\|_{H^1(Y)}=\|u_{E_1}(p)-\phi_{E_1}\|_{H^1(Y)}\le C\epsilon(p).
\end{equation}

The four functions $u_A(p),u_B(p),u_{E_1}(p),u_{E_2}(p)$ are mutually orthogonal with respect to $\mathfrak b_p$. Indeed, $u_A(p)$ and $u_B(p)$ are even in both coordinates. Since $b_p$ is invariant under coordinate reflections, the parity of $u_{E_1}(p)$ and $u_{E_2}(p)$ gives
\begin{equation}
\begin{aligned}
    \mathfrak b_p(u_A(p),u_{E_1}(p))&=\mathfrak b_p(u_A(p),u_{E_2}(p))=\mathfrak b_p(u_B(p),u_{E_1}(p)) \\
    &=\mathfrak b_p(u_B(p),u_{E_2}(p))=\mathfrak b_p(u_{E_1}(p),u_{E_2}(p))=0.
\end{aligned}
\end{equation}
For the remaining pair, rotational symmetry yields
\begin{equation}
    \mathfrak b_p(u_A(p),u_B(p))=\mathfrak b_p\bigl(u_A(p)(R^\top\cdot),u_B(p)(R^\top\cdot)\bigr)=-\mathfrak b_p(u_A(p),u_B(p)),
\end{equation}
and therefore, $\mathfrak b_p(u_A(p),u_B(p))=0$. Together with normalization, this proves \eqref{normalizebpuMuN}.

\textit{Step 3.} Now, we prove part (c). Let $\lambda_n(p)$ denote the eigenvalues of problem \eqref{perturedpeigenprob}, ordered increasingly and counted with multiplicity. The pointwise coefficient bounds from \eqref{eq:weakjointfamily} and \eqref{defofepsp} give, for $v\in H^1_\#(Y;\mathbb{R})$,
\begin{equation}
    \begin{aligned}
    (1-\epsilon(p))\mathfrak a_0(v,v)&\le\mathfrak a_p(v,v)\le(1+\epsilon(p))\mathfrak a_0(v,v),\\
    (1-\epsilon(p))\mathfrak b_0(v,v)&\le\mathfrak b_p(v,v)\le(1+\epsilon(p))\mathfrak b_0(v,v).
\end{aligned}
\end{equation}
The min--max principle therefore gives
\begin{equation}\label{eq:weak-minmax-comparison}
    \frac{1-\epsilon(p)}{1+\epsilon(p)}\lambda_n(0)\le\lambda_n(p)\le\frac{1+\epsilon(p)}{1-\epsilon(p)}\lambda_n(0).
\end{equation}
For the constant medium with $A_0=aI_2$ and $b=1$,
\begin{equation}
    \lambda_1(0)=0,\qquad\lambda_2(0)=\lambda_3(0)=\lambda_4(0)=\lambda_5(0)=\lambda_0,\qquad\lambda_6(0)=2\lambda_0.
\end{equation}
Since $\epsilon(p)<\frac{1}{1000}$, \eqref{eq:weak-minmax-comparison} implies that $\lambda_2(p),\lambda_3(p),\lambda_4(p),\lambda_5(p)$ belong to $(\lambda_0/2,3\lambda_0/2)$,
whereas $\lambda_1(p)=0$ and $\lambda_n(p)>3\lambda_0/2$ for $n\ge6$. Thus, $(\lambda_0/2,3\lambda_0/2)$ contains exactly four eigenvalues, counting multiplicity.

Finally, \eqref{eq:quantitative-mu-range} gives
\begin{equation}
    |\lambda_\iota(p) -\lambda_0|\leq Ka\epsilon(p) < \frac{K}{2000(1+K)}a < \frac{1}{2}\lambda_0. 
\end{equation}
Thus, $\lambda_A(p),\lambda_B(p),\lambda_E(p)$ all belong to $(\lambda_0/2,3\lambda_0/2)$. This completes the proof of part (c).
\end{proof}

We define some quantities related to the perturbation. For $j=1,2$, 
\begin{equation}\label{eq:weakcombinedmoments}
\begin{aligned}
     \mathscr{M}_j&:=\int_Y (\beta_j - \sigma_j) \,dx, \\
     \mathscr{U}_j&:= \int_Y (\beta_j + \sigma_j) \cos(4\pi x_1) \,dx = \int_Y (\beta_j + \sigma_j) \cos(4\pi x_2) \,dx, \\
     \mathscr{Z}_j&:= 2\int_Y \beta_j\cos(2\pi x_1)\cos(2\pi x_2)\,dx ,
\end{aligned}
\end{equation}
and
\begin{equation}\label{eq:weakjointweights}
 m_{A}^j:=\mathscr{M}_j+\mathscr{U}_j+\mathscr{Z}_j,\qquad
 m_{B}^j:=\mathscr{M}_j+\mathscr{U}_j-\mathscr{Z}_j,\qquad
 m_E^j=m_{E_1}^j=m_{E_2}^j:=\mathscr{M}_j-\mathscr{U}_j.
\end{equation}

By a direct computation, one has, for $M\in \{A,B,E_1,E_2\}$, 
\begin{equation}
    \int_Y\sigma_j|\nabla\phi_M|^2\,dx-4\pi^2\int_Y\beta_j\phi_M^2\,dx = -4\pi^2 m_{M}^j.
\end{equation}
Thus, by \eqref{eq:weak-branch-expansion}, we have, for $p\in \mathcal{P}_{\epsilon_0}$,
\begin{equation}\label{estimateoflambdaMp}
     \lambda_M(p)=\lambda_0-\lambda_0\sum_{j=1}^2p_jm_{M}^j+O(a\epsilon(p)^2).
\end{equation}
In particular,
\begin{equation}\label{eq:weak-singlet-doublet-differences}
\begin{aligned}
    \lambda_A(p)-\lambda_E(p)&=-\lambda_0\bigl[p_1(2\mathscr U_1+\mathscr Z_1)+p_2(2\mathscr U_2+\mathscr Z_2)\bigr]+O(a\epsilon(p)^2),\\
    \lambda_B(p)-\lambda_E(p)&=-\lambda_0\bigl[p_1(2\mathscr U_1-\mathscr Z_1)+p_2(2\mathscr U_2-\mathscr Z_2)\bigr]+O(a\epsilon(p)^2).
\end{aligned}
\end{equation}
Here, the remainder constant is universal.

The following theorem gives moment conditions under which a suitable tuning of the perturbation parameters produces accidental threefold Dirac points at $\Gamma$.

\begin{theorem}\label{thm:weak-dirac}
Assume that
\begin{equation}\label{momentcondition}
    \mathfrak d:=\mathscr U_1\mathscr Z_2-\mathscr U_2\mathscr Z_1\ne0,
    \qquad2\mathscr U_2+\mathscr Z_2\ne0,\qquad2\mathscr U_2-\mathscr Z_2\ne0.
\end{equation}
We call \eqref{momentcondition} the moment condition. Set
\begin{equation}
    \mathcal{M}_\infty:=\sum_{j=1}^2\bigl(\|\sigma_j\|_\infty+\|\beta_j\|_\infty\bigr),
\qquad\eta_{\rm m}:=\min\bigl\{ |\mathfrak d|, |2\mathscr U_2+\mathscr Z_2|, |2\mathscr U_2-\mathscr Z_2|\bigr\}>0,
\end{equation}
and define
\begin{equation}
    \mathfrak g_A:=\frac{4|\mathfrak d|}{|2\mathscr U_2+\mathscr Z_2|},\qquad 
    \mathfrak g_B:=\frac{4|\mathfrak d|}{|2\mathscr U_2-\mathscr Z_2|}.
\end{equation}
There exist constants $p_0,K>0$, depending only on $\mathcal{M}_\infty$ and $\eta_{\rm m}$, and a universal constant $C>0$, with the following properties. 
\begin{itemize}
    \item[(a)] There exist real analytic functions $p_{2}^A,p_{2}^B:(-p_0,p_0)\to\mathbb R$, independent of $a$, such that $p_{2}^A(0)=p_{2}^B(0)=0$, and
\begin{equation}\label{derivativeofp2X}
    (p_{2}^A)'(0)=-\frac{2\mathscr U_1+\mathscr Z_1}{2\mathscr U_2+\mathscr Z_2},
    \qquad (p_{2}^B)'(0)=-\frac{2\mathscr U_1-\mathscr Z_1}{2\mathscr U_2-\mathscr Z_2}.
\end{equation}
Write
\begin{equation}
    p^\iota(p_1):=(p_1,p_2^\iota(p_1)),\qquad \epsilon_\iota(p_1):=\epsilon(p^{\iota}(p_1)), \qquad \iota\in\{A,B\}.
\end{equation}
For each $\iota\in\{A,B\}$ and $0<|p_1|<p_0$, the perturbed medium associated with parameters $p=p^\iota(p_1)$ has an isolated eigenvalue
\begin{equation}
    \lambda_*^\iota(p_1):=\lambda_\iota(p^\iota(p_1)) =\lambda_E(p^\iota(p_1))
\end{equation}
of multiplicity three at $\Gamma$, satisfying
\begin{equation}\label{eq:weak-tuned-eigenvalue-bound}
    |\lambda_*^\iota(p_1)-\lambda_0|\le Ca\epsilon_\iota(p_1).
\end{equation}
    \item[(b)] For each $\iota\in\{A,B\}$ and $0<|p_1|<p_0$, the three-dimensional eigenspace of $\lambda_{*}^\iota(p_1)$
    admits a real basis $\{f_0,f_1,f_2\}\subset H^1_\#(Y)$ satisfying
\begin{equation}\label{normalweakperb}
    \int_Y b_{p^\iota(p_1)} f_rf_s\,dx=\delta_{rs},\qquad r,s\in\{0,1,2\},
\end{equation}
for which the coupling matrices take the form
\begin{equation}
    C_1=\gamma_\iota(p_1)\begin{pmatrix}
        0&1&0\\
        -1&0&0\\
        0&0&0
    \end{pmatrix},\qquad C_2=\gamma_\iota(p_1)
    \begin{pmatrix}
        0&0&1\\
        0&0&0\\
        -1&0&0
    \end{pmatrix}.
\end{equation}
Here, $\gamma_\iota$ satisfies
\begin{equation}\label{gammaXp1}
    \bigl|\gamma_\iota(p_1)-2\sqrt2\pi a\bigr|\le Ca\epsilon_\iota(p_1),\qquad \gamma_\iota(p_1)>0.
\end{equation}
\item[(c)] Let $\sigma_{p_1}^\iota$ denote the spectrum of the perturbed eigenvalue problem \eqref{perturedpeigenprob} with $p=p^\iota(p_1)$. The gap to the remaining eigenvalues at $\Gamma$ satisfies
\begin{equation}\label{twosidegapestimate}
    \left|\frac{ \operatorname{dist}\!\left(  \lambda_*^\iota(p_1),  \sigma_{p_1}^\iota\setminus\{\lambda_*^\iota(p_1)\}
 \right)}{\lambda_0\mathfrak{g}_\iota|p_1|}-1\right|\le K|p_1|.
\end{equation}
In particular, after decreasing $p_0$ if necessary,
\begin{equation}\label{twosidegapestimate2}
    \frac12\lambda_0\mathfrak{g}_\iota|p_1|\le\operatorname{dist}\!\left( \lambda_*^\iota(p_1),
 \sigma_{p_1}^\iota\setminus\{\lambda_*^\iota(p_1)\}\right)\le\frac32\lambda_0\mathfrak{g}_\iota|p_1|.
\end{equation}
\end{itemize}
\end{theorem}

\begin{remark}
The perturbation profiles $\sigma_j,\beta_j$ may be arbitrary bounded functions subject to their symmetries and \eqref{momentcondition}; symmetry alone does not imply this moment condition. For example, constant profiles only rescale the constant medium and the eigenvalue remains fourfold. The Jacobian determinant of the gap map $p\mapsto ((\lambda_A-\lambda_E)/\lambda_0,(\lambda_B-\lambda_E)/\lambda_0)$ at $p=0$ is $-4\mathfrak d$, so the two curves meet transversely. The additional conditions $2\mathscr U_2\pm\mathscr Z_2\ne0$ allow the path $\iota$ to be
expressed as the graph $p_{2}^\iota(p_1)$.

If the density profiles vanish, then $b_p\equiv1$ and $\mathscr Z_1=\mathscr Z_2=0$, and hence $\mathfrak d=0$. Thus, the moment condition \eqref{momentcondition} does not apply.
\end{remark}

\begin{proof}[Proof of Theorem~\ref{thm:weak-dirac}]

\textit{Step 1. Construction of the tuning curves.} For $\iota\in\{A,B\}$, define the scalar function
\begin{equation}
    h_\iota(p):=\frac{\lambda_\iota(p)-\lambda_E(p)}{\lambda_0}.
\end{equation}
By \eqref{eq:weak-singlet-doublet-differences}, $h_\iota(0)=0$, and, for $j=1,2$,
\begin{equation}
    \partial_{p_j}h_A(0)=-(2\mathscr U_j+\mathscr Z_j),\qquad
    \partial_{p_j}h_B(0)=-(2\mathscr U_j-\mathscr Z_j).
\end{equation}
In particular, 
\begin{equation}\label{upperboundofdpjhX}
    |\partial_{p_j}h_\iota(0)| \leq 4\mathcal{M}_\infty,
\end{equation}
and the moment condition \eqref{momentcondition} gives
\begin{equation}\label{loweeboundetapfdp2h}
    |\partial_{p_2}h_\iota(0)|\ge\eta_{\rm m}.
\end{equation}

Dividing the eigenvalue equation by $a$ shows that each ratio $\lambda_\iota(p)/\lambda_0$, and therefore each tuning curve, is independent of $a$.

We first explain why the tuning interval can be chosen using only $\mathcal M_\infty$ and $\eta_{\rm m}$. The construction in Lemma \ref{lem:perturbationsplitweak} extends to sufficiently small complex parameters, with uniform convergence of its iterations. The coefficient bounds used there depend on the profiles only through $\epsilon(p)$. Consequently, there is a complex parameter neighborhood, whose size depends only on $\mathcal M_\infty$, on which
\begin{equation}
    \frac{\lambda_A(p)}{\lambda_0},\qquad\frac{\lambda_B(p)}{\lambda_0},\qquad
\frac{\lambda_E(p)}{\lambda_0}
\end{equation}
are analytic and uniformly bounded. Cauchy's estimates on a smaller neighborhood therefore give constants $\rho,H>0$, depending only on $\mathcal M_\infty$, such that, for real parameters with $|p_1|+|p_2|<\rho$,
\begin{equation}\label{taylorestimateH}
    \left|h_\iota(p)-\partial_{p_1}h_\iota(0)p_1-\partial_{p_2}h_\iota(0)p_2\right|\le H(|p_1|+|p_2|)^2,
\end{equation}
and
\begin{equation}\label{lipschitzHofhX}
    |\partial_{p_2}h_\iota(p)-\partial_{p_2}h_\iota(0)| \le H(|p_1|+|p_2|).
\end{equation}
We may assume $H\ge1$.

Set
\begin{equation}
    L:=2\left(1+\frac{4\mathcal M_\infty}{\eta_{\rm m}}\right).
\end{equation}
For fixed $p_1$, consider the map
\[
T_{\iota,p_1}(p_2):=p_2-\frac{h_\iota(p_1,p_2)}{\partial_{p_2}h_\iota(0)}
\]
on the interval $I_{p_1}:=[-L|p_1|,L|p_1|]$. Choose $p_0>0$ so that
\begin{equation}
    p_0\le
    \min\left\{ \frac{\rho}{2(1+L)}, \frac{\eta_{\rm m}}{4H(1+L)}, \frac{\epsilon_0}{2(1+\mathcal M_\infty L)}\right\},
\end{equation}
where $\epsilon_0$ is the universal constant in Lemma \ref{lem:perturbationsplitweak}.
For $|p_1|<p_0$ and $p_2\in I_{p_1}$, we then have
\begin{equation}\label{estimateofdp2Tp1}
\begin{aligned}
    |\partial_{p_2}T_{\iota,p_1}(p_2)|& =\frac{|\partial_{p_2} h_\iota(0,0)-\partial_{p_2}h_\iota(p_1,p_2)|}{|\partial_{p_2} h_\iota(0,0)|}\\
    &  \le\frac{H(1+L)p_0}{\eta_{\rm m}} \\
    &\le\frac14,
\end{aligned}
\end{equation}
where we used \eqref{loweeboundetapfdp2h} and \eqref{lipschitzHofhX} in the first inequality, and $p_0 \leq \frac{\eta_{\rm m}}{4H(1+L)}$ in the second inequality. The restriction $p_0\leq\rho/[2(1+L)]$ ensures that $|p_1|+|p_2|<\rho$ in \eqref{lipschitzHofhX}. Moreover, 
\begin{equation}\label{estiofTp10}
\begin{aligned}
    |T_{\iota,p_1}(0)|&=\frac{|h_\iota(p_1,0)|}{|\partial_{p_2}h_\iota(0)|} \leq \frac{|\partial_{p_1}h_\iota(0) p_1 | + H|p_1|^2}{\eta_{\rm m}}\\
    &\le\left( \frac{4\mathcal M_\infty}{\eta_{\rm m}} +\frac{Hp_0}{\eta_{\rm m}}\right)|p_1| \\
    &\le\frac L2|p_1|,
\end{aligned}
\end{equation}
where we used \eqref{taylorestimateH} in the first inequality, and \eqref{upperboundofdpjhX} and $|p_1|<p_0$ in the second inequality. It follows from \eqref{estimateofdp2Tp1} and \eqref{estiofTp10} that
\begin{equation}
    |T_{\iota,p_1}(p_2)|\le\frac L2|p_1|+\frac14|p_2|\le\frac{3L}{4}|p_1|.
\end{equation}
Thus, $T_{\iota,p_1}$ maps $I_{p_1}$ into itself and is a contraction. Its unique fixed point, denoted $p_2^\iota(p_1)$, satisfies
\begin{equation}
    \lambda_\iota(p_1,p_2^\iota(p_1))=\lambda_E(p_1,p_2^\iota(p_1)),\qquad|p_2^\iota(p_1)|\le L|p_1|.
\end{equation}
On the graph $(p_1,p_2^\iota(p_1))$, $|p_1|<p_0$, 
\begin{equation}
    |\partial_{p_2}h_\iota(p_1,p_2^\iota(p_1))|\ge|\partial_{p_2}h_\iota(0)|   -H(|p_1|+|p_2^\iota(p_1)|)\ge\frac34\eta_{\rm m}.
\end{equation}
The analytic implicit function theorem, applied at each point of the graph, and uniqueness of the fixed point
therefore show that $p_2^\iota$ is real analytic on $(-p_0,p_0)$. Differentiating $h_\iota(p_1,p_2^\iota(p_1))=0$ at $p_1=0$ gives \eqref{derivativeofp2X}.

The estimate \eqref{taylorestimateH} for $p_2=p_2^\iota(p_1)$ gives
\begin{equation}\label{eq:weak-tuning-remainder}
\begin{aligned}
    \left|p_2^\iota(p_1)-(p_2^\iota)'(0)p_1\right|& = \frac{|\partial_{p_1} h_\iota(0) p_1 + \partial_{p_2}h_\iota(0) p_2^\iota(p_1)|}{|\partial_{p_2}h_\iota(0)|}\\
    &\le\frac{H(1+L)^2}{\eta_{\rm m}}\,p_1^2.
\end{aligned}
\end{equation}
In particular, this remainder bound depends only on
$\mathcal M_\infty$ and $\eta_{\rm m}$. Moreover, the restriction $p_0\leq \frac{\epsilon_0}{2(1+\mathcal M_\infty L)}$ implies that
\begin{equation}\label{boundofepsXp1}
    \epsilon_\iota(p_1)=\epsilon(p^\iota(p_1))\le\mathcal M_\infty L|p_1|<\frac{\epsilon_0}{2}.
\end{equation}
Therefore, Lemma \ref{lem:perturbationsplitweak} applies. Define
\begin{equation}\label{eq:weak-tuned-eigenvalues}
\lambda_*^\iota(p_1):=\lambda_\iota(p^\iota(p_1)) =\lambda_E(p^\iota(p_1)).
\end{equation}
Then, $\lambda_*^\iota$ is real analytic in $(-p_0,p_0)$, with $\lambda_*^\iota(0)=\lambda_0$, and Lemma \ref{lem:perturbationsplitweak} gives
\begin{equation}
|\lambda_*^\iota(p_1)-\lambda_0|\le Ca\epsilon_\iota(p_1).
\end{equation}

\textit{Step 2. Separation of the remaining singlet and
the spectral gap.} The expansions \eqref{eq:weak-singlet-doublet-differences} imply
\begin{equation}\label{differencelambdaBA}
    \lambda_B(p)-\lambda_A(p)=2\lambda_0(p_1\mathscr Z_1+p_2\mathscr Z_2) +O(a\epsilon(p)^2),
\end{equation}
where the remainder constant is universal.
Using \eqref{derivativeofp2X}, we obtain
\begin{equation}
    \mathscr Z_1+\mathscr Z_2(p_2^A)'(0)=-\frac{2\mathfrak d}{2\mathscr U_2+\mathscr Z_2},
\end{equation}
and
\begin{equation}\label{Z1Z2p2B}
    \mathscr Z_1+\mathscr Z_2(p_2^B)'(0)=-\frac{2\mathfrak d}{2\mathscr U_2-\mathscr Z_2}.
\end{equation}
Together with \eqref{eq:weak-tuning-remainder} and $\epsilon_\iota(p_1)\le\mathcal M_\infty L|p_1|$, \eqref{differencelambdaBA}--\eqref{Z1Z2p2B} give
\begin{equation}\label{eq:weak-remaining-singlet-gaps}
\begin{aligned}
\lambda_B(p^A(p_1))-\lambda_*^A(p_1)&=-\frac{4\lambda_0\mathfrak d}{2\mathscr U_2+\mathscr Z_2}\,p_1
  +O_{\mathcal M_\infty,\eta_{\rm m}}( a p_1^2),\\
\lambda_A(p^B(p_1))-\lambda_*^B(p_1)&=\frac{4\lambda_0\mathfrak d}{2\mathscr U_2-\mathscr Z_2}\,p_1
  +O_{\mathcal M_\infty,\eta_{\rm m}}(ap_1^2).
\end{aligned}
\end{equation}
Here and below, this subscript means that the remainder constant depends only on the displayed quantities.
The leading coefficients are nonzero by the moment condition \eqref{momentcondition}. Furthermore,
\begin{equation}
    \frac{\eta_{\rm m}}{1+\mathcal M_\infty}\le\mathfrak g_\iota\le\frac{8(1+\mathcal M_\infty)^2}{\eta_{\rm m}},
\qquad \iota\in\{A,B\}.
\end{equation}
Thus, after decreasing $p_0$ by an amount depending only on $\mathcal M_\infty$ and $\eta_{\rm m}$, both differences in
\eqref{eq:weak-remaining-singlet-gaps} are nonzero for $0<|p_1|<p_0$.

By \eqref{eq:weak-tuned-eigenvalue-bound} and \eqref{differencelambdaBA}, we may decrease $p_0$ so that the four eigenvalues $\lambda_A(p),\lambda_B(p),\lambda_E(p),\lambda_E(p)$ for $p=p^\iota(p_1)$ are contained in 
\begin{equation}\label{interval3454}
    \left( \frac{3}{4}\lambda_0,\frac{5}{4}\lambda_0 \right).
\end{equation}
Lemma \ref{lem:perturbationsplitweak} states that the interval $I=(\frac{1}{2}\lambda_0,\frac{3}{2}\lambda_0)$
contains exactly four eigenvalues, counted with multiplicity. The separation of the remaining singlet therefore proves that the common eigenvalue has multiplicity exactly three.

Next, we identify the distance between $\lambda_*^\iota$ and the remaining spectrum. By \eqref{interval3454}, every eigenvalue outside $I$ is then at distance at least $\frac{1}{4}\lambda_0$ from $\lambda_*^\iota(p_1)$.
On the other hand, the distances in \eqref{eq:weak-remaining-singlet-gaps} tend to zero. Their uniform bounds allow us to decrease $p_0$ further so that both distances are less than $\frac{1}{4}\lambda_0$. Consequently, the remaining singlet is the nearest eigenvalue outside the triple, and
\begin{equation}
    \begin{aligned}
&\operatorname{dist}\!\left( \lambda_*^A(p_1), \sigma_{p_1}^A\setminus\{\lambda_*^A(p_1)\}\right)
=\bigl|\lambda_B(p^A(p_1))-\lambda_*^A(p_1)\bigr|,\\
&\operatorname{dist}\!\left( \lambda_*^B(p_1), \sigma_{p_1}^B\setminus\{\lambda_*^B(p_1)\}\right)
=\bigl|\lambda_A(p^B(p_1))-\lambda_*^B(p_1)\bigr|.
\end{aligned}
\end{equation}

The estimate \eqref{eq:weak-remaining-singlet-gaps} now gives a constant $K_0>0$, depending only on $\mathcal M_\infty$ and $\eta_{\rm m}$, such that
\begin{equation}
    \left|\operatorname{dist}\!\left( \lambda_*^\iota(p_1), \sigma_{p_1}^\iota\setminus\{\lambda_*^\iota(p_1)\}\right)
-\lambda_0\mathfrak g_\iota|p_1|\right|\le K_0\lambda_0p_1^2.
\end{equation}
Since $\mathfrak g_\iota\ge \frac{\eta_{\rm m}}{1+\mathcal{M}_\infty}$, division by $\lambda_0\mathfrak g_\iota|p_1|$ proves
\begin{equation}
    \left|\frac{ \operatorname{dist}\!\left(  \lambda_*^\iota(p_1),  \sigma_{p_1}^\iota\setminus\{\lambda_*^\iota(p_1)\}
 \right)}{\lambda_0\mathfrak g_\iota|p_1|}-1\right|\le K|p_1|,
\end{equation}
where we may take $K=\frac{K_0(1+\mathcal M_\infty)}{\eta_{\rm m}}$. Decreasing $p_0$ so that $Kp_0\le1/2$ gives \eqref{twosidegapestimate2} from \eqref{twosidegapestimate}. It also proves that the threefold eigenvalue is isolated and positive.

\textit{Step 3. Coupling matrices and the Dirac condition.}
Fix $\iota\in\{A,B\}$ and $0<|p_1|<p_0$, and write
$p=p^\iota(p_1)$. Set
\begin{equation}
    \varsigma_A:=1,\qquad \varsigma_B:=-1,
\end{equation}
and choose
\begin{equation}
    (f_0,f_1,f_2)=\bigl(u_\iota(p),u_{E_1}(p),\varsigma_{\iota}u_{E_2}(p)\bigr).
\end{equation}
By Lemma \ref{lem:perturbationsplitweak} and Step~2, these functions form a real basis of the threefold eigenspace and satisfy the normalization \eqref{normalweakperb}.

Recall that
\begin{equation}
    (C_j)_{rs}=\int_Y \bigl(f_rA_p\nabla f_s-f_sA_p\nabla f_r\bigr) \cdot e_j\,dx.
\end{equation}
The function $f_0$ is even in both coordinates, $f_1$ is odd in $x_1$ and even in $x_2$, and $f_2$ is even
in $x_1$ and odd in $x_2$. Since the scalar coefficient in $A_p$ is even in both coordinates, these parities imply
\begin{equation}\label{vanishentry}
    (C_1)_{02}=(C_2)_{01}=(C_1)_{12}=(C_2)_{12}=0.
\end{equation}
The matrices are also antisymmetric.

Define
\begin{equation}
    \gamma_\iota(p_1):=(C_1)_{01}.
\end{equation}
We show that $(C_2)_{02}=\gamma_\iota(p_1)$. The chosen functions satisfy
\begin{equation}
    f_0(Rx)=\varsigma_{\iota}f_0(x),\qquad f_2(Rx)=\varsigma_{\iota}f_1(x).
\end{equation}
By the chain rule for weak derivatives,
\begin{equation}
    (\partial_2f_0)(Rx)=\varsigma_\iota\partial_1f_0(x),\qquad
(\partial_2f_2)(Rx)=\varsigma_\iota\partial_1f_1(x).
\end{equation}
Writing $A_p(x)=c_p(x)I_2$, we have $c_p(Rx)=c_p(x)$. Therefore, the change of variables $x=Ry$ gives
\begin{equation}
    \begin{aligned}
(C_2)_{02}&=\int_Y c_p(x) \bigl(f_0(x)\partial_2f_2(x)-f_2(x)\partial_2f_0(x)\bigr)\,dx\\
&=\int_Y c_p(y) \bigl(f_0(y)\partial_1f_1(y)-f_1(y)\partial_1f_0(y)\bigr)\,dy\\
&=\gamma_\iota(p_1).
\end{aligned}
\end{equation}
Together with \eqref{vanishentry} and antisymmetry, this proves
\begin{equation}\label{eq:weak-coupling-normal-form}
C_1=\gamma_\iota(p_1)\begin{pmatrix}
0&1&0\\
-1&0&0\\
0&0&0
\end{pmatrix},\qquad C_2=\gamma_\iota(p_1)\begin{pmatrix}
0&0&1\\
0&0&0\\
-1&0&0
\end{pmatrix}.
\end{equation}

For the unperturbed modes,
\begin{equation}
    a\bigl( \phi_\iota\partial_1\phi_{E_1} -\phi_{E_1}\partial_1\phi_\iota\bigr)
=2\sqrt2\pi a\bigl[ 1+\varsigma_\iota\cos(2\pi x_1)\cos(2\pi x_2)\bigr].
\end{equation}
The product of cosines has zero integral over $Y$. Hence, for both $\iota=A$ and $\iota=B$,
\begin{equation}
    a\int_Y \bigl(  \phi_\iota\partial_1\phi_{E_1}  -\phi_{E_1}\partial_1\phi_\iota \bigr)\,dx=2\sqrt2\pi a.
\end{equation}

Lemma \ref{lem:perturbationsplitweak} gives
\begin{equation}
    \|u_\iota(p)-\phi_\iota\|_{H^1(Y)}+\|u_{E_1}(p)-\phi_{E_1}\|_{H^1(Y)}\le C\epsilon_\iota(p_1).
\end{equation}
Since $\|A_p-aI_2\|_\infty\le a\epsilon_\iota(p_1)$, Cauchy--Schwarz yields
\begin{equation}
    \begin{aligned}
\bigl|\gamma_\iota(p_1)-2\sqrt2\pi a\bigr|\le{}&C\|A_p-aI_2\|_\infty\\
&+Ca\bigl( \|u_\iota(p)-\phi_\iota\|_{H^1(Y)} +\|u_{E_1}(p)-\phi_{E_1}\|_{H^1(Y)}\bigr)\\
\le{}&Ca\epsilon_\iota(p_1).
\end{aligned}
\end{equation}
Since $\epsilon_\iota(p_1)\le\mathcal M_\infty L|p_1|$, a further decrease of $p_0$, depending only on $\mathcal M_\infty$ and $\eta_{\rm m}$, ensures
\begin{equation}\label{eq:weak-coupling-size}
\gamma_\iota(p_1)\ge\sqrt2\pi a>0.
\end{equation}
This proves part~(b).
\end{proof}

\subsection{Effective dynamics near the tuned Dirac points}
\label{subsec:weakdynamics}

Assume the moment condition \eqref{momentcondition} and retain the notation of Theorem~\ref{thm:weak-dirac}. For $\iota\in\{A,B\}$ and $0<|p_1|<p_0$, define the tuned coefficients and the associated frequency by
\begin{equation}
    A_{p_1}^\iota:=A_{p^\iota(p_1)},\qquad
    b_{p_1}^\iota:=b_{p^\iota(p_1)},\qquad
    \omega_*^\iota(p_1):=\sqrt{\lambda_*^\iota(p_1)}.
\end{equation}
Choose the real normalized basis  $\{f_0,f_1,f_2\}$ from part~(b) of Theorem~\ref{thm:weak-dirac}.
Then, the effective envelope $\Env$ defined by \eqref{evolutionequationE} is given by
\begin{equation}\label{eq:weakeffectiveevolution}
    \widehat{\Env}(\xi,T)=e^{-\mathrm iT\mathsf D_{p_1}^\iota(\xi)}\widehat{\Env_0}(\xi),
\end{equation}
where
\begin{equation}\label{eq:weakexpliciteffectivesymbol}
    \mathsf D_{p_1}^\iota(\xi):=-\frac{\mathrm i}{2\omega_*^\iota(p_1)}(\xi_1C_1+\xi_2C_2)=-\frac{\mathrm i\gamma_\iota(p_1)}{2\omega_*^\iota(p_1)}\begin{pmatrix}
        0&\xi_1&\xi_2\\
        -\xi_1&0&0\\
        -\xi_2&0&0
    \end{pmatrix}.
\end{equation}
By \eqref{eq:weak-tuned-eigenvalue-bound} and \eqref{gammaXp1}, the propagation speed satisfies
\begin{equation}
    \frac{\gamma_\iota(p_1)}{2\omega_*^\iota(p_1)}=\sqrt{\frac a2}\,[1+O(\epsilon_\iota(p_1))],
\end{equation}
where the remainder constant is universal.

As in Section~\ref{sec:effective}, define
\begin{equation}
    \mathrm{WP}_\delta^{\Env}(x,t):=\delta\sum_{r=0}^2f_r(x)E^{(r)}(\delta x,\delta t),\qquad
    \mathrm{WP}_\delta^{\Env_0}:=\mathrm{WP}_\delta^{\Env}(\cdot,0).
\end{equation}
The approximation error is
\begin{equation}\label{eq:weakdiscrepancy}
    e_{\delta,p_1}^\iota(\cdot,t):=e^{-\mathrm it\sqrt{L_{p_1}^\iota}}\mathrm{WP}_\delta^{\Env_0}
    -e^{-\mathrm it\omega_*^\iota(p_1)}\mathrm{WP}_\delta^{\Env}(\cdot,t),
\end{equation}
where $L_{p_1}^\iota$ is the nonnegative self-adjoint operator in $L^2_{b_{p_1}^\iota}(\mathbb R^2)$ given formally by
\begin{equation}
    L_{p_1}^\iota=-\frac{1}{b_{p_1}^\iota}\nabla\cdot\bigl(A_{p_1}^\iota\nabla\bigr).
\end{equation}

\begin{theorem}\label{cor:weakdynamics}
Under the assumptions of Theorem~\ref{thm:weak-dirac}, after decreasing $p_0$ by an amount depending only on
$\mathcal M_\infty$ and $\eta_{\rm m}$, there is a universal constant $C>0$ and, for each integer $N\ge2$, a constant
$C_N>0$ depending only on $N$, such that
\begin{equation}\label{eq:weakprofiledynamics}
\begin{aligned}
    \|e_{\delta,p_1}^\iota(\cdot,t)\|_{L^2_{b_{p_1}^\iota}(\mathbb R^2)}&\le C\frac{\delta}{\mathfrak g_\iota|p_1|}\|\Env_0\|_{\dot H^1}+C\sqrt a\,\frac{\delta^2t}{\mathfrak g_\iota|p_1|}\|\Env_0\|_{\dot H^2} \\
    &\quad +C_N\left(\frac{\delta}{\mathfrak g_\iota|p_1|}\right)^N\|\Env_0\|_{\dot H^N}.
\end{aligned}
\end{equation}
This estimate holds for every $a>0$, $\iota\in\{A,B\}$,
$0<|p_1|<p_0$, $\delta>0$, $t\ge0$, and
$\Env_0\in\mathcal S(\mathbb R^2;\mathbb C^3)$.
The constant $p_0$ depends only on $\mathcal M_\infty$
and $\eta_{\rm m}$, and is independent of $a,\delta,t,\Env_0$, and $N$.
\end{theorem}

\begin{remark}
For a fixed initial datum $\Env_0$, the approximation converges uniformly on $0\le t\le t_\delta$ provided
\begin{equation}\label{eq:weakvalidtimes}
    \delta=o(\mathfrak{g}_\iota|p_1|),\qquad t_\delta =o\left(\frac{\mathfrak{g}_\iota|p_1|}{\sqrt a\,\delta^2}\right).
\end{equation}
In particular, for fixed $a>0$ and fixed nonzero $p_1$, the approximation converges on $0\le t\le C\delta^{-2+\nu}$ for every fixed $C>0$ and $\nu>0$.
\end{remark}

\begin{proof}[Proof of Theorem~\ref{cor:weakdynamics}]
Theorem~\ref{thm:weak-dirac} gives an isolated threefold eigenvalue at $\Gamma$, and $\gamma_\iota(p_1)>0$.
The coupling matrices give $G=\gamma_\iota(p_1)^2 I_2>0$. Thus, Proposition~\ref{prop:generalcriteria} gives a threefold Dirac point, and the general effective dynamics
theorem (Theorem \ref{thm:effective-dynamics}) applies. It remains to specify and estimate the quantities entering Theorem \ref{thm:effective-dynamics}. 

Using the notation $\sigma_{p_1}^\iota$ for the periodic spectrum from Theorem~\ref{thm:weak-dirac}, define the spectral gap at $\Gamma$ and its relative value by
\begin{equation}\label{eq:weak-dynamic-gap}
\begin{aligned}
    \Delta^\iota(p_1)&:=\operatorname{dist}\!\left(\lambda_*^\iota(p_1),\sigma_{p_1}^\iota\setminus\{\lambda_*^\iota(p_1)\}
    \right),\\
    d^\iota(p_1)&:=\min\left\{1,\frac{\Delta^\iota(p_1)}{\lambda_*^\iota(p_1)}\right\}.
\end{aligned}
\end{equation}
In particular, $\Delta^\iota(p_1)$ is the actual spectral gap, whereas $\mathfrak g_\iota$ is the leading dimensionless
coefficient, characterized in \eqref{twosidegapestimate}.

Define the coefficient quantity $s^\iota(p_1)$ as in \eqref{coefficientbound}, \textit{i.e.},
\begin{equation}\label{eq:weakactuals}
\begin{aligned}
    \bigl(s^\iota(p_1)\bigr)^2:=\sup_{0\ne v\in H^1_\#(Y)}\frac{\displaystyle\int_Y \lambda_{\max}(A_{p_1}^\iota)|v|^2\,dx}{\displaystyle\frac{1}{\lambda_*^\iota(p_1)}\int_Y A_{p_1}^\iota\nabla v\cdot
    \overline{\nabla v}\,dx+\int_Y b_{p_1}^\iota|v|^2\,dx}.
\end{aligned}
\end{equation}
For the absolute constant $c=1/100$ used in Lemma~\ref{componentk0est}, set
\begin{equation}\label{eq:weakprofileradius}
    r^\iota(p_1):=\min\left\{1,\frac{c\omega_*^\iota(p_1)d^\iota(p_1)}{s^\iota(p_1)}\right\}.
\end{equation}

The definition of
$\epsilon_\iota(p_1)$ gives
\begin{equation}\label{eq:weak-coefficient-comparison}
\begin{aligned}
    a(1-\epsilon_\iota(p_1))I_2
    &\le A_{p_1}^\iota
    \le a(1+\epsilon_\iota(p_1))I_2,\\
    1-\epsilon_\iota(p_1)
    &\le b_{p_1}^\iota
    \le 1+\epsilon_\iota(p_1).
\end{aligned}
\end{equation}
Testing \eqref{eq:weakactuals} on $v=1$ gives the lower bound below. Discarding the nonnegative gradient term
in its denominator gives the upper bound. Consequently,
\begin{equation}\label{eq:weakprofilesbound}
    a\frac{1-\epsilon_\iota(p_1)}{1+\epsilon_\iota(p_1)}\le \bigl(s^\iota(p_1)\bigr)^2\le
    a\frac{1+\epsilon_\iota(p_1)}{1-\epsilon_\iota(p_1)}.
\end{equation}

By \eqref{eq:weak-tuned-eigenvalue-bound} and \eqref{boundofepsXp1} in Theorem~\ref{thm:weak-dirac} and its proof,
\begin{equation}
    |\lambda_*^\iota(p_1)-4\pi^2 a|\le Ca\epsilon_\iota(p_1),\qquad \epsilon_\iota(p_1)\le C_{\mathcal M_\infty,\eta_{\rm m}}|p_1|,
\end{equation}
where $C$ is universal. Hence, after decreasing $p_0$ by an amount depending only on $\mathcal M_\infty$ and $\eta_{\rm m}$,
\begin{equation}\label{eq:weak-frequency-s-comparison}
\begin{aligned}
    C^{-1}\sqrt a&\le \omega_*^\iota(p_1)\le C\sqrt a,\\
    C^{-1}\sqrt a&\le s^\iota(p_1)\le C\sqrt a,
\end{aligned}
\end{equation}
with universal constants $C>0$.

Combining \eqref{twosidegapestimate2} with $\frac{1}{2}\lambda_0 \leq \lambda_*^\iota(p_1)\leq \frac{3}{2}\lambda_0$ gives
\begin{equation}\label{relativeboundbygX}
    C^{-1}\mathfrak g_\iota|p_1|\le\frac{\Delta^\iota(p_1)}{\lambda_*^\iota(p_1)}\le C\mathfrak g_\iota|p_1|.
\end{equation}
Since the upper bound for $\mathfrak g_\iota$ depends only on $\mathcal M_\infty$ and $\eta_{\rm m}$, we may therefore decrease $p_0$, with the same dependence, so that 
\begin{equation}
    d^\iota(p_1)=\frac{\Delta^\iota(p_1)}{\lambda_*^\iota(p_1)}, \qquad  r^\iota(p_1)=\frac{c\omega_*^\iota(p_1)d^\iota(p_1)}{s^\iota(p_1)}.
\end{equation}
Using \eqref{eq:weak-frequency-s-comparison} and \eqref{relativeboundbygX}, we obtain
\begin{equation}\label{eq:weak-gap-radius-comparison}
\begin{aligned}
    C^{-1}\mathfrak g_\iota|p_1|&\le d^\iota(p_1)\le C\mathfrak g_\iota|p_1|,\\
    C^{-1}\mathfrak g_\iota|p_1|&\le r^\iota(p_1)\le C\mathfrak g_\iota|p_1|.
\end{aligned}
\end{equation}
All constants $C$ here are universal.

Substituting the estimates \eqref{eq:weak-frequency-s-comparison} and \eqref{eq:weak-gap-radius-comparison}  into Theorem~\ref{thm:effective-dynamics} gives the desired estimate \eqref{eq:weakprofiledynamics}.
\end{proof}

\begin{proposition}[Axial dispersion of the central band]\label{prop:weak-axial-curvature}
Fix $a>0$ and either tuning curve $\iota\in\{A,B\}$ of Theorem~\ref{thm:weak-dirac}. For each sufficiently small $p_1\ne0$ there is a coefficient $c_\iota(p_1)>0$ such that for $j=1,2$,
\begin{equation}\label{eq:weak-axial-curvature}
 \lambda_c(k e_j)=\lambda_*^\iota(p_1)+c_\iota(p_1)k^2+O_{p_1}(k^4),
 \qquad c_\iota(p_1)\longrightarrow a\quad (p_1\to0).
\end{equation}
The expansion in $k$ is taken at fixed $p_1$.
\end{proposition}
\begin{proof}
For $\alpha=k e_1$, the Bloch fiber preserves the subspace of functions odd in $x_2$. At $(p,k)=(0,0)$ the eigenvalue $4\pi^2 a$ is simple in this subspace, with eigenfunction $\phi_{E_2}=\sqrt{2}\sin(2\pi x_2)$ from \eqref{eq:weakconstantquartet}. For $p=0$ the same function is an eigenfunction at $k e_1$ with eigenvalue $4\pi^2a+ak^2$. Analytic perturbation of the restricted form pencil thus gives a simple real analytic eigenvalue near $(p,k)=(0,0)$. Reflection in $x_1$ preserves the subspace and takes $k$ to $-k$. The eigenvalue is therefore even in $k$, and its quadratic coefficient is continuous in $p$ and equals $a$ at $p=0$.

At a fixed tuned $p=p^\iota(p_1)$ with $p_1\ne0$, this branch passes through $\lambda_*^\iota(p_1)$ with zero slope. By \eqref{diracconegeneralcriteria} and \eqref{eq:weakexpliciteffectivesymbol}, the two other branches have nonzero one-sided slopes. Hence, the restricted branch is the middle band for all sufficiently small $|k|$ at that $p$. This proves \eqref{eq:weak-axial-curvature} for $j=1$. A quarter-turn gives the result for $j=2$.
\end{proof}
In frequency, the coefficient of $k^2$ is $c_\iota(p_1)/(2\omega_*^\iota(p_1))$. This positive curvature is a second-order correction to the horizontal line in Figure~\ref{fig:threefold-cone}(b).

\subsection{Examples verifying the moment conditions}\label{subsec:weakexamples}

We verify \eqref{momentcondition} for trigonometric profiles and for piecewise-constant inclusions. We also examine the coupling to a uniform source in \eqref{eq:dz-forced-cell}. On the $B_1$ curve the singlet changes sign under the quarter-turn $R$, so its cell mean vanishes. On the $A_1$ curve the mean depends on the choice of profiles.

\subsubsection{Trigonometric perturbations}

Set
\begin{equation}\label{perturbAequal0}
    \sigma_1 = \sigma_2 = 0, 
\end{equation}
and
\begin{equation}\label{perturbBtrigono}
    \beta_1(x) = \cos(4\pi x_1)+\cos(4\pi x_2), \qquad \beta_2(x)=\cos(2\pi x_1)\cos(2\pi x_2).
\end{equation}
Fourier orthogonality gives
\begin{equation}
   (\mathscr M_1,\mathscr U_1,\mathscr Z_1)=\left(0,\frac{1}{2},0\right),\qquad (\mathscr M_2,\mathscr U_2,\mathscr Z_2)=\left(0,0,\frac{1}{2}\right),\qquad \mathfrak d=\frac{1}{4}.
\end{equation}
Thus, the moment condition \eqref{momentcondition} holds. Moreover, the exact tuning graphs satisfy
\begin{equation}\label{eq:weakfouriertuning}
 p_{2}^A(p_1)=-2p_1+O(p_1^2),\qquad p_{2}^B(p_1)=2p_1+O(p_1^2).
\end{equation}
Throughout the small parameter neighborhood, \eqref{estimateoflambdaMp} and Theorem \ref{thm:weak-dirac} give
\begin{equation}\label{eq:weakfourierparameters}
 \begin{aligned}
\lambda_A(p)&=\lambda_0\left[1-\frac{p_1+p_2}{2}\right]+O(a|p|^2),\\
\lambda_B(p)&=\lambda_0\left[1-\frac{p_1-p_2}{2}\right]+O(a|p|^2),\\
\lambda_E(p)&=\lambda_0\left[1+\frac{p_1}{2}\right]+O(a|p|^2).
 \end{aligned}
\end{equation}

For the profiles \eqref{perturbAequal0}--\eqref{perturbBtrigono}, the translations $t_+=(1/2,1/2)$ and $t_-=(1/2,-1/2)$ preserve both coefficients. On either curve $p_2\ne0$ for sufficiently small $p_1\ne0$ by \eqref{eq:weakfouriertuning}. The nonzero Fourier coefficients of $b_p$ at $(1,1)$ and $(1,-1)$ then force every period $s$ to satisfy $s_1+s_2,s_1-s_2\in\mathbb Z$. Thus, $t_+$ and $t_-$ generate the primitive lattice. All four unperturbed modes in \eqref{eq:weakconstantquartet} are antiperiodic under both translations. The projection onto their perturbed cluster commutes with each translation and, for small $p$, maps the unperturbed eigenspace isomorphically onto that cluster. The perturbed modes therefore have the same translation character. The crossing at $\Gamma$ for the unit square is the image under band folding of the primitive-zone corner represented by $(2\pi,0)$. The uniform cell source consequently has zero overlap with the triple on either curve.

A small change of profile gives a crossing at the primitive zone center. Keep $\beta_2$ in \eqref{perturbBtrigono} and replace $\beta_1$ by $\beta_1+\kappa\phi_A$, where $\kappa\ne0$ is fixed and $\phi_A$ is given in \eqref{eq:weakconstantquartet}. By Fourier orthogonality the added term contributes zero to all three moments in \eqref{eq:weakcombinedmoments}, so the leading slopes in \eqref{eq:weakfouriertuning} persist. For $p_1\ne0$ the new Fourier coefficients at $(1,0)$ and $(0,1)$ force both coordinates of a period to be integers. On the $A_1$ curve, choose $f_0=u_A(p)$ as in Lemma~\ref{lem:perturbationsplitweak}. The weighted orthogonality to constants and $f_0=\phi_A+O_{H^1(Y)}(|p_1|)$ give
\begin{equation}\label{eq:trig-primitive-mean}
  m_*=-\int_Y (b_p-1)f_0\,dx=-\kappa p_1+O(p_1^2).
\end{equation}
Here $\|\phi_A\|_{L^2(Y)}=1$ and the original profiles in \eqref{perturbBtrigono} are orthogonal to $\phi_A$. Therefore, Proposition~\ref{prop:dz-bulk} applies to the modified $A_1$ crossing for all sufficiently small nonzero $p_1$.

\begin{remark}
    The following generalization of \eqref{perturbAequal0}--\eqref{perturbBtrigono} concerns the moment condition. The translation lattice must be determined from the individual profiles. 
    Let $\psi,\varphi\in L^\infty_\#(0,1;\mathbb R)$ have even, one-periodic extensions satisfying
\begin{equation}
    \int_0^1\varphi(x)\,dx=0,\qquad a_4:=\int_0^1\psi(x)\cos(4\pi x)\,dx\ne0,\qquad
    b_2:=\int_0^1\varphi(x)\cos(2\pi x)\,dx\ne0.
\end{equation}
Set
\begin{equation}
    \sigma_1=\sigma_2=0,\qquad \beta_1(x)=\psi(x_1)+\psi(x_2),\qquad \beta_2(x)=\varphi(x_1)\varphi(x_2).
\end{equation}
Then, we have
\begin{equation}
   (\mathscr M_1,\mathscr U_1,\mathscr Z_1)=\left(2\int_0^1 \psi(x)\,dx\,,a_4\,,0\right),\qquad (\mathscr M_2,\mathscr U_2,\mathscr Z_2)=\left(0\,,0\,,2b_2^2\right),\qquad \mathfrak d=2a_4b_2^2.
\end{equation}
Thus, the moment condition \eqref{momentcondition} holds.
\end{remark}

\subsubsection{Piecewise-constant perturbations}
\label{subsec:weak-piecewise}

Let $D\subset Y$ be a nonempty $D_4$-invariant Lipschitz inclusion such that
\begin{equation}\label{Din110}
    D\subset \left(-\frac{1}{10},\frac{1}{10}\right)^2.
\end{equation}
Write $\mathbf 1_D$ for the indicator of $D$, and set
\begin{equation}\label{eq:piecewise-profiles}
    \sigma_1=\mathbf 1_D,\qquad \sigma_2=0,\qquad \beta_1=0,\qquad \beta_2=\mathbf 1_D.
\end{equation}
Define the geometric moments
\begin{equation}\label{eq:piecewise-geometric-moments}
    M_D:=|D|,\qquad
    U_D:=\int_D\cos(4\pi x_1)\,dx,\qquad
    Z_D:=2\int_D
        \cos(2\pi x_1)\cos(2\pi x_2)\,dx.
\end{equation}
Then, it follows that
\begin{equation}\label{eq:piecewise-combined-moments}
    (\mathscr M_1,\mathscr U_1,\mathscr Z_1)=(-M_D,U_D,0),\qquad 
    (\mathscr M_2,\mathscr U_2,\mathscr Z_2)=(M_D,U_D,Z_D).
\end{equation}
Consequently,
\begin{equation}
    \mathfrak d=U_DZ_D,\qquad 2\mathscr U_2+\mathscr Z_2=2U_D+Z_D,\qquad 2\mathscr U_2-\mathscr Z_2=2U_D-Z_D.
\end{equation}

We verify the moment condition \eqref{momentcondition}. By \eqref{Din110}, the integrands defining $U_D$ and $Z_D$ are strictly positive on $D$. Hence,
\begin{equation}
    U_D>0,\qquad Z_D>0,
\end{equation}
which gives $\mathfrak d>0$ and $2U_D+Z_D>0$. To determine the sign of $2U_D-Z_D$, the $D_4$ symmetry of $D$ gives
\begin{equation}
\begin{aligned}
    2U_D-Z_D&=\int_D\bigl[\cos(4\pi x_1)+\cos(4\pi x_2)-2\cos(2\pi x_1)\cos(2\pi x_2)\bigr]\,dx \\
    &= 2 \int_D (c_1^2+c_2^2-c_1c_2-1) \,dx,
\end{aligned}
\end{equation}
where $c_j=\cos(2\pi x_j)$, $j=1,2$. If $c_1\ge c_2$, then
\begin{equation}
    c_1^2+c_2^2-c_1c_2 =c_1^2-c_2(c_1-c_2)\le c_1^2\le1,
\end{equation}
with equality only when $c_1=c_2=1$. The case $c_2\ge c_1$ follows by symmetry. Thus $2U_D-Z_D <0$.
The moment condition \eqref{momentcondition} therefore holds.

Theorem~\ref{thm:weak-dirac} now gives two real analytic tuning curves $p_2^A$ and $p_2^B$, independent of $a$,
with
\begin{equation}\label{eq:piecewise-tuning-curves}
    \begin{aligned}
        p_2^A(p_1)&=-\frac{2U_D}{2U_D+Z_D}\,p_1+O_D(p_1^2),\\
        p_2^B(p_1)&=\frac{2U_D}{Z_D-2U_D}\,p_1+O_D(p_1^2).
    \end{aligned}
\end{equation}
Here, the constants in $O_D(\cdot)$ may depend on the fixed inclusion $D$. The first-order expansions of the three eigenvalue branches are
\begin{equation}\label{eq:piecewise-branch-expansions}
    \begin{aligned}
        \lambda_A(p)&=\lambda_0\bigl[1+p_1(M_D-U_D)-p_2(M_D+U_D+Z_D)\bigr]+O\bigl(a|p|^2\bigr),\\
        \lambda_B(p)&=\lambda_0\bigl[1+p_1(M_D-U_D)-p_2(M_D+U_D-Z_D)\bigr]+O\bigl(a|p|^2\bigr),\\
        \lambda_E(p)&=\lambda_0\bigl[1+p_1(M_D+U_D)-p_2(M_D-U_D)\bigr]+O\bigl(a|p|^2\bigr).
    \end{aligned}
\end{equation}

For the piecewise-constant example, the $A_1$ singlet has nonzero cell mean. Indeed, set $p=(p_1,p_2^A(p_1))$ and choose $f_0=u_A(p)$ with the sign in Lemma~\ref{lem:perturbationsplitweak}. Testing its positive-frequency eigenvalue equation with the constant function gives $\int_Y b_p f_0\,dx=0$. Since $b_p=1+p_2^A(p_1)\mathbf 1_D$, we obtain
\begin{equation}\label{eq:weak-mean-piecewise}
 m_*=-p_2^A(p_1)\int_D f_0\,dx.
\end{equation}
By Lemma~\ref{lem:perturbationsplitweak}, $f_0\to\phi_A$ in $H^1(Y)$ as $p_1\to0$, and \eqref{eq:weakconstantquartet} and \eqref{Din110} give $\int_D\phi_A\,dx>0$. The linear coefficient in \eqref{eq:piecewise-tuning-curves} is nonzero. Thus $|m_*|\asymp |p_1|$. Proposition~\ref{prop:dz-bulk} applies for all sufficiently small $p_1\ne0$, and $Z_*\asymp p_1^2$ since the cone speed tends to $\sqrt{a/2}$.

\section{The dilute high-contrast model}\label{sec:dilute-verification}

We consider a dilute inclusion with high contrast in the coefficient $b$. Let $B\subset\mathbb R^2$ be a bounded, connected, $D_4$-invariant Lipschitz domain containing the origin. For sufficiently small $\eta>0$, set
\begin{equation}\label{eq:dilute-weight}
    D_\eta:=\eta B\Subset Y,\qquad A:=I_2,\qquad
    b=b_{\eta,\tau}:=1+(\tau-1)\mathbf{1}_{D_\eta},\qquad \ell_\eta:=\log(1/\eta).
\end{equation}
The contrast is $\tau>1$, and $\mathbf 1_D$ denotes the indicator of a set $D$. The cell and its translates are shown in Figure~\ref{fig:crystal-geometry}. The coefficients are extended periodically. The condition $QB=B$ for every $Q\in D_4$ gives $b_{\eta,\tau}(Qx)=b_{\eta,\tau}(x)$, and $A=I_2$ satisfies \eqref{symmetryk=2}. Thus one centered inclusion gives the symmetry required to preserve the dipole doublet. Remark~\ref{rem:bubble-comparison} explains why nonconstant internal modes remain subwavelength in this high-permittivity scaling. At $\Gamma$, the eigenvalue problem \eqref{prob:blocheigenpair} becomes
\begin{equation}\label{eq:dilute-periodic-problem}
    \int_Y\nabla u\cdot\nabla v\,dx =\lambda\int_Yb_{\eta,\tau}uv\,dx, \qquad v\in H^1_\#(Y;\mathbb R).
\end{equation}
All eigenfunctions used in the construction below are real.

We retain the symmetry spaces $\mathcal V_A$ and $\mathcal V_{E_1}$ from \eqref{symmetryspaceAB1}. Let $\lambda_A(\eta,\tau)$ be the first positive eigenvalue in $\mathcal V_A$, \textit{i.e.}, $\lambda_A(\eta,\tau)$ is the smallest positive number such that 
\begin{equation}\label{eq:dilute-periodic-probleminVA}
    \int_Y\nabla u_A\cdot\nabla v\,dx =\lambda_A(\eta,\tau)\int_Yb_{\eta,\tau}u_Av\,dx, \qquad v\in \mathcal V_A,
\end{equation}
for some nonzero $u_A\in\mathcal V_A$. Let $\lambda_E(\eta,\tau)$ be the smallest eigenvalue for which 
\begin{equation}\label{eq:dilute-periodic-probleminVE1}
    \int_Y\nabla u_{E_1}\cdot\nabla v\,dx =\lambda_E(\eta,\tau)\int_Yb_{\eta,\tau}u_{E_1}v\,dx, \qquad v\in \mathcal V_{E_1},
\end{equation}
for some nonzero real function $u_{E_1}\in \mathcal{V}_{E_1}$. Applying the projections in \eqref{symmetrizationAB} and \eqref{symmetrizationE1} shows that the pairs $(\lambda_A,u_A)$ and $(\lambda_E,u_{E_1})$ also solve \eqref{eq:dilute-periodic-problem} for all periodic test functions. 

Define the finite-energy space
\begin{equation}\label{eq:dilute-E-limit-space}
\mathcal V_{E_1}^{\infty}:=\left\{
\begin{aligned}
U\in H^1_{\mathrm{loc}}(\mathbb R^2;\mathbb R):\;&\nabla U\in L^2(\mathbb R^2),\\
&U(-x_1,x_2)=-U(x_1,x_2),\ U(x_1,-x_2)=U(x_1,x_2)
\end{aligned}
\right\},
\end{equation}
and the corresponding variational constant
\begin{equation}\label{eq:LambdaE-dilute}
    \Lambda_E(B):=\inf_{\substack{U\in\mathcal V_{E_1}^{\infty}\\ \int_BU^2\,dx\neq 0}}
    \frac{\displaystyle\int_{\mathbb R^2}|\nabla U|^2\,dx}{\displaystyle\int_BU^2\,dx}.
\end{equation}

The first positive Neumann eigenvalue in the fully invariant sector of $B$ is
\begin{equation}\label{eq:dilute-Neumann-constant}
    \mu_A^{\mathrm N}(B):=
    \inf_{\substack{0\ne w\in H^1(B;\mathbb R),\,\int_Bw\,dx=0\\
         w\ \text{is }D_4\text{-invariant}}}
       \frac{\int_B|\nabla w|^2\,dx}{\int_Bw^2\,dx}.
\end{equation}
We shall separate this interior threshold from the first dipole eigenvalue.

\begin{assumption}\label{assump:neumanngap}
The reference inclusion $B$ satisfies
\begin{equation}\label{eq:A1-Neumann-gap}
    \mu_A^{\mathrm N}(B)>\Lambda_E(B).
\end{equation}
\end{assumption}

Fix $0<\rho<1$ such that $(1+\rho)\Lambda_E(B)<\mu_A^{\mathrm{N}}(B)$, and set
\begin{equation}\label{eq:critical-contrast-dilute}
    \tau_0(\eta):=\frac{\Lambda_E(B)}{2\pi\eta^2}\ell_\eta,
    \qquad I_\eta:=[(1-\rho)\tau_0(\eta),(1+\rho)\tau_0(\eta)].
\end{equation}
We write $f\asymp_B g$ when $c_Bg\le f\le C_Bg$ for positive shape-dependent constants, and use $O_B$ for estimates with this dependence. Throughout this section, constants may depend on the fixed shape $B$ and the chosen $\rho$, but not on $\eta$ or $\tau$. Once $\rho$ is chosen in terms of $B$, we write this dependence as $C_B$.

\subsection{Limits of the dipole doublet and the capacitary singlet}\label{subsec:dilute-E}

We first determine the limit of the axial eigenfunction after rescaling the inclusion to $B$. A quarter-turn gives the second member of the doublet.

\Needspace{6\baselineskip}
\begin{lemma}\label{lem:dilute-E-convergence}
Let $\Lambda_E(B)$ and $\lambda_E(\eta,\tau)$ be given by \eqref{eq:LambdaE-dilute} and \eqref{eq:dilute-periodic-probleminVE1}, respectively.
\begin{itemize}
    \item[(a)] The infimum \eqref{eq:LambdaE-dilute} has a minimizer unique up to multiplication by a nonzero scalar. We choose the representative $U_1$ satisfying
\begin{equation}\label{U1property}
    \int_BU_1^2\,dx=1,\qquad U_1>0\quad\mathrm{\ on\ }\{x_1>0\}.
\end{equation}
Moreover,
\begin{equation}\label{U1equation}
    -\Delta U_1=\Lambda_E(B)\mathbf{1}_BU_1 \qquad \mathrm{in} \ \mathbb{R}^2,
\end{equation}
and
\begin{equation}\label{eq:dilute-E-dipole}
    U_1(x)=d_B\frac{x_1}{|x|^2}+O(|x|^{-2}) \quad \mathrm{as} \ x\rightarrow \infty,\qquad d_B:=\frac{\Lambda_E(B)}{2\pi}\int_B x_1U_1(x)\,dx>0.
\end{equation}
\item[(b)] The eigenvalue $\lambda_E(\eta,\tau)$ is simple in $\mathcal V_{E_1}$, and
\begin{equation}\label{eq:E-branch-dilute}
    \lambda_E(\eta,\tau)=\frac{\Lambda_E(B)}{\tau\eta^2}[1+o(1)] \qquad \mathrm{as} \ \eta\to0
\end{equation}
uniformly in $\tau\in I_\eta$.
\item[(c)] If the eigenfunction $u_{E_1}$ of $\lambda_E(\eta,\tau)$ is normalized in $L^2_{b_{\eta,\tau}}(Y)$ and
chosen positive on $Y\cap\{x_1>0\}$, then
\begin{equation}\label{eq:dilute-E-mode-convergence}
    q_{\eta,\tau}(x):=\eta\sqrt{\tau}\,u_{E_1}(\eta x)\longrightarrow U_1\quad\text{in }L^2(B)\quad(\eta\to0),
\end{equation}
uniformly for $\tau\in I_\eta$. The rotated function $u_{E_2}(x):=u_{E_1}(R^\top x)$ has the same eigenvalue $\lambda_E(\eta,\tau)$.
\end{itemize}
\end{lemma}

\begin{proof}
Every $U\in\mathcal V_{E_1}^{\infty}$ is odd in $x_1$, so the average of its angular slice on $|x|=r$ vanishes for almost every $r>0$. The polar coordinate representation and the one-dimensional Poincar\'{e}--Wirtinger inequality give the Hardy inequality for $U\in \mathcal{V}_{E_1}^{\infty}$:
\begin{equation}\label{eq:dilute-angular-Hardy}
    \int_{\mathbb R^2}\frac{|U(x)|^2}{|x|^2}\,dx\le\int_{\mathbb R^2}|\nabla U|^2\,dx,
    \qquad U\in\mathcal V_{E_1}^{\infty}.
\end{equation}
For each $R>0$ the inequality bounds the local $L^2(B_R(0))$ norm by $R\|\nabla U\|_{L^2(\mathbb R^2)}$, where $B_R(0)$ is the disk of radius $R$. A minimizing sequence normalized in $L^2(B)$ therefore has a subsequence converging weakly in $H^1_{\rm loc}$ and strongly in $L^2_{\rm loc}$. Its limit has the prescribed parities and unit mass on $B$. Lower semicontinuity of the energy gives a minimizer $\widetilde U$ \cite{evans2022partial}.

The first variation of the Rayleigh quotient yields
\begin{equation}\label{eq:E-proof-variational-equation}
    \int_{\mathbb R^2}\nabla \widetilde{U}\cdot\nabla v\,dx=\Lambda_E(B)\int_B \widetilde{U} v\,dx,
    \qquad v\in\mathcal V_{E_1}^{\infty}.
\end{equation}
Applying the parity projection in \eqref{symmetrizationE1} to a test function shows that \eqref{eq:E-proof-variational-equation} holds for every $v\in C_c^\infty(\mathbb R^2)$. Therefore,
\begin{equation}\label{eq:E-proof-whole-plane-equation}
    -\Delta \widetilde{U}=\Lambda_E(B)\mathbf{1}_B \widetilde{U} \quad\mathrm{ in }\ \mathbb R^2.
\end{equation}

Let $H=\{x_1>0\}$ and $B_+=B\cap H$. For $z=(z_1,z_2)$, set $z^*=(-z_1,z_2)$ and let $\delta_z$ denote the Dirac mass at $z$. The half-plane Dirichlet Green function on $H$ is
\begin{equation}
    G_H(x,z)=\frac1{2\pi}\log\frac{|x-z^*|}{|x-z|}.
\end{equation}
This function satisfies $-\Delta_xG_H(x,z)=\delta_z$ in $H$ and vanishes on $\partial H$. To obtain an integral representation, set $f=\Lambda_E(B)\mathbf1_B\widetilde U$. This compactly supported $L^2$ function has zero integral by oddness. Its Newton potential
\[
 Nf(x)=-\frac{1}{2\pi}\int_B\log|x-z|\, f(z)\,dz
\]
is locally $H^2$, is odd in $x_1$, and satisfies $\nabla Nf(x)=O(|x|^{-2})$ at infinity by the zero integral of $f$. Thus, $\widetilde U-Nf$ is entire harmonic with gradient in $L^2(\mathbb R^2)$. Each derivative vanishes by the mean-value inequality on expanding disks. The difference is an odd constant and therefore zero. Splitting the integral into the two half-planes gives
\begin{equation}\label{eq:E-proof-half-plane-representation}
    \widetilde{U}(x)=\Lambda_E(B)\int_{B_+}G_H(x,z)\widetilde{U}(z)\,dz,
    \qquad x\in H.
\end{equation}

Define the integral operator
\begin{equation}
    (Kg)(x):=\int_{B_+}G_H(x,z)g(z)\,dz, \qquad x\in B_+.
\end{equation}
For $x,z\in H$, we have $|x-z^*|> |x-z|$. Hence, $G_H(x,z)>0$. Moreover, since $B$ is bounded,
\begin{equation}
    0<G_H(x,z)\le C_B\bigl(1+|\log|x-z||\bigr) \qquad \mathrm{on} \ B_+\times B_+
\end{equation}
away from the diagonal. The logarithmic singularity is square-integrable in two dimensions. Hence $K$ is compact
and self-adjoint on $L^2(B_+)$. By the min--max principle, the largest eigenvalue $\kappa_1$ of $K$ is given by
\begin{equation}
    \kappa_1=\max_{0\ne g\in L^2(B_+)}\frac{\langle Kg,g\rangle_{L^2(B_+)}}{\|g\|_{L^2(B_+)}^2}.
\end{equation}
Since $\langle K1,1\rangle_{L^2(B_+)}>0$, we have $\kappa_1>0$. For real $g \in L^2(B_+)$,
\begin{equation}\label{Kgg-Kgg}
    \langle K|g|,|g|\rangle_{L^2(B_+)}-\langle Kg,g\rangle_{L^2(B_+)}
    =\int_{B_+}\int_{B_+}G_H(x,z) \bigl(|g(x)g(z)|-g(x)g(z)\bigr)\,dzdx.
\end{equation}
If $g$ has positive and negative parts of positive measure, the right-hand side of \eqref{Kgg-Kgg} is strictly positive. Thus, every eigenfunction corresponding to $\kappa_1$ has one sign. Choose a nonnegative one, denoted by
\begin{equation}
    g_1 \geq 0 \qquad \mathrm{on} \ B_+.
\end{equation}
The identity $g_1=\kappa_1^{-1}Kg_1$ and the strict positivity of $G_H(x,z)$ imply that $g_1>0$ on $B_+$. If the eigenspace of $\kappa_1$ has dimension at least two, it would contain two orthogonal real eigenfunctions. Both could be chosen
strictly positive, contradicting their orthogonality. Therefore, $\kappa_1$ is simple.

Reflection in $x_2$ preserves $B_+$ and the kernel $G_H$. Hence, the reflection in $x_2$ transforms a positive eigenfunction into another positive eigenfunction with the same norm. Simplicity of $\kappa_1$ gives
\begin{equation}
    g_1(x_1,-x_2)=g_1(x_1,x_2).
\end{equation}
Normalize $g_1$ by $\int_{B_+}g_1^2\,dx=\frac{1}{2}$. Extend it to a function on $B$ by
\begin{equation}\label{symmetryg1odd}
    g_1^{\rm odd}(x):=\begin{cases}
       g_1(x),&x\in B_+,\\
       -g_1(x^*),&x\in B\cap\{x_1<0\},
    \end{cases}
\end{equation}
and set it equal to zero outside $B$. Define
\begin{equation}\label{defofU1byG}
    U_1(x):=-\frac1{2\pi\kappa_1} \int_B\log|x-z|\,g_1^{\rm odd}(z)\,dz, \qquad x\in \mathbb{R}^2.
\end{equation}
The compactly supported odd density in \eqref{defofU1byG} has zero mean. The same decay and regularity estimates as for $Nf$ show that $U_1\in\mathcal V_{E_1}^{\infty}$. On $H$ it has the representation
\begin{equation}
    U_1(x)=\frac1{\kappa_1} \int_{B_+}G_H(x,z)g_1(z)\,dz>0.
\end{equation}
On $B_+$, the equation $Kg_1=\kappa_1g_1$ gives $U_1=g_1$; oddness gives $U_1=g_1^{\rm odd}$ on $B$.
Consequently,
\begin{equation}\label{equationforU1}
    \int_BU_1^2\,dx=1,\qquad -\Delta U_1=\kappa_1^{-1}\mathbf{1}_BU_1.
\end{equation}

We now show that $\Lambda_E(B)=\kappa_1^{-1}$. Test \eqref{equationforU1} with a radial cutoff of $U_1$ equal to $U_1$ on $B_R(0)$ and zero outside $B_{2R}(0)$, with cutoff gradient bounded by $C/R$. By \eqref{eq:dilute-angular-Hardy} the cutoff error tends to zero as $R\to\infty$. This gives
\begin{equation}
    \int_{\mathbb R^2}|\nabla U_1|^2\,dx =\kappa_1^{-1}\int_BU_1^2\,dx =\kappa_1^{-1}.
\end{equation}
By the definition \eqref{eq:LambdaE-dilute} of $\Lambda_E(B)$, we get $\Lambda_E(B)\le\kappa_1^{-1}$.
Conversely, the restriction $g:=\widetilde{U}|_{B_+}$ of the minimizer already constructed is nonzero. Equation
\eqref{eq:E-proof-half-plane-representation} gives $Kg=\Lambda_E(B)^{-1}g$. Since $\kappa_1$ is the largest
eigenvalue, $\Lambda_E(B)^{-1}\le\kappa_1$, giving the reverse inequality. Therefore $\Lambda_E(B)=\kappa_1^{-1}$. Simplicity of $\kappa_1$ and the integral representation show that every minimizer is proportional to $U_1$. This proves \eqref{U1property} and \eqref{U1equation} in part (a).

For the far field expansion of $U_1$, using the definition \eqref{defofU1byG} for $U_1$, the fact that $U_1=g_1^{\rm odd}$ on $B$, and the expansion 
\begin{equation}
    \log|x-z|=\log|x|-\frac{x\cdot z}{|x|^2}+O_B(|x|^{-2}),
\end{equation}
uniformly in $z\in B$, we obtain
\begin{equation}
    U_1(x)=\frac{\Lambda_E(B)}{2\pi}\frac{x\cdot\int_BzU_1(z)\,dz}{|x|^2}+O_B(|x|^{-2}).
\end{equation}
Here $\int_BU_1=0$ by oddness, and $\int_Bz_2U_1(z)\,dz=0$ by evenness in $z_2$. Moreover,
\begin{equation}
    \int_Bz_1U_1(z)\,dz=2\int_{B_+}z_1U_1(z)\,dz>0.
\end{equation}
This proves the far field expansion \eqref{eq:dilute-E-dipole} in part~(a).

We now prove part (b). Restrict a function in $\mathcal V_{E_1}$ to the quarter-cell $\Omega=(0,\frac{1}{2})^2$. Oddness and periodicity give zero trace on both vertical sides $\{0\}\times [0,\frac{1}{2}]$ and $\{\frac{1}{2}\}\times [0,\frac{1}{2}]$. For an eigenfunction, even reflection and the weak equation give natural Neumann conditions on both horizontal sides $[0,\frac{1}{2}]\times \{0\}$ and $[0,\frac{1}{2}]\times \{\frac{1}{2}\}$. Thus, the first eigenvalue $\lambda_E(\eta,\tau)$ of problem \eqref{eq:dilute-periodic-probleminVE1} in $\mathcal V_{E_1}$ is the first eigenvalue of
\begin{equation}\label{eigenvalueprobonOmega}
    -\Delta u=\lambda b_{\eta,\tau}u\quad\text{in }\Omega, \qquad
    u=0\ \text{on the vertical sides},\qquad \partial_{x_2}u=0\ \text{on the horizontal sides}.
\end{equation}
The Rayleigh principle gives
\begin{equation}\label{Rayleighprinciple}
    \lambda_E(\eta,\tau)=\inf_{v\in X\setminus \{0\} } \frac{\int_\Omega |\nabla v|^2\,dx}{\int_\Omega b_{\eta,\tau}v^2\,dx},
\end{equation}
where
\begin{equation}
    X:=\left\{v\in H^1(\Omega): v\text{ has zero trace on both vertical sides} \right\}.
\end{equation}
An eigenfunction $w$ of the eigenvalue problem \eqref{eigenvalueprobonOmega} associated with the eigenvalue $\lambda_E(\eta,\tau)$ can be chosen nonnegative, since the absolute value $|w|$ is still in $X$ and preserves the Rayleigh quotient in the right-hand side of \eqref{Rayleighprinciple}. Consequently, \eqref{eigenvalueprobonOmega} implies that
\begin{equation}
    -\Delta w=\lambda b_{\eta,\tau}w\geq 0 \quad\text{in }\Omega.
\end{equation}
The strong maximum principle for nonnegative superharmonic functions gives $w>0$ in $\Omega$.
For any real ground state, its absolute value is also a ground state and is strictly positive in the interior. Local continuity thus forces every ground state to have one sign. Two independent ground states could be chosen to be orthogonal in the weighted inner product, contradicting this sign property. Hence $\lambda_E(\eta,\tau)$ is simple in $\mathcal V_{E_1}$.

Set $Y_\eta:=\eta^{-1}Y$. Changing variables gives
\begin{equation}\label{eq:E-proof-scaled-quotient}
    \tau\eta^2\lambda_E(\eta,\tau)=\inf_{0\ne V}
    \frac{\displaystyle\int_{Y_\eta}|\nabla V|^2\,dx}{\displaystyle\int_BV^2\,dx+\tau^{-1}\int_{Y_\eta\setminus B}V^2\,dx},
\end{equation}
where the infimum is over $V\in H^1_\#(Y_\eta;\mathbb R)$ with period lattice $\eta^{-1}\mathbb Z^2$, odd in $x_1$ and even in $x_2$. Such a function has zero cell average, and the scaled Poincar\'e inequality for periodic functions gives
\begin{equation}\label{eq:E-proof-exterior-bound}
    \tau^{-1}\int_{Y_\eta\setminus B}V^2\,dx\le\frac1{4\pi^2\tau\eta^2}\int_{Y_\eta}|\nabla V|^2\,dx.
\end{equation}

Let $\chi_M \in C_0^\infty(\mathbb{R}^2;[0,1])$ be a smooth radial cutoff satisfying
\begin{equation}
    \chi_M=1 \quad \mathrm{on}\ B_M(0), \qquad \mathrm{supp}\,\chi_M\subset B_{2M}(0),\qquad |\nabla\chi_M|\le \frac{C}{M}.
\end{equation}
For $U_M:=\chi_MU_1$, the Hardy inequality \eqref{eq:dilute-angular-Hardy} gives
\begin{equation}
    \|\nabla U_M-\nabla U_1\|_{L^2(\mathbb R^2)}^2\le2\int_{\{|x|>M\}}|\nabla U_1|^2\,dx
    +C\int_{\{M<|x|<2M\}}\frac{U_1^2}{|x|^2}\,dx\rightarrow0
\end{equation}
as $M\rightarrow \infty$. Assume that $B\subset B_{M_B}(0)$. If $M>M_B$, then $\int_BU_M^2\,dx=1$. For fixed $M$ and sufficiently small $\eta$,
$B_{2M}(0)\Subset Y_\eta$, and so the periodic extension of $U_M$ is admissible in \eqref{eq:E-proof-scaled-quotient}. Therefore,
\begin{equation}
    \tau\eta^2\lambda_E(\eta,\tau) \le\int_{\mathbb R^2}|\nabla U_M|^2\,dx
\end{equation}
for every $\tau \in I_\eta$. Letting $\eta\to0$ and then $M\to\infty$ gives
\begin{equation}\label{eq:E-proof-uniform-limsup}
    \limsup_{\eta\to0}\sup_{\tau\in I_\eta} \tau\eta^2\lambda_E(\eta,\tau)\le\Lambda_E(B).
\end{equation}

Take arbitrary sequences $\eta_j\to0$ and $\tau_j\in I_{\eta_j}$. Write
\begin{equation}
    a_j=\tau_j\eta_j^2\lambda_E(\eta_j,\tau_j), \qquad L:=\liminf_{j\rightarrow \infty} a_j.
\end{equation}
The estimate \eqref{eq:E-proof-uniform-limsup} shows that $\{a_j\}_{j}$ is uniformly bounded, so $L$ is finite. By the definition of liminf, there is a subsequence along which $a_j \rightarrow L $, after relabelling. Let $u_j$ be the principal axial eigenfunction at $(\eta_j,\tau_j)$, positive on the right half-cell and normalized by
\begin{equation}
    \tau_j \int_{D_{\eta_j}} u_j(x)^2\,dx + \int_{Y\setminus D_{\eta_j}} u_j(x)^2\,dx =1.
\end{equation}
Define
\begin{equation}
    q_j(x) := \eta_j \sqrt{\tau_j} u_j(\eta_j x).
\end{equation}
By a change of variables, we have
\begin{equation}\label{normofqj}
    \int_B q_j^2 \,dx + \tau_j^{-1} \int_{Y_{\eta_j}\setminus B} q_j^2\,dx=1,
\end{equation}
and
\begin{equation}\label{valueofnablaqj}
    \int_{Y_{\eta_j}} |\nabla q_j|^2 \,dx= a_j.
\end{equation}

By \eqref{valueofnablaqj} and the Poincar\'{e} inequality \eqref{eq:E-proof-exterior-bound}, we get
\begin{equation}
    0\leq \tau_j^{-1} \int_{Y_{\eta_j}\setminus B} q_j^2\,dx \leq \frac{a_j}{4\pi^2\tau_j \eta_j^2} \leq \frac{C}{\tau_j \eta_j^2}.
\end{equation}
Using $\tau_j \eta_j^2 \geq C\ell_{\eta_j}$ and \eqref{normofqj}, we get
\begin{equation}\label{eq:E-proof-inclusion-mass}
    \tau_j^{-1}\int_{Y_{\eta_j}\setminus B}q_j^2\,dx=O_B(\ell_{\eta_j}^{-1}),
    \qquad\int_Bq_j^2\,dx=1+o(1).
\end{equation}

For each fixed $M$, the ball $B_M(0)$ lies in $Y_{\eta_j}$ for all sufficiently large $j$. Since $q_j$ is odd in $x_1$, the Poincar\'e inequality yields
\begin{equation}
    \|q_j\|_{H^1(B_M(0))}^2 \le C_M\int_{Y_{\eta_j}}|\nabla q_j|^2\,dx = C_M a_j \leq C_M.
\end{equation}
The standard diagonal argument gives a subsequence and a function $Q$ such that
\begin{equation}
    q_j\rightharpoonup Q\quad\text{in }H^1_{\rm loc}(\mathbb R^2), \qquad
    q_j\rightarrow Q\quad\text{in }L^2_{\rm loc}(\mathbb R^2).
\end{equation}
It follows that $Q$ is still odd in $x_1$ and even in $x_2$, and \eqref{eq:E-proof-inclusion-mass} gives $\int_BQ^2\,dx=1$. Lower semicontinuity and \eqref{valueofnablaqj} give
\begin{equation}
    \int_{\mathbb R^2}|\nabla Q|^2\,dx \le \liminf_{j\to\infty} a_j = L.
\end{equation}
Thus, $Q\in\mathcal V_{E_1}^{\infty}$, and $\Lambda_E(B)\le L $. Since the original sequences were arbitrary, this lower bound is uniform in $\tau\in I_\eta$. Combining it with \eqref{eq:E-proof-uniform-limsup} gives
\begin{equation}
    \sup_{\tau\in I_\eta}\bigl|\tau\eta^2\lambda_E(\eta,\tau)-\Lambda_E(B)\bigr|\rightarrow0.
\end{equation}
This is precisely \eqref{eq:E-branch-dilute}.

We now prove part (c). For the arbitrary sequences considered above, part~(b) and lower semicontinuity show
\begin{equation}
    \int_{\mathbb R^2}|\nabla Q|^2\,dx=\Lambda_E(B), \qquad \int_BQ^2\,dx=1.
\end{equation}
Hence every subsequential limit is a normalized minimizer. On each fixed compact subset of $H$, the functions $q_j$ are positive for sufficiently large $j$, because the positive rescaled half-cell expands to $H$. Therefore, strong local $L^2$ convergence gives $Q\ge0$ on $H$. By the simplicity and sign choice proved in part~(a),
this limit is $U_1$.

This also proves the uniform convergence on $I_\eta$. Indeed, if it failed, there would be an $\varepsilon>0$ and sequences $\eta_j\to0$, $\tau_j\in I_{\eta_j}$ with $\|q_{\eta_j,\tau_j}-U_1\|_{L^2(B)}\ge\varepsilon$. The preceding compactness argument would produce a subsequence converging strongly in $L^2(B)$ to $U_1$, a contradiction. Thus,
\begin{equation}
    \sup_{\tau\in I_\eta} \|q_{\eta,\tau}-U_1\|_{L^2(B)}\longrightarrow0,
\end{equation}
which proves \eqref{eq:dilute-E-mode-convergence}.

Finally, a direct calculation shows that $u_{E_2}(x):=u_{E_1}(R^\top x)$ is an eigenfunction with the same eigenvalue $\lambda_E(\eta,\tau)$, as claimed.
\end{proof}

We now study the asymptotics of the eigenvalue $\lambda_A(\eta,\tau)$ and the eigenfunction $u_A$ in the dilute high-contrast regime.

The constants span the zero eigenspace of \eqref{eq:dilute-periodic-problem} in $\mathcal V_A$. Consequently,
\begin{equation}\label{eq:A1-Rayleigh-dilute}
    \lambda_A(\eta,\tau)=\min_{\substack{0\ne u\in\mathcal V_A\\\int_Yb_{\eta,\tau}u\,dx=0}}
    \frac{\displaystyle\int_Y|\nabla u|^2\,dx}{\displaystyle\int_Yb_{\eta,\tau}u^2\,dx}.
\end{equation}
For an inclusion $D\Subset Y$, define the periodic capacity
\begin{equation}\label{eq:periodic-capacity}
    \operatorname{Cap}_Y(D):=\inf_{\substack{v\in H^1_\#(Y;\mathbb R)\\v=0\ {\rm on}\ D\\\int_Yv\,dx=1}}
    \int_Y|\nabla v|^2\,dx.
\end{equation}

\begin{lemma}\label{lem:dilute-A-convergence}
Suppose that Assumption~\ref{assump:neumanngap} holds. Then
\begin{equation}\label{eq:A1-branch-dilute}
    \lambda_A(\eta,\tau)=\frac{2\pi}{\ell_\eta}[1+o(1)] \qquad \mathrm{as\ } \eta\to0 
\end{equation}
uniformly in $\tau\in I_\eta$. This eigenvalue is simple in $\mathcal V_A$ for sufficiently
small $\eta$. Its $L^2_{b_{\eta,\tau}}(Y)$-normalized eigenfunction $u_A$, chosen with
positive cell average, satisfies
\begin{equation}\label{eq:dilute-A-mode-convergence}
    u_A\rightarrow1\quad\text{in }H^1(Y),\qquad \int_{D_\eta}b_{\eta,\tau}u_A^2\,dx\rightarrow0,
\end{equation}
uniformly on $I_\eta$.
\end{lemma}

\begin{proof}
The mean constraint gives coercivity by Poincar\'e's inequality and makes the energy strictly convex. Let $v_\eta$ be the unique minimizer in \eqref{eq:periodic-capacity} for $D_\eta$.
Uniqueness and symmetry imply that $v_\eta$ is $D_4$-invariant. We claim that, as $\eta \rightarrow 0$,
\begin{equation}\label{eq:dilute-capacity-asymptotic}
    \operatorname{Cap}_Y(D_\eta)=\frac{2\pi}{\ell_\eta}[1+o(1)],\qquad
    v_\eta\rightarrow1\quad\text{in }H^1(Y).
\end{equation}
For completeness, choose $r_-,r_+>0$ with $B_{r_-}(0)\subset B\subset B_{r_+}(0)$, and fix $R\in(0,1/2)$. Use the trial function equal to zero for $|x|\le\eta r_+$, equal to $\log(|x|/(\eta r_+))/\log(R/(\eta r_+))$ for $\eta r_+<|x|<R$, and equal to one elsewhere. Its energy is $2\pi/\log(R/(\eta r_+))$ and its cell mean tends to one by dominated convergence. Normalizing the mean gives the upper bound in \eqref{eq:dilute-capacity-asymptotic}. Poincar\'e's inequality now gives $\|v_\eta-1\|_{H^1(Y)}=O(\ell_\eta^{-1/2})$. Denote its circular average at radius $r$ by $\overline v_\eta(r)$. Integrating in $r$ gives a radius $R_\eta\in[1/4,1/3]$ with $\overline v_\eta(R_\eta)=1+O(\ell_\eta^{-1/2})$. Since $v_\eta$ has zero trace on the circle of radius $\eta r_-$, Cauchy--Schwarz in the radial variable yields
\[
 \int_Y|\nabla v_\eta|^2\,dx
 \geq \frac{2\pi|\overline v_\eta(R_\eta)|^2}{\log(R_\eta/(\eta r_-))}
 =\frac{2\pi}{\ell_\eta}[1+o(1)] .
\]
This proves the capacity asymptotic with the normalization used here. Write 
\begin{equation}
    K:=(\tau-1)|D_\eta|=(\tau-1)\eta^2 |B|, \qquad \widetilde{v}_{\eta,\tau}:=v_\eta - \frac{1}{1+K}.
\end{equation}
Since $\tau \in I_\eta$, $\tau \eta^2 \asymp_B \ell_\eta $. Hence, $K\rightarrow \infty$ as $\eta \rightarrow 0$.
Since $v_\eta=0$ in $D_\eta$, and $\int_Y v_\eta \,dx=1$, we get
\begin{equation}
    \int_Y b_{\eta, \tau} \widetilde{v}_{\eta,\tau} \,dx = 1-  \frac{\int_Y b_{\eta, \tau} \,dx}{1+K}=0.
\end{equation}
Moreover, since $v_\eta\rightarrow1$ in $H^1(Y)$, we have
\begin{equation}
    \int_Yb_{\eta,\tau} \widetilde{v}_{\eta,\tau} ^2\,dx=\int_Yv_\eta^2\,dx- \frac{1}{\int_Y b_{\eta, \tau}\,dx} =1+o(1) \qquad \mathrm{as\ } \eta\to0.
\end{equation}
Consequently,
\begin{equation}\label{upperboundoflambdaA}
    \lambda_A \leq \frac{\mathrm{Cap}_Y(D_\eta)}{ \int_Yb_{\eta,\tau} \widetilde{v}_{\eta,\tau} ^2\,dx } =\frac{2\pi}{\ell_\eta}[1+o(1)] \qquad \mathrm{as\ } \eta\to0.
\end{equation}

We now prove a lower bound for $\lambda_A$. Let $u$ be a minimizer in \eqref{eq:A1-Rayleigh-dilute}, normalized by
\begin{equation}\label{normalizbuA}
    \int_Y b_{\eta,\tau} u^2\,dx =1.
\end{equation}
Then
\begin{equation}
    \lambda_A = \int_Y|\nabla u|^2\,dx.
\end{equation}
Set
\begin{equation}
    m:=\int_Yu\,dx,\qquad
    a_B:=\frac1{|D_\eta|}\int_{D_\eta}u\,dx, \qquad q:=(\tau-1)\int_{D_\eta}|u-a_B|^2\,dx.
\end{equation}
The condition $\int_Yb_{\eta,\tau}u\,dx=0$ is equivalent to
\begin{equation}\label{normalizationmKaB}
    m+Ka_B=0.
\end{equation}
Since $|m|\le1$, we have $a_B=O_B(\ell_\eta^{-1})$. The Poincar\'{e} inequality and the normalization \eqref{normalizbuA} yield
\begin{equation}\label{eq:dilute-A-mass-splitting}
\begin{aligned}
    1 &= \int_Y u^2\,dx + (\tau-1) \int_{D_\eta} u^2\,dx \\
      &= m^2 + \int_Y (u-m)^2\,dx + q + K a_B^2 \\
      &= m^2+\int_Y|u-m|^2\,dx+q+\frac{m^2}{K} \\
      &= m^2+q+o(1),
\end{aligned}
\end{equation}
where we used \eqref{normalizationmKaB} in the third equality, and $|m|\leq 1$, $K\asymp_B \ell_\eta$, and $\lambda_A=O_B(\ell_\eta^{-1})$ in the last equality.

Applying \eqref{eq:dilute-Neumann-constant} to $w(y)=u(\eta y)-a_B$ on $B$ gives
\begin{equation}\label{eq:dilute-A-interior-bound}
    \ell_\eta\int_{D_\eta}|\nabla u|^2\,dx\ge\frac{\mu_A^{\mathrm{N}}(B)\ell_\eta}{(\tau-1)\eta^2}\,q.
\end{equation}

To estimate the exterior energy, choose $R_B>0$ with $B\Subset B_{R_B/2}(0)$, and let $\overline u(r)$ denote the average of $u$ on the circle $\{|x|=r\}$. Let
\begin{equation}
    W(x)=u(\eta x), \qquad x\in B_{R_B}(0).
\end{equation}
Let $c$ be the average of $W$ over $B_{R_B}(0)$. The Poincar\'e inequality gives
\begin{equation}
    \| W-c\|_{L^2(B_{R_B})} \leq C_B \| \nabla W\|_{L^2(B_{R_B})}.
\end{equation}
Also,
\begin{equation}
    |c-a_B| \leq |B|^{-1/2} \|W-c\|_{L^2(B)} \leq C_B \|\nabla W\|_{L^2(B_{R_B})}.
\end{equation}
Hence,
\begin{equation}
    \| W-a_B\|_{H^1(B_{R_B})} \leq C_B \| \nabla W\|_{L^2(B_{R_B})}.
\end{equation}
The trace theorem, followed by averaging on the boundary circle, yields
\begin{equation}\label{estimateubareta}
    \overline u(R_B\eta)=a_B+O_B(\lambda_A^{1/2})=o(1).
\end{equation}

Now, Jensen's inequality and the cell Poincar\'e inequality give
\begin{equation}
    \int_{1/4}^{1/3} 2\pi r|\overline{u}(r)-m|^2\,dr \leq \int_{\{1/4<|x|<1/3\} } |u-m|^2\,dx \leq C\lambda_A.
\end{equation}
Averaging in $r$ gives a radius $R_{\eta,\tau}\in[1/4,1/3]$ such that 
\begin{equation}\label{estimateubarRetatau}
    \overline u(R_{\eta,\tau})=m+O(\lambda_A^{1/2}) = m+ O_B(\ell_{\eta}^{-1/2}).
\end{equation}
Using polar coordinates, Jensen's inequality and the Cauchy--Schwarz inequality, we get
\begin{equation}\label{eq:dilute-A-exterior-bound}
\begin{aligned}
    \ell_\eta\int_{\{R_B\eta<|x|<R_{\eta,\tau}\}}|\nabla u|^2\,dx
    &\ge\frac{2\pi\ell_\eta|\overline u(R_{\eta,\tau})-\overline u(R_B\eta)|^2}{\log(R_{\eta,\tau}/(R_B\eta))}\\
    &=2\pi m^2+o(1),
\end{aligned}
\end{equation}
where we used \eqref{estimateubareta} and \eqref{estimateubarRetatau} in the last equality. All errors here are uniform on $\tau \in I_\eta$.

Set 
\begin{equation}
    \kappa=\frac{\mu_A^{\mathrm{N}}(B)}{(1+\rho)\Lambda_E(B)} >1.
\end{equation}
Since
\begin{equation}\label{intervalIetaimplication}
    \frac{(\tau-1)\eta^2}{\ell_\eta} \le(1+\rho) \frac{\Lambda_E(B)}{2\pi},
\end{equation}
adding the energy bounds on the two disjoint regions $D_\eta$ and ${\{R_B\eta<|x|<R_{\eta,\tau}\}}$ gives
\begin{equation}
    \ell_\eta \lambda_A\ge2\pi m^2+2\pi\kappa q+o(1)=2\pi[1+(\kappa-1)q]+o(1),
\end{equation}
where we used \eqref{eq:dilute-A-interior-bound}, and \eqref{eq:dilute-A-exterior-bound}--\eqref{intervalIetaimplication} in the inequality and \eqref{eq:dilute-A-mass-splitting} in the second equality. The upper bound \eqref{upperboundoflambdaA} forces $q\to0$, $m^2\to1$, and $\ell_\eta\lambda_A\to2\pi$. This proves \eqref{eq:A1-branch-dilute}.

Choose the sign of $u$ such that $m>0$. Then, $m\rightarrow 1$, and
\begin{equation}
    \| u -1 \|_{L^2(Y)}^2 = \| u -m \|_{L^2(Y)}^2 + (m-1)^2 \rightarrow 0.
\end{equation}
Also, $\| \nabla u \|_{L^2(Y)}^2 = \lambda_A\rightarrow 0$. Hence $u\rightarrow 1$ in $H^1(Y)$. Moreover,
\begin{equation}\label{vanishmassestimate}
    \int_{D_\eta}b_{\eta,\tau}u^2\,dx =\frac{\tau}{\tau-1}\left(q+\frac{m^2}{K}\right)\to0.
\end{equation}
This completes the proof of \eqref{eq:dilute-A-mode-convergence}.

Finally, if simplicity failed for arbitrarily small $\eta$, there would be sequences $\eta_j\to0$ and $\tau_j\in I_{\eta_j}$ for which we could choose two real, $b_{\eta_j,\tau_j}$-orthonormal eigenfunctions $u_j,w_j$ of $\lambda_A$ satisfying the above estimates. Choose their signs so that both cell averages are positive. Then
\begin{equation}
    u_j\to1,\qquad w_j\to1 \quad\text{in }L^2(Y),
\end{equation}
and consequently
\begin{equation}
    \int_Yu_jw_j\,dx\rightarrow1.
\end{equation}
The Cauchy--Schwarz inequality and the estimate \eqref{vanishmassestimate} give
\begin{equation}
    \left|(\tau_j-1)\int_{D_{\eta_j}}u_jw_j\,dx\right| \le\left(\int_{D_{\eta_j}} b_{\eta_j,\tau_j} u_j^2\,dx\right)^{1/2} \left(\int_{D_{\eta_j}} b_{\eta_j,\tau_j} w_j^2\,dx\right)^{1/2}\rightarrow0.
\end{equation}
Thus, their weighted inner product satisfies
\begin{equation}
    0=\int_Y b_{\eta_j,\tau_j} u_jw_j\,dx=\int_Yu_jw_j\,dx+(\tau_j-1)\int_{D_{\eta_j}}u_jw_j\,dx\rightarrow1,
\end{equation}
a contradiction. Since all estimates are uniform on $I_\eta$, this proves simplicity for all sufficiently small $\eta$, uniformly in $\tau\in I_\eta$.
\end{proof}

\subsection{Exact tuning and the threefold Dirac point}
\label{subsec:dilute-tuning}

We use the symmetry projections to separate the branches. We compute the signed coupling between the capacitary mode and a dipole mode, tune their eigenvalues by varying the contrast, and then exclude additional modes at the crossing. For each character $\chi$ of $D_4$, let $\mathcal V_\chi$ consist of the real periodic $H^1$ functions satisfying $u(Qx)=\chi(Q)u(x)$. Thus $\mathcal V_{A_1}=\mathcal V_A$. Denote by $\mathcal V_{E_2}$ the axial space of functions even in $x_1$ and odd in $x_2$.

\begin{lemma}\label{lem:dilute-symmetry-reduction}
Let $b\in L^\infty(Y)$ be positive almost everywhere and $D_4$-invariant. For $Q\in D_4$ set $(T_Qu)(x)=u(Q^\top x)$. If a finite-dimensional space $\mathcal M\subset H^1_\#(Y;\mathbb R)$ is invariant under all $T_Q$, the four character projections and the two axial projections from Lemma~\ref{rep-prelim:D4-six-projections} preserve its orthogonal complement $\mathcal M^{\perp_b}$ in $L_b^2(Y)$ and every periodic eigenspace of $-\Delta u=\lambda bu$. Every nonzero eigenfunction in $\mathcal M^{\perp_b}$ has a nonzero component in one of the six symmetry spaces with the same eigenvalue. In each space except $\mathcal V_A$, every function has zero cell average.
\end{lemma}
\begin{proof}
The transformations $T_Q$ are unitary on the real weighted space and preserve the Dirichlet energy. Hence, they preserve every periodic eigenspace. They preserve $\mathcal M^{\perp_b}$ because $T_Q^*=T_{Q^{-1}}$ and $\mathcal M$ is invariant. The projections are finite linear combinations of these transformations and sum to the identity, so at least one projected component is nonzero. In every nontrivial character space some symmetry changes the sign of the function. A coordinate reflection does so in each axial space. Integration over $Y$ then gives zero.
\end{proof}

\begin{lemma}\label{lem:uniform-dipole-coupling}
Under Assumption~\ref{assump:neumanngap}, let $u_A$ and $u_{E_1}$ have the sign and normalization of Lemmas~\ref{lem:dilute-A-convergence} and \ref{lem:dilute-E-convergence}. Define the signed coupling
\[
 \Gamma_{\eta,\tau}=\int_Y(u_A\partial_1u_{E_1}-u_{E_1}\partial_1u_A)\,dx .
\]
Then $\sqrt\tau\,\Gamma_{\eta,\tau}\to-2\pi d_B$ as $\eta\to0$, uniformly for $\tau\in I_\eta$.
\end{lemma}
\begin{proof}
Use the normalized
eigenfunctions $u_A$ and $u_{E_1}$ as in
Lemmas~\ref{lem:dilute-A-convergence} and
\ref{lem:dilute-E-convergence}, respectively, and write
\[
 T:=\tau-1,\qquad m_A:=\int_Yu_A(x)\,dx,\qquad
 q_{\eta,\tau}(y):=\eta\sqrt\tau\,u_{E_1}(\eta y),\qquad
 a_{\eta,\tau}(z):=u_A(\eta z).
\]
In particular, $m_A>0$. The two lemmas and the definition of $I_\eta$
give, uniformly for $\tau\in I_\eta$,
\[
\begin{aligned}
 \lambda_A&=O_B(\ell_\eta^{-1}),\qquad
 \lambda_E=O_B(\ell_\eta^{-1}),\qquad
 \tau\eta^2\asymp_B\ell_\eta,\\
 \lambda_A\ell_\eta&\longrightarrow2\pi,\qquad
 \lambda_ET\eta^2\longrightarrow\Lambda_E(B),\qquad
 m_A\longrightarrow1.
\end{aligned}
\]
Here and below $\lambda_A$ and $\lambda_E$ denote
$\lambda_A(\eta,\tau)$ and $\lambda_E(\eta,\tau)$. Moreover,
$q_{\eta,\tau}\to U_1$ in $L^2(B)$ uniformly on $I_\eta$,
as established in Lemma~\ref{lem:dilute-E-convergence}.

Periodic integration by parts gives
\[
 \Gamma_{\eta,\tau}=-2\int_Yu_{E_1}\partial_1u_A\,dx.
\]

\textit{The periodic resolvent representations.}
For $0\le\lambda\le2\pi^2$, define the mean-zero periodic kernel
\[
 G_\lambda(x):=
 \sum_{k\in\mathbb Z^2\setminus\{0\}}
 \frac{e^{2\pi\mathrm i k\cdot x}}
      {4\pi^2|k|^2-\lambda},
\]
where the series is understood distributionally. This is the inverse
of $-\Delta-\lambda$ on the nonconstant Fourier modes. In fact,
the eigenvalue equation can be written as
\[
 (-\Delta-\lambda_A)u_A
      =\lambda_AT\mathbf{1}_{D_\eta}u_A,
 \qquad
 (-\Delta-\lambda_E)u_{E_1}
      =\lambda_ET\mathbf{1}_{D_\eta}u_{E_1}.
\]
Taking their nonzero Fourier coefficients therefore yields
\[
 u_A-m_A=\lambda_ATG_{\lambda_A}*
                   (\mathbf{1}_{D_\eta}u_A),
 \qquad
 u_{E_1}=\lambda_ETG_{\lambda_E}*
                   (\mathbf{1}_{D_\eta}u_{E_1}).
\]
The second equality uses $\int_Yu_{E_1}=0$, which follows from
oddness in $x_1$. Convolution is on the unit torus $Y$.
For sufficiently small $\eta$, both eigenvalues lie in
$[0,2\pi^2]$, uniformly on $I_\eta$.

\textit{The local convolution-kernel expansion.}
Set
\[
 K:=G_{\lambda_E}*G_{\lambda_A},\qquad
 K_{0,0}:=G_0*G_0.
\]
Their nonzero Fourier coefficients are
\[
 \widehat K(k)=
 \frac{1}{(4\pi^2|k|^2-\lambda_E)
               (4\pi^2|k|^2-\lambda_A)},
 \qquad
 \widehat K_{0,0}(k)=\frac{1}{(4\pi^2|k|^2)^2}.
\]
The Fourier coefficients of $K$ are $O(|k|^{-4})$, so $K$ has a
continuous first derivative. Also,
\[
 \left|\widehat K(k)-\widehat K_{0,0}(k)\right|
      \le C_B|k|^{-6},\qquad k\ne0,
\]
uniformly for the eigenvalues under consideration. Since
$\sum_{k\ne0}|k|^{-4}<\infty$ in dimension two, the difference
$K-K_{0,0}$ has uniformly bounded second derivatives. All these
kernels are even, so $\nabla(K-K_{0,0})(0)=0$. Consequently,
\[
 |\nabla(K-K_{0,0})(x)|\le C_B|x|
\]
near the origin.

The periodic kernel $K_{0,0}$ satisfies
$\Delta^2K_{0,0}=\delta_0-1$. In dimension two,
\[
 \Delta^2\left(\frac{|x|^2\log|x|}{8\pi}\right)=\delta_0.
\]
It follows by local elliptic regularity that, in a neighborhood of
the origin,
\[
 K_{0,0}(x)=\frac{|x|^2\log|x|}{8\pi}+H(x),
\]
where $H$ is smooth and even. Hence $\nabla H(0)=0$, and
\[
 \partial_1K(x)
      =\frac{x_1\log|x|}{4\pi}+O_B(|x|).
\]
For $h=y-z$, $y,z\in B$, the function $h_1\log|h|$, extended by
zero at $h=0$, is bounded on the bounded set $B-B$. Therefore
\begin{equation}\label{eq:dilute-kernel-expansion}
 \partial_1K\bigl(\eta(y-z)\bigr)
      =\frac{\eta(y_1-z_1)\log\eta}{4\pi}+O_B(\eta),
      \qquad y,z\in B,
\end{equation}
where the error is uniform in both variables and in $\tau\in I_\eta$.

\textit{Reduction to the dipole moment.}
Substituting the two resolvent representations into the expression
for $\Gamma_{\eta,\tau}$ gives
\[
 \Gamma_{\eta,\tau}
   =-2\lambda_A\lambda_ET^2
     \int_{D_\eta}\int_{D_\eta}
       u_{E_1}(x)u_A(z)\partial_1K(x-z)\,dz\,dx.
\]
Indeed, the evenness of $G_{\lambda_E}$ and periodic integration by
parts give the convolution identity
\[
 \int_YG_{\lambda_E}(w-x)
         \partial_{w_1}G_{\lambda_A}(w-z)\,dw
        =\partial_1K(x-z).
\]
These identities may first be established with finite Fourier sums
and then passed to the limit; the resolvents map $L^2(Y)$ into
$H^2(Y)$, and $\partial_1K$ is continuous.

Rescaling both integration variables to $B$ gives the exact identity
\[
 \sqrt\tau\,\Gamma_{\eta,\tau}
    =-2\lambda_A\lambda_ET^2\eta^3
        \int_B\int_Bq_{\eta,\tau}(y)a_{\eta,\tau}(z)
          \partial_1K\bigl(\eta(y-z)\bigr)\,dz\,dy.
\]
There are two exact integral identities that determine the leading
coefficient. First, the oddness gives
\[
 \int_Bq_{\eta,\tau}(y)\,dy=0.
\]
Second, testing the equation for $u_A$ against the constant function
$1$ gives weighted orthogonality to the zero mode:
\[
 0=\int_Yb_{\eta,\tau}u_A
      =m_A+T\int_{D_\eta}u_A.
\]
Consequently,
\[
 \int_Ba_{\eta,\tau}(z)\,dz
       =-\frac{m_A}{T\eta^2}.
\]
It follows that
\[
 \begin{aligned}
 &\int_B\int_Bq_{\eta,\tau}(y)a_{\eta,\tau}(z)
                         (y_1-z_1)\,dz\,dy\\
 &\quad=
    \left(\int_By_1q_{\eta,\tau}(y)\,dy\right)
    \left(\int_Ba_{\eta,\tau}(z)\,dz\right)
   -\left(\int_Bq_{\eta,\tau}(y)\,dy\right)
    \left(\int_Bz_1a_{\eta,\tau}(z)\,dz\right)\\
 &\quad=-\frac{m_A}{T\eta^2}
              \int_By_1q_{\eta,\tau}(y)\,dy.
 \end{aligned}
\]
Using \eqref{eq:dilute-kernel-expansion}, the preceding exact
integral identity, and $\log\eta=-\ell_\eta$, we obtain
\begin{equation}\label{eq:gamma-dilute}
 \begin{aligned}
 \sqrt\tau\,\Gamma_{\eta,\tau}
    &=-\frac{\lambda_A\ell_\eta}{2\pi}
       (\lambda_ET\eta^2)m_A
          \int_By_1q_{\eta,\tau}(y)\,dy
        +O_B(\ell_\eta^{-1/2})\\
    &\longrightarrow-\Lambda_E(B)\int_By_1U_1(y)\,dy
       =-2\pi d_B.
 \end{aligned}
\end{equation}
The remainder from the kernel expansion is controlled by weighted normalization, which gives
\[
 \|q_{\eta,\tau}\|_{L^2(B)}\le1,
 \qquad
 \|a_{\eta,\tau}\|_{L^2(B)}
          \le\frac{1}{\eta\sqrt\tau}.
\]
Thus, the contribution of the $O_B(\eta)$ kernel error to
$\sqrt\tau\,\Gamma_{\eta,\tau}$ is bounded by
\[
 \begin{aligned}
 C_B\lambda_A\lambda_ET^2\eta^4
       \|q_{\eta,\tau}\|_{L^2(B)}
       \|a_{\eta,\tau}\|_{L^2(B)}
 &\le C_B\frac{\lambda_A\lambda_ET^2\eta^4}
                     {\eta\sqrt\tau}\\
 &=O_B(\ell_\eta^{-1/2}),
 \end{aligned}
\]
because $\lambda_A\lambda_E=O_B(\ell_\eta^{-2})$,
$T^2\eta^4=O_B(\ell_\eta^2)$, and
$\eta\sqrt\tau\asymp_B\sqrt{\ell_\eta}$. The $O_B(\ell_\eta^{-1/2})$ estimate controls the error in the convolution kernel. The remaining factors in \eqref{eq:gamma-dilute} converge uniformly on $I_\eta$, giving $\sqrt\tau\,\Gamma_{\eta,\tau}=-2\pi d_B+o(1)$. Since $d_B>0$, it also shows
$\Gamma_{\eta,\tau}<0$ for all sufficiently small $\eta$,
uniformly on $I_\eta$.

\end{proof}

\begin{theorem}\label{thm:dilute-dirac}
Suppose that Assumption~\ref{assump:neumanngap} holds. There exist $\eta_0>0$ and $c_B>0$, depending only on the fixed shape $B$ and the chosen $\rho$, such that, for every $0<\eta<\eta_0$, there is a unique contrast $\tau_*(\eta)\in I_\eta$ for which the branches coincide. Their common value
\begin{equation}\label{defoflambdastareta}
    \lambda_*(\eta):=\lambda_A(\eta,\tau_*(\eta))=\lambda_E(\eta,\tau_*(\eta))
\end{equation}
has multiplicity exactly three and occupies bands two, three and four at $\Gamma$. The crossing is transverse, with
\begin{equation}\label{eq:tuning-transverse}
 \left.\partial_\tau(\lambda_E-\lambda_A)\right|_{\tau=\tau_*(\eta)}
   =-\frac{\lambda_*(\eta)}{\tau_*(\eta)}[1+o(1)]
   \qquad(\eta\to0).
\end{equation}
If $\sigma_\eta$ is the spectrum of the eigenvalue problem \eqref{eq:dilute-periodic-problem} for these tuned coefficients, then
\begin{equation}\label{eq:dilute-uniform-relative-gap}
    \operatorname{dist}\bigl(\lambda_*(\eta),\sigma_\eta\setminus\{\lambda_*(\eta)\}
    \bigr)\ge c_B\lambda_*(\eta).
\end{equation}
The eigenspace has a real basis $\{f_0,f_1,f_2\}$ satisfying
\eqref{b-orthogonality}, for which
\begin{equation}\label{eq:dilute-coupling-matrices}
    C_1=\gamma_\eta \begin{pmatrix}0&1&0\\-1&0&0\\0&0&0\end{pmatrix},
    \qquad C_2=\gamma_\eta \begin{pmatrix}0&0&1\\0&0&0\\-1&0&0\end{pmatrix}, \qquad \gamma_\eta>0.
\end{equation}
As $\eta\to0$,
\begin{equation}\label{eq:dilute-tuned-asymptotics}
\begin{aligned}
    \tau_*(\eta)&=\frac{\Lambda_E(B)}{2\pi\eta^2}\ell_\eta[1+o(1)],\\
    \lambda_*(\eta)&=\frac{2\pi}{\ell_\eta}[1+o(1)],\\
    \gamma_\eta&=\frac{(2\pi)^{3/2}d_B}{\sqrt{\Lambda_E(B)}}\frac{\eta}{\sqrt{\ell_\eta}}[1+o(1)].
\end{aligned}
\end{equation}
In particular, $G=\gamma_\eta^2I_2>0$, and Proposition~\ref{prop:generalcriteria} gives
\begin{equation}\label{eq:dilute-cones}
    \lambda_\pm(\alpha)=\lambda_*(\eta)\pm\gamma_\eta|\alpha| +O_\eta(|\alpha|^2),\qquad
    \lambda_c(\alpha)=\lambda_*(\eta)+O_\eta(|\alpha|^2).
\end{equation}
The two conical frequency branches have common speed
\begin{equation}\label{eq:vD-dilute}
    \frac{\gamma_\eta}{2\omega_*(\eta)} =\frac{\pi d_B}{\sqrt{\Lambda_E(B)}}\eta[1+o(1)],
    \qquad \omega_*(\eta):=\sqrt{\lambda_*(\eta)}.
\end{equation}
\end{theorem}

\begin{proof}[Proof of Theorem~\ref{thm:dilute-dirac}]

\medskip
\noindent\textit{Step 1: exact tuning.}
For fixed $\eta$, apply min--max on the fixed symmetry spaces, with the second ordered eigenvalue in $\mathcal V_A$ and the first in $\mathcal V_{E_1}$. Since the positive weights vary continuously in the sense that $\|b_{\eta,\tau'}/b_{\eta,\tau}-1\|_\infty\to0$ as $\tau'\to\tau$, this gives continuity of
$\lambda_A(\eta,\tau)$ and $\lambda_E(\eta,\tau)$ in $\tau$.
Lemmas~\ref{lem:dilute-E-convergence} and
\ref{lem:dilute-A-convergence} give, uniformly for $\tau\in I_\eta$,
\begin{equation}\label{expansionoflambdaEA}
    \lambda_E(\eta,\tau)
       =\frac{\Lambda_E(B)}{\tau\eta^2}[1+o(1)],
    \qquad
    \lambda_A(\eta,\tau)
       =\frac{2\pi}{\ell_\eta}[1+o(1)].
\end{equation}
At $\tau_\pm=(1\pm\rho)\tau_0(\eta)$, these formulas imply
\[
    \frac{\lambda_E(\eta,\tau_\pm)}
         {\lambda_A(\eta,\tau_\pm)}
       \longrightarrow\frac{1}{1\pm\rho}.
\]
Since $0<\rho<1$, for all sufficiently small $\eta$, we have
\[
    \lambda_E(\eta,\tau_-)>
         \lambda_A(\eta,\tau_-),
    \qquad
    \lambda_E(\eta,\tau_+)<
         \lambda_A(\eta,\tau_+).
\]
The intermediate value theorem gives
$\tau_*(\eta)\in(\tau_-,\tau_+)$ such that
\[
    \lambda_A(\eta,\tau_*(\eta))
       =\lambda_E(\eta,\tau_*(\eta))=:\lambda_*(\eta).
\]
We show that the crossing in $I_\eta$ is unique. By Lemmas~\ref{lem:dilute-E-convergence} and \ref{lem:dilute-A-convergence}, each branch is simple in its own symmetry space for all sufficiently small $\eta$, including at the common eigenvalue. For fixed $\eta$ the mass form is affine in $\tau$. Analytic perturbation of a simple eigenvalue in each fixed sector therefore gives locally differentiable branches and normalized eigenfunctions \cite[Ch.~VII]{Kato1995}. Equivalently, one may reduce to the compact operator defined on the form domain by the mass form relative to the positive form $\int_Y|\nabla u|^2+\int_Y b_{\eta,\tau}|u|^2$. Its nonzero eigenvalues are $(1+\lambda)^{-1}$ and its dependence on $\tau$ is analytic near each real parameter.

For $X=A,E$, use $u_X=u_A,u_{E_1}$, respectively, with $\int_Y b_{\eta,\tau}u_X^2=1$, and put $P_X(\eta,\tau)=\tau\int_{D_\eta}u_X^2$. Differentiate the variational equation in $\tau$ and test against $u_X$. The terms containing $\partial_\tau u_X$ cancel using the original eigenvalue equation, leaving $\partial_\tau\log\lambda_X=-\int_{D_\eta}u_X^2$. By \eqref{eq:dilute-E-mode-convergence}, $P_E=\int_Bq_{\eta,\tau}^2\to1$. By \eqref{eq:dilute-A-mode-convergence}, $P_A\to0$. Both limits are uniform over $I_\eta$. Hence,
\[
 \tau\,\partial_\tau\log(\lambda_E/\lambda_A)=P_A-P_E=-1+o(1)
 \quad\text{uniformly on }I_\eta .
\]
For sufficiently small $\eta$, the ratio is strictly decreasing throughout $I_\eta$. Its endpoint values straddle one, so its crossing is unique. Evaluation at $\tau_*(\eta)$ gives \eqref{eq:tuning-transverse}.

The normalized functions $u_A,u_{E_1},u_{E_2}$ at this contrast
are mutually orthogonal in $L^2_{b_{\eta,\tau_*}}(Y)$ by their
reflection parities. Thus the multiplicity is at least three.
Equating the uniform asymptotics in \eqref{expansionoflambdaEA}
gives the first two formulas in \eqref{eq:dilute-tuned-asymptotics}.
In particular,
\begin{equation}\label{eq:dilute-crossing-scaled-limit}
    \tau_*(\eta)\eta^2\lambda_*(\eta)
       \longrightarrow\Lambda_E(B).
\end{equation}

\medskip
\noindent\textit{Step 2: exact multiplicity and spectral isolation.}
Define the second variational threshold in the first axial parity
space by
\begin{equation}\label{defofLambdaE2}
    \Lambda_{E,2}(B):=
    \inf_{\substack{U\in\mathcal V_{E_1}^{\infty}\\
                    \int_BUU_1\,dx=0,\ \int_BU^2\,dx>0}}
       \frac{\int_{\mathbb R^2}|\nabla U|^2\,dx}
            {\int_BU^2\,dx}.
\end{equation}
The angular Hardy inequality \eqref{eq:dilute-angular-Hardy}
and local compactness give attainment, exactly as in
Lemma~\ref{lem:dilute-E-convergence}.
If $\Lambda_{E,2}(B)=\Lambda_E(B)$, a normalized minimizer
would also minimize the first quotient. Its simplicity would
force $U=\pm U_1$, contradicting $\int_BUU_1=0$.
Consequently,
\[
    \Lambda_{E,2}(B)>\Lambda_E(B).
\]

For each nontrivial one-dimensional character $\chi$, set
\[
    \mathcal V_\chi^{\infty}:=
    \{U\in H^1_{\mathrm{loc}}(\mathbb R^2;\mathbb R):
      \nabla U\in L^2(\mathbb R^2),\
      U(Qy)=\chi(Q)U(y)\text{ for all }Q\in D_4\},
\]
and define
\[
    \Lambda_\chi(B):=
    \inf_{\substack{U\in\mathcal V_\chi^{\infty}\\
                    \int_BU^2\,dx>0}}
       \frac{\int_{\mathbb R^2}|\nabla U|^2\,dx}
            {\int_BU^2\,dx}.
\]
Such a character is odd under some reflection. Its functions
therefore have zero angular average for almost every radius about the
origin, and satisfy the same angular Hardy inequality.
The direct method again gives attainment.
We next prove that
\begin{equation}\label{eq:dilute-character-strict-gap}
    \Lambda_\chi(B)>\Lambda_E(B)
    \qquad(\chi\ne1).
\end{equation}

For a reflection $Q\in D_4$, let $\mathcal H_Q^{\infty}$ be
the finite-energy space of functions satisfying $U(Qy)=-U(y)$,
and let $\nu_Q$ be its first Rayleigh quotient with denominator
$\int_BU^2$.
The positive half-plane Green-kernel argument in
Lemma~\ref{lem:dilute-E-convergence}, applied to the half-plane
bounded by the reflection axis of $Q$, shows that $\nu_Q$ has
a simple normalized minimizer $P_Q$, positive on that half-plane.
The perpendicular reflection $Q'=-Q$ preserves both the half-plane
and $B$. Simplicity and positivity imply
$P_Q(Q'y)=P_Q(y)$, and hence $P_Q(-y)=-P_Q(y)$.

Define $\mathcal V_{E_2}^{\infty}=\{U(R^\top\cdot):U\in\mathcal V_{E_1}^{\infty}\}$. The centrally odd finite-energy space is the energy- and
$L^2(B)$-orthogonal sum of the two axial parity spaces.
Their first values are both $\Lambda_E(B)$, with simple
normalized minimizers $U_1$ and
$U_2(y):=U_1(R^\top y)$.
Thus, every centrally odd function has Rayleigh quotient at least
$\Lambda_E(B)$. Applying this bound to the centrally odd minimizer $P_Q$ gives $\nu_Q\ge\Lambda_E(B)$.
Moreover, the space $\operatorname{span}\{U_1,U_2\}$ transforms
according to the standard representation of $D_4$:
\[
    \begin{pmatrix}U_1(Qy)\\U_2(Qy)\end{pmatrix}
       =Q\begin{pmatrix}U_1(y)\\U_2(y)\end{pmatrix}.
\]
Every reflection has a $-1$ eigenvector in this representation,
so this space contains a nonzero $Q$-odd function with quotient
$\Lambda_E(B)$. Therefore $\nu_Q=\Lambda_E(B)$.

For a nontrivial character $\chi$, choose a reflection $Q$ with
$\chi(Q)=-1$. Since $\chi(-I_2)=\chi(R)^2=1$, we also have
$\chi(Q')=-1$.
Its normalized minimizer belongs to $\mathcal H_Q^{\infty}$,
so $\Lambda_\chi(B)\ge\nu_Q=\Lambda_E(B)$.
Equality would force this minimizer to be proportional to $P_Q$.
This is impossible because it is odd under $Q'$, while $P_Q$
is even under $Q'$. This proves
\eqref{eq:dilute-character-strict-gap}.

By Assumption~\ref{assump:neumanngap}, the number
\begin{equation}\label{eq:dilute-fixed-shape-gap}
    \kappa_B:=\min\left\{
       \frac{\mu_A^{\mathrm N}(B)}{\Lambda_E(B)}-1,\
       \frac{\Lambda_{E,2}(B)}{\Lambda_E(B)}-1,\
       \min_{\chi\ne1}
          \left(\frac{\Lambda_\chi(B)}{\Lambda_E(B)}-1\right)
       \right\}
\end{equation}
is strictly positive.

At the tuned contrast, write
\[
    b_\eta:=b_{\eta,\tau_*(\eta)},\qquad
    \mathcal M_\eta:=
       \operatorname{span}\{1,u_A,u_{E_1},u_{E_2}\},\qquad
    \mathcal W_\eta:=\mathcal M_\eta^{\perp_{b_\eta}}.
\]
We claim that, for all sufficiently small $\eta$, the spectrum
on $\mathcal W_\eta$ is contained in
\begin{equation}\label{eq:dilute-complement-threshold}
    [(1+\kappa_B/2)\lambda_*(\eta),\infty).
\end{equation}
Suppose otherwise. Then there exist $\eta_j\to0$ and real
eigenfunctions $v_j\in\mathcal W_{\eta_j}$ with eigenvalues
$\lambda_j<(1+\kappa_B/2)\lambda_*(\eta_j)$, normalized by
\[
    \int_Yb_jv_j^2\,dx=1,
    \qquad b_j:=b_{\eta_j,\tau_j},\quad
    \tau_j:=\tau_*(\eta_j),\quad T_j:=\tau_j-1.
\]
The space $\mathcal M_\eta$ is $D_4$-invariant: the first two basis functions are invariant, while the last two are exchanged,
up to signs, by rotations and reflections.
Lemma~\ref{lem:dilute-symmetry-reduction} therefore allows us
to replace each $v_j$ by a normalized nonzero symmetry component
with the same eigenvalue, still in $\mathcal W_{\eta_j}$.
After passing to a subsequence, all $v_j$ lie in one fixed
character space or one fixed axial parity space.
The eigenvalue equation and \eqref{eq:dilute-crossing-scaled-limit}
give
\begin{equation}\label{eq:dilute-contradiction-energy}
    \int_Y|\nabla v_j|^2\,dx=\lambda_j
       =O_B(\ell_{\eta_j}^{-1}),\qquad
    \limsup_{j\to\infty}\tau_j\eta_j^2\lambda_j
       \le(1+\kappa_B/2)\Lambda_E(B).
\end{equation}

\smallskip
\noindent\textit{Case 1: a nontrivial character or an axial parity space.}
Lemma~\ref{lem:dilute-symmetry-reduction} gives
$\int_Yv_j\,dx=0$. Hence,
\[
    \int_Yv_j^2\,dx
       \le\frac1{4\pi^2}\int_Y|\nabla v_j|^2\,dx
       =o(1).
\]
Set $Y_j:=\eta_j^{-1}Y$ and
$V_j(y):=\eta_j\sqrt{\tau_j}\,v_j(\eta_jy)$.
Normalization and a change of variables give
\[
    \int_BV_j^2\,dy
       =\frac{\tau_j}{T_j}\left(1-\int_Yv_j^2\,dx\right)
       =1+o(1),\qquad
    \int_{Y_j}|\nabla V_j|^2\,dy
       =\tau_j\eta_j^2\lambda_j=O_B(1).
\]
For every fixed ball centered at the origin, the reflection
oddness gives a zero average there. Therefore, using Poincar\'e's inequality bounds
$V_j$ in local $H^1$ uniformly in $j$.
After taking a subsequence along which
$\tau_j\eta_j^2\lambda_j$ converges to its liminf,
local compactness and a diagonal extraction give
\[
    V_j\rightharpoonup V\text{ in }H^1_{\mathrm{loc}}(\mathbb R^2),
    \qquad V_j\to V\text{ in }L^2_{\mathrm{loc}}(\mathbb R^2),
\]
with the same symmetry, $\int_BV^2\,dy=1$, and
\[
    \int_{\mathbb R^2}|\nabla V|^2\,dy
       \le\liminf_{j\to\infty}\tau_j\eta_j^2\lambda_j.
\]
If the selected space is $\mathcal V_\chi$, $\chi\ne1$, this yields
\[
    \liminf_{j\to\infty}\tau_j\eta_j^2\lambda_j
       \ge\Lambda_\chi(B)
       \ge(1+\kappa_B)\Lambda_E(B),
\]
contradicting \eqref{eq:dilute-contradiction-energy}.

If the selected space is $\mathcal V_{E_1}$, let
$q_j(y):=\eta_j\sqrt{\tau_j}\,u_{E_1}(\eta_jy)$, with the mode
evaluated at $(\eta_j,\tau_j)$.
Weighted orthogonality gives the exact identity
\[
    \int_BV_jq_j\,dy
       =-\int_{Y\setminus D_{\eta_j}}v_ju_{E_1}\,dx=o(1),
\]
where Cauchy--Schwarz, $\|v_j\|_{L^2(Y)}\to0$, and
$\|u_{E_1}\|_{L^2(Y)}\le1$ prove the last assertion.
Since $q_j\to U_1$ in $L^2(B)$, it follows that
$\int_BVU_1\,dy=0$. Thus,
\[
    \liminf_{j\to\infty}\tau_j\eta_j^2\lambda_j
       \ge\Lambda_{E,2}(B)
       \ge(1+\kappa_B)\Lambda_E(B),
\]
again a contradiction. If the selected space is
$\mathcal V_{E_2}$, rotation by $R^\top$ reduces it to the
$E_1$ case; the energy, inclusion norm, and orthogonality to
the first mode are preserved.

\smallskip
\noindent\textit{Case 2: the fully invariant space.}
Let $u_{A,j}$ denote the normalized $A$-mode at $(\eta_j,\tau_j)$,
and put $m_j:=\int_Yv_j\,dx$.
Weighted orthogonality to $u_{A,j}$ gives
\[
    m_j=\int_Yv_j(1-u_{A,j})\,dx
          -T_j\int_{D_{\eta_j}}v_ju_{A,j}\,dx.
\]
Since $b_j\ge1$, normalization implies
$\|v_j\|_{L^2(Y)}\le1$ and $T_j\int_{D_{\eta_j}}v_j^2\le1$.
Consequently,
\[
    |m_j|\le\|1-u_{A,j}\|_{L^2(Y)}
       +\left(T_j\int_{D_{\eta_j}}u_{A,j}^2\,dx\right)^{1/2}
       =o(1)
\]
by \eqref{eq:dilute-A-mode-convergence}. Therefore,
\[
    \int_Yv_j^2\,dx
       =m_j^2+\int_Y|v_j-m_j|^2\,dx
       \le m_j^2+\frac{\lambda_j}{4\pi^2}=o(1).
\]
Define
\[
    a_j:=\frac1{|D_{\eta_j}|}\int_{D_{\eta_j}}v_j\,dx,
    \qquad K_j:=T_j|D_{\eta_j}|\asymp_B\ell_{\eta_j},
    \qquad p_j:=T_j\int_{D_{\eta_j}}|v_j-a_j|^2\,dx.
\]
Orthogonality to $1$ yields $m_j+K_ja_j=0$, and then normalization
gives the exact mass decomposition
\[
    1=\int_Yv_j^2\,dx+p_j+K_ja_j^2
       =\int_Yv_j^2\,dx+p_j+\frac{m_j^2}{K_j}.
\]
Hence $p_j=1+o(1)$.
The function $w_j(y):=v_j(\eta_jy)-a_j$ is invariant and has
zero average on $B$. Applying the invariant Neumann inequality
and changing variables gives
\[
\begin{aligned}
    \lambda_j
       &\ge\int_{D_{\eta_j}}|\nabla v_j|^2\,dx
        =\int_B|\nabla w_j|^2\,dy\\
       &\ge\mu_A^{\mathrm N}(B)\int_Bw_j^2\,dy
        =\frac{\mu_A^{\mathrm N}(B)}{T_j\eta_j^2}\,p_j.
\end{aligned}
\]
It follows that
\[
    \liminf_{j\to\infty}\tau_j\eta_j^2\lambda_j
       \ge\mu_A^{\mathrm N}(B)
       \ge(1+\kappa_B)\Lambda_E(B),
\]
contradicting \eqref{eq:dilute-contradiction-energy}.
This proves \eqref{eq:dilute-complement-threshold}.

An additional eigenfunction at $\lambda_*(\eta)$ could be chosen
orthogonal to the three constructed modes and is automatically
orthogonal to the constant zero mode. It would therefore belong
to $\mathcal W_\eta$, contradicting
\eqref{eq:dilute-complement-threshold}.
Thus, the multiplicity is exactly three. Since the zero eigenvalue is simple and \eqref{eq:dilute-complement-threshold} places every other eigenvalue above the triple, it occupies precisely bands two through four at $\Gamma$.
The distance to zero is $\lambda_*(\eta)$, and the distance to the spectrum on $\mathcal W_\eta$ is at least $\kappa_B\lambda_*(\eta)/2$. This proves \eqref{eq:dilute-uniform-relative-gap} with
$c_B=\min\{1,\kappa_B/2\}$, after decreasing $\eta_0$.

\medskip
\noindent\textit{Step 3. The coupling matrices and conical branches.} By Lemma~\ref{lem:uniform-dipole-coupling}, the signed coupling is strictly negative for sufficiently small $\eta$.
At the crossing choose
\[
 f_0:=u_A,\qquad f_1:=-u_{E_1},\qquad f_2:=-u_{E_2},
 \qquad u_{E_2}(x):=u_{E_1}(R^\top x).
\]
The parity of $f_0$ is even in both variables; the parities of
$f_1$ and $f_2$ are, respectively, odd--even and even--odd.
Because $b_{\eta,\tau_*}$ is even in both variables, these parities
give all three weighted cross inner products equal to zero.
Rotation preserves the weighted norm of $u_{E_1}$, so the three
functions satisfy \eqref{b-orthogonality} and form a basis of the
three-dimensional eigenspace established in Step~2.

Since $A=I_2$, the definition \eqref{defofCj} becomes
\[
 (C_j)_{rs}=\int_Y
       (f_r\partial_jf_s-f_s\partial_jf_r)\,dx.
\]
Parity gives
\[
 (C_1)_{02}=(C_1)_{12}=(C_2)_{01}=(C_2)_{12}=0.
\]
The only possibly nonzero upper-triangular entries are thus
$(C_1)_{01}$ and $(C_2)_{02}$. Moreover,
\[
 f_0(Rx)=f_0(x),\qquad f_2(Rx)=f_1(x).
\]
Differentiation with respect to $x_1$ gives
\[
 (\partial_2f_0)(Rx)=\partial_1f_0(x),\qquad
 (\partial_2f_2)(Rx)=\partial_1f_1(x).
\]
Changing the variables by $R$ in the integral for $(C_2)_{02}$ therefore
proves
\[
 (C_2)_{02}=(C_1)_{01}
      =-\Gamma_{\eta,\tau_*(\eta)}=:\gamma_\eta>0.
\]
Since each $C_j$ is skew-symmetric, these identities give exactly
\eqref{eq:dilute-coupling-matrices}. In particular,
$G=\gamma_\eta^2I_2>0$.

By \eqref{eq:gamma-dilute} and the already proved asymptotic for
$\tau_*(\eta)$,
\[
 \begin{aligned}
 \gamma_\eta
   &=\frac{2\pi d_B}{\sqrt{\tau_*(\eta)}}[1+o(1)]\\
   &=\frac{(2\pi)^{3/2}d_B}{\sqrt{\Lambda_E(B)}}
        \frac{\eta}{\sqrt{\ell_\eta}}[1+o(1)].
 \end{aligned}
\]
This is the third formula in \eqref{eq:dilute-tuned-asymptotics}.
For each fixed sufficiently small $\eta$, Proposition~
\ref{prop:generalcriteria} applies to the isolated triple
eigenvalue, and $\sqrt{\alpha^\top G\alpha}=\gamma_\eta|\alpha|$.
It therefore gives \eqref{eq:dilute-cones}. Finally, taking the square
root of the two conical eigenvalue branches shows that their linear
frequency slopes have magnitude $\gamma_\eta/(2\omega_*(\eta))$.
Using $\omega_*(\eta)=\sqrt{2\pi/\ell_\eta}\,[1+o(1)]$, we find
\[
 \frac{\gamma_\eta}{2\omega_*(\eta)}
    =\frac{\pi d_B}{\sqrt{\Lambda_E(B)}}\eta[1+o(1)],
\]
which proves \eqref{eq:vD-dilute} and completes the proof.

\end{proof}

The choice $A=I_2$ in \eqref{eq:dilute-weight} identifies the scalar problem with the axial electric polarization of a dielectric crystal with microscopic permittivity $b_{\eta,\tau}$ and unit permeability. The capacitary singlet has nonzero mean, so the forced cell problem determines both bulk parameters at the tuned contrast.

\begin{corollary}\label{cor:dilute-double-zero}
Under the assumptions of Theorem~\ref{thm:dilute-dirac}, let $\varepsilon_{\mathrm{eff},\eta}$ and $\mu_{\mathrm{eff},\eta}$ be the parameters in \eqref{eq:dz-definition} at $\tau=\tau_*(\eta)$, and let $Z_{*,\eta}$ be the limiting bulk impedance. For each sufficiently small $\eta$ both parameters vanish simply at $\omega_*(\eta)$. Writing $\omega_*=\omega_*(\eta)$, we have, as $\eta\to0$,
\begin{equation}\label{eq:dilute-double-zero}
 \omega_*\varepsilon_{\mathrm{eff},\eta}'(\omega_*)\sim 2,\qquad
 \omega_*\mu_{\mathrm{eff},\eta}'(\omega_*)\sim
       \frac{\Lambda_E(B)}{2\pi^2 d_B^2\eta^2},\qquad
 Z_{*,\eta}\sim\frac{\sqrt{\Lambda_E(B)}}{2\pi d_B\eta}.
\end{equation}
The derivatives are taken with respect to frequency at fixed $\eta$.
\end{corollary}
\begin{proof}
At the tuned contrast, \eqref{eq:dilute-A-mode-convergence} gives $m_*^2=(\int_Y u_A\,dx)^2\to1$. Substitute this limit and the cone speed \eqref{eq:vD-dilute} into \eqref{eq:dz-parameters}.
\end{proof}

\subsection{Effective dynamics near the tuned Dirac points}\label{subsec:dilute-dynamics}

Fix the crossing given by Theorem~\ref{thm:dilute-dirac}, and write
$b_\eta:=b_{\eta,\tau_*(\eta)}$.
Let $L_\eta=-b_\eta^{-1}\Delta$ be the nonnegative self-adjoint
operator in $L^2_{b_\eta}(\mathbb R^2)$ associated with the
form $\int_{\mathbb R^2}|\nabla u|^2\,dx$.
In the basis of \eqref{eq:dilute-coupling-matrices}, the effective
symbol from Section~\ref{sec:effective} is
\begin{equation}\label{eq:dilute-effective-symbol}
    \mathsf D_\eta(\xi)
       =-\frac{\mathrm i\gamma_\eta}{2\omega_*(\eta)}
        \begin{pmatrix}
        0&\xi_1&\xi_2\\-\xi_1&0&0\\-\xi_2&0&0
        \end{pmatrix}.
\end{equation}
Given $\Env_0\in\mathcal S(\mathbb R^2;\mathbb C^3)$, let $\Env$
solve \eqref{effectiveequationofE} with this symbol, equivalently
$\widehat{\Env}(\xi,T)=e^{-\mathrm iT\mathsf D_\eta(\xi)}
\widehat{\Env_0}(\xi)$ as in \eqref{evolutionequationE}.
Using the wave-packet notation of Section~\ref{sec:effective}, set
\[
    \mathrm{WP}_\delta^{\Env}(x,t)
       :=\delta\sum_{r=0}^2f_r(x)E^{(r)}(\delta x,\delta t),
    \qquad
    \mathrm{WP}_\delta^{\Env_0}:=\mathrm{WP}_\delta^{\Env}(\cdot,0),
\]
and define
\[
    e_{\delta,\eta}(\cdot,t):=
    e^{-\mathrm it\sqrt{L_\eta}}\mathrm{WP}_\delta^{\Env_0}
      -e^{-\mathrm it\omega_*(\eta)}
                    \mathrm{WP}_\delta^{\Env}(\cdot,t).
\]

\begin{corollary}\label{cor:dilute-dynamics}
Under \eqref{eq:A1-Neumann-gap}, fix $\rho$ in terms of $B$ as in \eqref{eq:critical-contrast-dilute}. There are constants $C_B>0$
and $C_{B,N}>0$, depending only on the fixed shape $B$
(and on $N$ for the latter), such that, for sufficiently small
$\eta$, every $\delta>0$, $t\ge0$, and integer $N\ge2$,
\begin{equation}\label{eq:dilute-dynamics-error}
\begin{aligned}
    \|e_{\delta,\eta}(\cdot,t)\|_{L^2_{b_\eta}(\mathbb R^2)}
    &\le C_B\delta\sqrt{\ell_\eta}\,\|\Env_0\|_{\dot H^1}
       +C_B\delta^2t\sqrt{\ell_\eta}\,\|\Env_0\|_{\dot H^2}\\
    &\quad+C_{B,N}(\delta\sqrt{\ell_\eta})^N
                                      \|\Env_0\|_{\dot H^N}.
\end{aligned}
\end{equation}
In particular, for fixed $\Env_0$, the error tends to zero uniformly
on $0\le t\le t_\delta$ in any joint limit for which
\begin{equation}\label{eq:dilute-valid-times}
    \delta\sqrt{\ell_\eta}\longrightarrow0,\qquad
    \delta^2t_\delta\sqrt{\ell_\eta}\longrightarrow0.
\end{equation}
\end{corollary}

\begin{proof}
Let $\Delta_\eta$ be the distance on the left of \eqref{eq:dilute-uniform-relative-gap}. The quantities in \eqref{coefficientbound} and Lemma~\ref{componentk0est} become
\[
    s_\eta^2=
    \sup_{0\ne v\in H^1_\#(Y)}
       \frac{\int_Y|v|^2\,dx}
       {\lambda_*^{-1}(\eta)\int_Y|\nabla v|^2\,dx+
                             \int_Yb_\eta|v|^2\,dx},
    \qquad
    d_\eta=\min\{1,\Delta_\eta/\lambda_*(\eta)\}.
\]

Since $b_\eta\ge1$, we have $s_\eta\le1$.
Testing with the normalized $u_A$ gives
$s_\eta^2\ge\frac12\int_Yu_A^2\,dx=\frac12+o(1)$.
Theorem~\ref{thm:dilute-dirac} gives $c_B\le d_\eta\le1$
and $\omega_*(\eta)\asymp\ell_\eta^{-1/2}$.
Hence the radius in \eqref{radiusassump} satisfies
\[
    r_\eta=\min\{1,c\omega_*(\eta)d_\eta/s_\eta\}
          \asymp_B\ell_\eta^{-1/2}.
\]
Substitution in \eqref{eq:effective-dynamics-error} proves
\eqref{eq:dilute-dynamics-error} and
\eqref{eq:dilute-valid-times}. In particular, the increasing
upper bound $\|b_\eta\|_\infty=\tau_*(\eta)$ causes no additional
loss beyond the parameters displayed above.
\end{proof}

\begin{remark}
The speed in \eqref{eq:vD-dilute} is of order $\eta$.
For a propagating component, an order-one displacement in the envelope variable $X=\delta x$
therefore occurs on the physical time scale
$t\asymp(\delta\eta)^{-1}$. The sufficient conditions
\eqref{eq:dilute-valid-times} cover this scale if
$\delta\sqrt{\ell_\eta}=o(\eta)$.
The Bloch remainders in \eqref{eq:dilute-cones} are stated for
each fixed medium; the uniform dynamical estimate follows
directly from Theorem~\ref{thm:effective-dynamics}.
\end{remark}

\begin{remark}[Comparison with the bubble framework]\label{rem:bubble-comparison}
The number of subwavelength bands depends on the material scaling as well as on the number of connected inclusions. In the bubble framework of \cite[\S4.2, Lemma~4.3]{AmmariDaviesHiltunen2026}, the density ratio, denoted here by $\delta_{\rm b}$, tends to zero with fixed geometry and wave speeds of order one. The zero Neumann mode in each connected bubble gives one leading-order amplitude, whereas the nonzero interior Neumann frequencies remain bounded away from zero. Consequently, one bubble per unit cell gives one subwavelength Bloch band. The corresponding modes are constant to leading order inside each bubble.

In \eqref{eq:dilute-weight}, the contrast enters only the mass coefficient $b$, and the interior wave speed is $\tau^{-1/2}$. The field and its normal derivative are continuous across the inclusion boundary, in the trace and weak flux senses, respectively. At the tuned contrast, define the interior wavenumber by $k_{\rm in}(\eta):=\omega_*(\eta)\sqrt{\tau_*(\eta)}$. Equation~\eqref{eq:dilute-tuned-asymptotics} gives
\begin{equation}\label{eq:dilute-interior-wavenumber}
 \omega_*(\eta)\longrightarrow0,\qquad
 \eta k_{\rm in}(\eta)\longrightarrow\sqrt{\Lambda_E(B)}>0
 \quad(\eta\to0).
\end{equation}
The background wavelength is therefore large compared with the unit period, while the interior wavelength is comparable to the inclusion size. This permits nonconstant internal modes, including the dipole doublet, at subwavelength frequencies. The distinction between the bubble and high-refractive-index mechanisms is also reflected in \cite[Appendix~C.3]{AmmariDaviesHiltunen2026}. In the present crystal, $D_4$ symmetry preserves the dipole degeneracy, and tuning $\tau$ brings the capacitary singlet to the same eigenvalue. The nonzero coupling proved in Theorem~\ref{thm:dilute-dirac} then produces the cone.
\end{remark}

\begin{remark}[Number of subwavelength bands]\label{rem:dilute-band-count}
At the tuned contrast, write $\lambda_n(\alpha;\eta)$ for the ordered Bloch eigenvalues, counted with multiplicity, and set $\omega_n(\alpha;\eta):=\sqrt{\lambda_n(\alpha;\eta)}$. Theorem~\ref{thm:dilute-dirac} gives
\begin{equation}\label{eq:dilute-first-four-bands}
\begin{gathered}
 \lambda_1(0;\eta)=0,\qquad
 \lambda_2(0;\eta)=\lambda_3(0;\eta)=\lambda_4(0;\eta)=\lambda_*(\eta),\\
 \lambda_5(0;\eta)\ge(1+c_B)\lambda_*(\eta).
\end{gathered}
\end{equation}
The first band has the constant mode at $\Gamma$, and the next three form the crossing. The relative gap at $\Gamma$ isolates this triple even though higher bands also become subwavelength.

To see the latter point, let $\mu_n^{\mathrm D}(B)$ be the $n$th Dirichlet eigenvalue of $-\Delta$ on $B$. For every $n\ge1$ and $\alpha\in Y^*$, the min--max principle gives
\begin{equation}\label{eq:dilute-band-count-bound}
 0\le\lambda_n(\alpha;\eta)
 \le\frac{\mu_n^{\mathrm D}(B)}{\tau_*(\eta)\eta^2}
 =\frac{2\pi\mu_n^{\mathrm D}(B)}
        {\Lambda_E(B)\ell_\eta}[1+o(1)]
 \quad(\eta\to0).
\end{equation}
Indeed, let $\varphi$ range over the span of the first $n$ Dirichlet eigenfunctions on $B$. Extend
$u_\alpha(x):=e^{-\ii\alpha\cdot x}\varphi(x/\eta)$ by zero outside $D_\eta$. Since $D_\eta\Subset Y$ and $\varphi\in H^1_0(B)$, this extension belongs to $H^1_\#(Y)$. For nonzero $\varphi$,
\[
 \frac{\int_Y|(\nabla+\ii\alpha)u_\alpha|^2\,dx}
      {\int_Y b_{\eta,\tau_*(\eta)}|u_\alpha|^2\,dx}
 =\frac{1}{\tau_*(\eta)\eta^2}
      \frac{\int_B|\nabla\varphi|^2\,dy}{\int_B|\varphi|^2\,dy}
 \le\frac{\mu_n^{\mathrm D}(B)}{\tau_*(\eta)\eta^2}.
\]
The resulting $n$-dimensional trial space proves \eqref{eq:dilute-band-count-bound}, uniformly in $\alpha$. It follows that
\begin{equation}\label{eq:dilute-all-fixed-bands}
 \sup_{\alpha\in Y^*}\omega_n(\alpha;\eta)\longrightarrow0
 \quad\text{for every fixed }n\ge1.
\end{equation}
Thus, the number of asymptotically subwavelength bands is not bounded by the number of inclusions per cell. For each fixed $\eta>0$ and $\alpha$, compactness of the fiber resolvent still gives only finitely many eigenfrequencies below any fixed finite cutoff. The limit in \eqref{eq:dilute-all-fixed-bands} is taken at fixed band index. The effective Dirac equation describes the three bands near $\omega_*(\eta)$, whose separation is quantified by \eqref{eq:dilute-uniform-relative-gap}.
\end{remark}

\subsection{Examples verifying the spectral-gap condition}
\label{subsec:dilute-shape-example}
We verify Assumption~\ref{assump:neumanngap} first by a variational bound for a square and then by separation of variables for a disk.

\subsubsection{Square inclusions}

The spectral-gap condition \eqref{eq:A1-Neumann-gap} can be verified without computing the transmission spectrum. Let $B=(-R,R)^2$. Separation of variables for the Neumann problem, restricted to invariant functions of mean zero, gives
\[
    \mu_A^{\mathrm{N}}(B)=\frac{\pi^2}{R^2},
\]
attained by $\cos(\pi y_1/R)+\cos(\pi y_2/R)$.
For \eqref{eq:LambdaE-dilute}, use the finite-energy trial function
\[
    U(y)=
    \begin{cases}
       y_1,& |y|\le R,\\
       R^2y_1/|y|^2,& |y|>R.
    \end{cases}
\]
It has the required parities and energy $2\pi R^2$.
Since $B$ contains the radius-$R$ disk,
\[
    \int_BU^2\,dx\ge\frac{\pi R^4}{4},
    \qquad
    \Lambda_E(B)\le\frac8{R^2}
       <\frac{\pi^2}{R^2}=\mu_A^{\mathrm{N}}(B).
\]
Thus, all preceding conclusions apply to a square inclusion;
one may choose $\rho=1/10$.

The strict inequality also persists under sufficiently small
$D_4$-equivariant bi-Lipschitz deformations of the square.
To see this, pull back the energy and mass integrals by a
$D_4$-equivariant map of $\mathbb R^2$, equal to the identity
outside a fixed ball, whose derivative and inverse derivative
are close to $I_2$. The Rayleigh quotients change by
multiplicative factors tending to one. For the Neumann quotient,
write the mass as $\inf_{c\in\mathbb R}\int_B|w-c|^2$ to account
for the changing mean constraint. This gives an open class of
noncircular reference inclusions satisfying
\eqref{eq:A1-Neumann-gap}.

\subsubsection{Disk inclusions}\label{sec:disk}
For a disk one can evaluate both variational constants and the normalized dipole moment. Let $B=B_R(0)$ with $R>0$. For $\nu\ge0$, denote by $J_\nu$ the Bessel function of the first kind of order $\nu$ and by $j_{\nu,1}$ its first positive zero. We claim that
\begin{equation}\label{eq:disk-constants}
 \Lambda_E(B)=\frac{j_{0,1}^2}{R^2},\qquad
 \mu_A^{\mathrm N}(B)=\frac{j_{1,1}^2}{R^2},\qquad
 d_B=\sqrt{\frac{2}{\pi}} .
\end{equation}
In polar coordinates $(r,\theta)$, seek a dipole with angular factor $\cos\theta$ and wave number $k>0$. With constants $c_{\rm in}>0$ and $c_{\rm out}>0$, take
\[
 U(r,\theta)=
 \begin{cases}
 c_{\rm in}J_1(kr)\cos\theta,&r\le R,\\
 c_{\rm out}\,r^{-1}\cos\theta,&r\ge R.
 \end{cases}
\]
The continuity of $U$ at $R$ gives $c_{\rm out}=c_{\rm in}R J_1(kR)$. The continuity of its normal derivative reduces to $kR J_1'(kR)+J_1(kR)=kR J_0(kR)=0$. Choose $k=j_{0,1}/R$. Interlacing of the Bessel zeros implies $J_1(kr)>0$ for $0<r\le R$. Thus, $U$ is a finite-energy solution of $-\Delta U=k^2\mathbf1_B U$ with the required parities and is positive on the right half-plane. Its restriction to the half-disk must belong to the principal eigenspace of the positive self-adjoint kernel in the proof of Lemma~\ref{lem:dilute-E-convergence}, since it has positive inner product with positive first eigenfunction. This proves $\Lambda_E(B)=k^2$. The normalization \eqref{U1property} gives
\[
 1=\pi(c_{\rm in})^2\int_0^R rJ_1(kr)^2\,dr
   =\frac{\pi(c_{\rm in})^2R^2}{2}J_1(j_{0,1})^2 .
\]
It follows that $c_{\rm out}=d_B=\sqrt{2/\pi}$. The recurrence, interlacing and integral identities used here are given in \cite[\S\S10.6, 10.21, 10.22]{DLMF}.

The invariant Neumann modes have angular factors $\cos(4m\theta)$, $m=0,1,\ldots$. In the radial sector the first positive eigenvalue satisfies $J_0'(kR)=-J_1(kR)=0$. For $m\ge1$ the angular part of the Dirichlet energy alone bounds the Rayleigh quotient below by $16/R^2$. Since $j_{1,1}\simeq3.831706<4$, the first positive invariant Neumann eigenvalue is radial. This completes \eqref{eq:disk-constants}. Since $j_{0,1}<j_{1,1}$, the disk satisfies Assumption~\ref{assump:neumanngap}. The tuned contrast and the cone speed are, therefore,
\begin{equation}\label{eq:disk-speed}
 \tau_*(\eta)\sim\frac{j_{0,1}^2}{2\pi R^2}\frac{\ell_\eta}{\eta^2},\qquad
 \frac{\gamma_\eta}{2\omega_*(\eta)}
 \sim\frac{\sqrt{2\pi}}{j_{0,1}}R\eta
 \quad (\eta\to0).
\end{equation}
The dimensionless speed coefficient $\sqrt{2\pi}/j_{0,1}$ is approximately $1.04233$.

\section{Concluding remarks}\label{sec:conclusion}
A square-symmetric medium can have an isotropic conical crossing of three bands at $\Gamma$ when a singlet and a doublet are tuned to the same eigenvalue and their linear coupling is nonzero. The two constructions yield distinct velocity scales. For weak perturbations of a medium with stiffness $aI_2$ and unit density, either singlet gives a cone speed tending to $\sqrt{a/2}$. For a dilute inclusion, the speed tends to zero linearly in its size. Its coefficient depends on the normalized dipole moment of the reference inclusion, whereas the leading squared crossing frequency, $2\pi/\log(1/\eta)$, is independent of the shape. In the dilute family, the tuning is transverse because nearly all of the doublet's weighted mass lies inside the inclusion, while the singlet's weighted mass there tends to zero. Although every fixed Bloch band becomes subwavelength in this limit, the relative spectral gap at $\Gamma$ isolates the three modes used in the effective equation.

The effective equation separates propagation along the two conical branches from the stationary transverse component. For a fixed medium, it approximates the positive-frequency evolution in the weighted $L^2$ norm on times $t=o(\delta^{-2})$, where $\delta$ is the quasimomentum bandwidth. For a propagating component in the dilute family, an order-one displacement in the slow variable occurs on the time scale $t\asymp(\delta\eta)^{-1}$. Corollary~\ref{cor:dilute-dynamics} justifies the approximation on this scale if $\delta\sqrt{\log(1/\eta)}=o(\eta)$. For the weak-contrast families, Proposition~\ref{prop:weak-axial-curvature} gives a positive quadratic coefficient of the central band along both coordinate axes. Further questions include the evaluation of the directional coefficient in \eqref{eq:central-directional-quadratic} for the dilute crystal and the justification of a second-order wave-packet approximation on the scale $t\asymp\delta^{-2}$.

For the dielectric interpretation, the cone speed determines the product of the two constitutive slopes. A uniform cell source couples to the singlet through its unweighted mean. When this mean is nonzero, the forced response determines each slope separately. Corollary~\ref{cor:dilute-double-zero} shows that the limiting bulk impedance of the dilute crystal grows as $\eta^{-1}$. Determining reflection and transmission for a finite sample requires an analysis of the fields at its interface \cite{markel}.

\section*{Statements and declarations}

\paragraph{Funding.}

This work was partially supported by the City University of Hong Kong start-up fund 7200843, Hong Kong RGC Early Career Scheme grant 21301326.

\paragraph{Competing interests.}
The authors have no competing interests to declare.

\paragraph{Data availability.}
No external datasets were used. 

\subsection*{Declaration of AI use}
The authors interacted with and guided ChatGPT~6 whilst developing the proofs and then wrote the paper. The authors take full responsibility for the accuracy and integrity of all content in this paper.

\appendix
\section{Real representations of finite groups}\label{app:rep-proofs}
We give the general decomposition and projection results used for the square group in Section~\ref{subsec:representation-preliminaries}. All spaces and representations in this appendix are real.

\subsection{Representations and intertwiners}

Let $G$ be a finite group with identity $e$, and $V$ be either a finite-dimensional real vector space, or a real Hilbert space. Write $\operatorname{GL}(V)$ for the group of bounded invertible linear maps on $V$. A real representation of $G$ on $V$ is a map $\rho:G\to\operatorname{GL}(V)$ satisfying
\begin{equation}
    \rho(gh)=\rho(g)\rho(h),\qquad \rho(e)= I_V.
\end{equation}
If $\rho(g)=I_V$ for every $g\in G$, $\rho$ is said to be trivial. If $V$ is finite-dimensional, the dimension of $\rho$ is $\dim V$. A one-dimensional real representation is called a character.

A subspace $W\subset V$ is invariant if $W$ is closed and $\rho(g)W\subset W$ for
every $g$. A representation on $V$ is irreducible if $V\ne\{0\}$ and its only invariant subspaces are $\{0\}$ and $V$. Every irreducible representation of the finite group $G$ is finite-dimensional. Indeed, for $v\ne0$, the orbit span $\operatorname{span}\{\rho(g)v:g\in G\}$ is a nonzero invariant subspace of dimension at most $|G|$, hence is equal to the whole space. 

An intertwiner from $\rho$ on $V$ to $\sigma$ on $W$ is a bounded linear map $A:V\to W$ satisfying
\begin{equation}
    A\rho(g)=\sigma(g)A, \qquad g\in G.
\end{equation}
The two representations are equivalent, written $\rho\simeq\sigma$, if they admit an invertible intertwiner. With this equivalence relation, we write $\widehat G$ for a set containing one representative of each equivalence class of irreducible real representations of $G$.

Direct sums and tensor products are defined by
\begin{equation}
    (\rho\oplus\sigma)(g)=\rho(g)\oplus\sigma(g), \qquad (\rho\otimes\sigma)(g)=\rho(g)\otimes\sigma(g).
\end{equation}

If $V=H$ is a real Hilbert space, a representation is orthogonal if it preserves the given inner
product, equivalently $\rho(g)^*=\rho(g^{-1})$ for every $g\in G$. We write $\mathrm O(H)$ for the group of surjective linear isometries of $H$. Thus, an orthogonal representation on $H$ is precisely a homomorphism from $G$ to $\mathrm O(H)$.

\begin{lemma}[Schur's lemma]\label{rep-prelim:schur}
A nonzero intertwiner between irreducible real representations is invertible. In particular, if $\rho$ is irreducible and orthogonal, every self-adjoint endomorphism $B$ commuting with all $\rho(g)$ is a real scalar multiple of the identity.
\end{lemma}

\begin{proof}
Let $A:V\to W$ be an intertwiner from $\rho$ on $V$ to $\sigma$ on $W$. Both $\ker A$ and $\operatorname{Ran}A$ are invariant. Therefore, irreducibility gives $\ker A=\{0\}$ and $\operatorname{Ran}A=W$ when $A\ne0$. Hence, $A$ is invertible. 

For the second assertion, self-adjointness and finite dimension give a real eigenvalue $\lambda$ of $B$. The intertwiner $B-\lambda I$ has a nonzero kernel, so the first assertion implies $B-\lambda I=0$.
\end{proof}

\subsection{Irreducible decompositions}

Let $T:G\to\mathrm O(H)$ be an orthogonal representation on a separable real Hilbert space $H$. Write $T_g=T(g)$. The following lemmas describe the orthogonal decomposition into irreducible subspaces. Their proofs are given in Section~\ref{app:rep-arguments}.

\begin{lemma}\label{lem:orthogonal-irreducible-decomposition}
There exists an at most countable family $\{U_j\}_{j\in J}$ of mutually orthogonal, nonzero, finite-dimensional irreducible invariant subspaces such that
\begin{equation}
    H=\overline{\operatorname{span} \Bigl(\bigcup_{j\in J}U_j\Bigr)}.
\end{equation}
\end{lemma}

\begin{lemma}\label{rep-prelim:complete-reducibility} 
For each $\rho\in\widehat G$, there exist an at most countable index set $A_\rho$ and nonzero finite-dimensional invariant subspaces $H_{\rho,a}\subset H$, $a\in A_\rho$, such that:
\begin{enumerate}
    \item[(a)] the representation $T$ on each $H_{\rho,a}$ is equivalent to $\rho$;
    \item[(b)] the subspaces $H_{\rho,a}$, over all pairs $(\rho,a)$, are mutually orthogonal;
    \item[(c)] every $h\in H$ has a unique expansion
    \begin{equation}
        h=\sum_{\rho\in\widehat G}\sum_{a\in A_\rho}h_{\rho,a},\qquad h_{\rho,a}\in H_{\rho,a},
    \end{equation}
    where the sum converges in the $\|\cdot\|_H$ norm. 
\end{enumerate}
These properties are expressed by the notation
\begin{equation}
    H=\widehat{\bigoplus}_{\rho\in\widehat G} \widehat{\bigoplus}_{a\in A_\rho}H_{\rho,a}.
\end{equation}
\end{lemma}

Fix a decomposition
\begin{equation}\label{decompofHeq}
    H=\widehat{\bigoplus}_{\rho\in\widehat G}\widehat{\bigoplus}_{a\in A_\rho}H_{\rho,a}
\end{equation}
as in Lemma~\ref{rep-prelim:complete-reducibility}, and let $\Pi_{\rho,a}$ be the orthogonal projection onto $H_{\rho,a}$. In the strong operator topology,
\begin{equation}\label{eq:irreducible-resolution-identity}
    I_H=\sum_{\rho\in\widehat G}\sum_{a\in A_\rho}\Pi_{\rho,a}.
\end{equation}

\begin{lemma}\label{lem:matrix-coefficient-projections}
Let $\rho\in\widehat G$ be an orthogonal irreducible representation on $V_\rho$, with
$d_\rho:=\dim V_\rho$. We assume that $\rho$ is absolutely irreducible, \textit{i.e.},
\begin{equation}\label{eq:scalar-commutant-hypothesis}
    \{L\in\operatorname{End}(V_\rho):L\rho(g)=\rho(g)L\text{ for all }g\in G\}=\mathbb R I_{V_\rho}.
\end{equation}
Choose an orthonormal basis $e_1,\ldots,e_{d_\rho}$ of $V_\rho$. Define
\begin{equation}\label{eq:matrix-coefficient-projections}
    P_{\rho,i}:=\frac{d_\rho}{|G|}\sum_{g\in G} \langle \rho(g) e_i ,e_i \rangle T_g, \qquad 1\le i\le d_\rho.
\end{equation}
For each $a\in A_\rho$ in \eqref{decompofHeq}, choose an orthogonal intertwiner $U_{\rho,a}:V_\rho\to H_{\rho,a}$, and let $\pi_{\rho,a,i}$ denote the orthogonal projection of $H$ onto the line $\mathbb R U_{\rho,a}e_i$.

Then $P_{\rho,i}$, $1\le i\le d_\rho$, are mutually orthogonal projections. More precisely,
\begin{equation}\label{eq:coordinate-copy-projections}
    P_{\rho,i}=\sum_{a\in A_\rho}\pi_{\rho,a,i}, \qquad
    \sum_{i=1}^{d_\rho}P_{\rho,i} =\sum_{a\in A_\rho}\Pi_{\rho,a}.
\end{equation}
If $\sigma\not\simeq\rho$ is another type satisfying \eqref{eq:scalar-commutant-hypothesis}, then $P_{\rho,i}$ and $P_{\sigma,j}$ have orthogonal ranges for every $i,j$.
\end{lemma}

\subsection{Proofs of the decomposition and projection lemmas}\label{app:rep-arguments}
We first construct an orthogonal decomposition into irreducible copies and then identify the coordinate projections on each copy.

\begin{proof}[Proof of Lemma~\ref{lem:orthogonal-irreducible-decomposition}]
The proof consists of two steps.

\textit{Step 1.} We show that every nonzero invariant subspace of $H$ contains a nonzero finite-dimensional irreducible invariant subspace of dimension at most $|G|$.

Let $W\subset H$ be a nonzero invariant subspace, and choose $0\ne h\in W$. The orbit span
\begin{equation}
    V_h:=\operatorname{span}\{T_gh:g\in G\}
\end{equation}
is a nonzero invariant subspace of $W$ with $\dim V_h\le |G|$. Among the nonzero invariant subspaces of $V_h$, choose one of smallest dimension, denoted by $U$. It has no proper nonzero invariant subspace and is therefore irreducible. 

\textit{Step 2.} Consider families of mutually orthogonal, nonzero, finite-dimensional irreducible invariant subspaces, ordered by inclusion. The union of any chain is again such a family. Zorn's lemma therefore supplies a maximal family $\{U_j\}_{j\in J}$. Set
\begin{equation}
    M:=\overline{\operatorname{span}\Bigl(\bigcup_{j\in J}U_j\Bigr)}.
\end{equation}
Since every $U_j$ is invariant and every $T_g$ is bounded, $M$ is invariant. Since $T_g^*=T_{g^{-1}}$, the orthogonal complement $M^\perp$ is also invariant: for $v\in M^\perp$ and $m\in M$,
\begin{equation}
    \langle T_gv,m \rangle_H= \langle v,T_{g^{-1}}m \rangle_H=0.
\end{equation}
If $M\ne H$, then $M^\perp$ is a nonzero invariant subspace. By the result obtained from Step 1, $M^\perp$ contains a nonzero finite-dimensional irreducible invariant subspace $U$. This $U$ is orthogonal to every $U_j$ and could be added to the family, contradicting maximality. Hence $M=H$.

Finally, choosing a unit vector in each $U_j$ gives an orthonormal set in the separable Hilbert space $H$. Such a set is at most countable, so $J$ is at most countable.
\end{proof}

\begin{proof}[Proof of Lemma~\ref{rep-prelim:complete-reducibility}]
Let $\{U_j\}_{j\in J}$ be the family obtained in  Lemma~\ref{lem:orthogonal-irreducible-decomposition}.
For each $\rho\in\widehat G$, set
\begin{equation}
    A_\rho:=\{j\in J:T|_{U_j}\simeq\rho\}, \qquad H_{\rho,a}:=U_a, \quad a\in A_\rho.
\end{equation}
Each irreducible representation is equivalent to exactly one representative in $\widehat G$, so the sets $A_\rho$ partition $J$. This proves parts (a) and (b). Since $J$ is at most countable, and $H$ is the closed span of the union of $U_j$, part (c) follows from orthogonality and the Pythagorean identity.
\end{proof}

\begin{proof}[Proof of Lemma~\ref{lem:matrix-coefficient-projections}]
The proof is divided into several steps.

\textit{Step 1: Existence of orthogonal intertwiners.}
Since the representation on $H_{\rho,a}$ is equivalent to $\rho$,
there is an invertible intertwiner $J:V_\rho\to H_{\rho,a}$.
Orthogonality of $T_g$ gives
\begin{equation}
    \rho(g)^*J^*J\rho(g)=(J\rho(g))^*(J\rho(g))=(T_gJ)^*(T_gJ)=J^*J.
\end{equation}
Hence $J^*J$ commutes with every $\rho(g)$.
By \eqref{eq:scalar-commutant-hypothesis},
$J^*J=cI_{V_\rho}$ for some $c>0$.
Thus $U_{\rho,a}:=c^{-1/2}J$ is an orthogonal intertwiner,
as required.

\textit{Step 2: The averaging identity.}
Let $\sigma:G\to\mathrm O(W)$ be either $\rho$ itself or an irreducible orthogonal representation which is not equivalent to $\rho$. For a linear map $A:V_\rho\to W$, define
\begin{equation}
    \mathcal M_{\sigma,\rho}(A):=\frac1{|G|}\sum_{g\in G}\sigma(g)A\rho(g)^{-1}.
\end{equation}
For every $h\in G$, reindexing the sum by $k=hg$ gives
\begin{equation}
    \sigma(h)\mathcal M_{\sigma,\rho}(A)\rho(h)^{-1}=\frac1{|G|}\sum_{g\in G}\sigma(hg)A\rho(hg)^{-1}
    =\mathcal M_{\sigma,\rho}(A).
\end{equation}
Therefore $\mathcal M_{\sigma,\rho}(A)$ is an intertwiner from $\rho$ to $\sigma$.

If $\sigma\not\simeq\rho$, Schur's lemma (Lemma \ref{rep-prelim:schur}) gives $\mathcal M_{\sigma,\rho}(A)=0$. If $\sigma=\rho$ on $V_\rho$, condition \eqref{eq:scalar-commutant-hypothesis} instead implies that $\mathcal M_{\rho,\rho}(A)=cI_{V_\rho}$. Taking traces, and using invariance of trace under conjugation, yields
\begin{equation}
    cd_\rho =\operatorname{Tr}\mathcal M_{\rho,\rho}(A) =\operatorname{Tr}A.
\end{equation}
Consequently,
\begin{equation}\label{eq:matrix-projection-averaging}
    \mathcal M_{\sigma,\rho}(A)=
       \begin{cases}
          0,&\sigma\not\simeq\rho,\\[2pt]
          \dfrac{\operatorname{Tr}A}{d_\rho}I_{V_\rho}, &\sigma=\rho.
       \end{cases}
\end{equation}

Choose an orthonormal basis $f_1,\ldots,f_{\dim W}$ of $W$, taking $f_k=e_k$ when $\sigma=\rho$. For fixed $i,\ell$, apply \eqref{eq:matrix-projection-averaging} to the rank-one map
\begin{equation}
    A v:=\langle v,e_i\rangle f_\ell.
\end{equation}
Since $\rho(g)^{-1}=\rho(g)^\top$, we have
\begin{equation}
    \begin{aligned}
    \left\langle\mathcal M_{\sigma,\rho}(A)e_i,f_k\right\rangle
       &=\frac1{|G|}\sum_{g\in G}\langle\rho(g)^{-1}e_i,e_i\rangle\langle\sigma(g)f_\ell,f_k\rangle\\
       &=\frac1{|G|}\sum_{g\in G} \langle\rho(g)e_i,e_i\rangle\langle\sigma(g)f_\ell,f_k\rangle .
\end{aligned}
\end{equation}
When $\sigma=\rho$, the rank-one map $A$ has trace $\delta_{i\ell}$. Thus \eqref{eq:matrix-projection-averaging}
gives
\begin{equation}\label{eq:diagonal-coefficient-orthogonality}
    \frac{ d_\rho}{|G|}\sum_{g\in G} \langle\rho(g)e_i,e_i\rangle\langle\sigma(g)f_\ell,f_k\rangle =
       \begin{cases}
          0,&\sigma\not\simeq\rho,\\
          \delta_{ik}\delta_{i\ell},&\sigma=\rho.
       \end{cases}
\end{equation}

\textit{Step 3: Identification of the projections.} Since $U_{\rho,a}$ is an intertwiner, $U_{\rho,a} \rho(g)=T_g U_{\rho,a}$. Hence,
\begin{equation}
    U_{\rho,a}^{-1}P_{\rho,i}U_{\rho,a}=\frac{ d_\rho}{|G|}\sum_{g\in G} \langle\rho(g)e_i,e_i\rangle \rho(g).
\end{equation}
By \eqref{eq:diagonal-coefficient-orthogonality}, the $(k,\ell)$ entry of this matrix is
$\delta_{ik}\delta_{i\ell}$. It is therefore the coordinate projection $v\mapsto\langle v,e_i\rangle e_i$.
Because $U_{\rho,a}$ is orthogonal, it follows that
\begin{equation}
    P_{\rho,i}|_{H_{\rho,a}}=\pi_{\rho,a,i}|_{H_{\rho,a}}.
\end{equation}
On an inequivalent copy $H_{\sigma,a}$, choose an orthonormal basis of that copy and apply the first case of
\eqref{eq:diagonal-coefficient-orthogonality} to its representation matrices. All matrix entries vanish, so
\begin{equation}
    P_{\rho,i}|_{H_{\sigma,a}}=0, \qquad \sigma\not\simeq\rho.
\end{equation}

For any $h\in H$, the fixed orthogonal decomposition \eqref{eq:irreducible-resolution-identity} yields
\begin{equation}
    h=\sum_{\sigma\in\widehat G}\sum_{a\in A_\sigma}\Pi_{\sigma,a}h \qquad\text{in }H.
\end{equation}
Applying the operator $P_{\rho,i}$ and using the restrictions just computed, we obtain
\begin{equation}\label{Pequapisuma}
    P_{\rho,i}h=\sum_{a\in A_\rho}\pi_{\rho,a,i}h.
\end{equation}
This proves the first identity in \eqref{eq:coordinate-copy-projections}. In particular, $P_{\rho,i}$ are mutually orthogonal projections.

For fixed $a$, the vectors $U_{\rho,a}e_1,\ldots,U_{\rho,a}e_{d_\rho}$ form an orthonormal basis of $H_{\rho,a}$. Hence
\begin{equation}\label{sumipiequalPi}
    \sum_{i=1}^{d_\rho}\pi_{\rho,a,i}=\Pi_{\rho,a}.
\end{equation}
Summing \eqref{Pequapisuma} over $1\leq i\leq d_\rho$ and using \eqref{sumipiequalPi}, we obtain
the second identity in \eqref{eq:coordinate-copy-projections}.

Finally, projections associated with inequivalent types have ranges contained in orthogonal sums of distinct irreducible copies. Their ranges are therefore orthogonal.
\end{proof}

\bibliographystyle{abbrv}
\bibliography{references}

\end{document}